%% file: paper.tex
\documentclass[a4paper,3p,10pt,times]{elsarticle}

\usepackage[utf8]{inputenc}
\usepackage{amsmath}
\usepackage{amsfonts}
\usepackage{amssymb}
\usepackage{amsthm}
\usepackage[dvipsnames,table]{xcolor}
\usepackage{graphicx}

\usepackage{tikz,tikz-3dplot}
\usetikzlibrary{spy}
\usetikzlibrary{matrix}
\usetikzlibrary{arrows}
\usetikzlibrary{pgfplots.groupplots}

\usepackage{booktabs}
\usepackage{multirow}
\usepackage{multicol}

\usepackage{pgfplots}
\usepgfplotslibrary{fillbetween}
\usepackage{subcaption}

\usepackage{bm}
\usepackage{dutchcal}
\usepackage{mathtools}
\usepackage{enumerate}
\usepackage{hyperref}

\usepackage{comment}

\hypersetup{
    colorlinks = true,
}
\usepackage[capitalize]{cleveref}

\usepackage{physics}

\newcommand{\VEC}[1]{\vb*{#1}}
\newcommand{\UVEC}[1]{\hat{\vb*{#1}}}
\newcommand{\MAT}[1]{#1}
\newcommand{\TEN}[1]{\vb{#1}}

\newcommand{\xvec}{\vb*{x}}
\newcommand{\Xvec}{\mathring{\vb*{x}}}

\newcommand{\avec}{\vb*{a}}

\newcommand{\Avec}{\mathring{\vb*{a}}}

\newcommand{\uvec}{\vb*{u}}

\newcommand{\epsten}{\bm{\varepsilon}}

\newcommand{\Ften}{\vb{F}}
\newcommand{\Eten}{\vb{E}}

\newcommand{\Nten}{\vb{N}}

\newcommand{\Sten}{\vb{S}}
\newcommand{\Cten}{\vb{C}}
\newcommand{\Ccten}{\mathbcal{C}}
\newcommand{\Ccoef}{\mathcal{C}}

\usepackage{lineno}

\biboptions{sort&compress}

\input{TUDColors}
\input{Triangle.tex}

\pgfplotsset{
	cycle list name = MyCyclelist,
	width = 0.9\linewidth,
    grid=major,
  	legend cell align={left},
    ylabel near ticks,
    xlabel near ticks,
    legend style={fill=white, fill opacity=0.6, draw opacity=1,draw=none,text opacity=1},
}

\theoremstyle{plain}
\newtheorem{remark}{Remark}
\crefname{rem}{Remark}{Remarks}
\newtheorem{example}{Example}
\crefname{ex}{Example}{Examples}

\newcommand{\fullcite}[1]{\cite{#1}}

\newcommand{\fullcitelist}[2]{\cite{#2}}

\usepackage[compact]{titlesec}

\tikzset{%
style0/.style = {thick,densely dashed,black,mark=none,mark size=2,every mark/.append style={solid,line width = 0.5}},
style1/.style = {thick,solid,col1,mark=*,mark size=2,every mark/.append style={solid,line width = 0.5}},
style2/.style = {thick,densely dashed,col2,mark=square*,mark size=2,every mark/.append style={solid,line width = 0.5}},
style3/.style = {thick,densely dotted,col3,mark=triangle*,mark size=2,every mark/.append style={solid,line width = 0.5}},
style4/.style = {thick,solid,col4,mark=*,mark size=2,every mark/.append style={solid,line width = 0.5}},
}

\begin{document}

\begin{frontmatter}

\title{A Wrinkling Model for General Hyperelastic Materials based on Tension Field Theory}

\author[label1]{\corref{cor1}H.M. Verhelst}
\ead{hugomaarten.verhelst@unipv.it}
\cortext[cor1]{Corresponding Author. }
\author[label3]{M. Möller}
\ead{m.moller@tudelft.nl}
\author[label2]{J.H. Den Besten}
\ead{henk.denbesten@tudelft.nl}
\address[label1]{University of Pavia, Department of Civil Engineering and Architecture, Via Adolfo Ferrata 3, 27100 PV, Pavia, Italy}
\address[label2]{Delft University of Technology, Department of Maritime and Transport Technology, Mekelweg 2, 2628 CD, Delft, The Netherlands}
\address[label3]{Delft University of Technology, Delft Institute of Applied Mathematics, Mekelweg 4, 2628 CD, Delft, The Netherlands}

\begin{abstract}
Wrinkling is the phenomenon of out-of-plane deformation patterns in thin walled structures, as a result of a local compressive (internal) loads in combination with a large membrane stiffness and a small but non-zero bending stiffness. Numerical modelling typically involves thin shell formulations. As the mesh resolution depends on the wrinkle wave lengths, the analysis can become computationally expensive for shorter ones. Implicitly modeling the wrinkles using a modified kinematic or constitutive relationship based on a taut, slack or wrinkled state derived from a so-called tension field, a simplification is introduced in order to reduce computational efforts. However, this model was restricted to linear elastic material models in previous works. Aiming to develop an implicit isogeometric wrinkling model for large strain and hyperelastic material applications, a modified deformation gradient has been assumed, which can be used for any strain energy density formulation. The model is an extension of a previously published model for linear elastic material behaviour and is generalized to other types of discretisation as well. The extension for hyperelastic materials requires the derivative of the material tensor, which can be computed numerically or derived analytically. The presented model relies on a combination of dynamic relaxation and a Newton-Raphson solver, because of divergence in early Newton-Raphson iterations as a result of a changing tension field, which is not included in the stress tensor variation. Using four benchmarks, the model performance is evaluated. Convergence with the expected order for Newton-Raphson iterations has been observed, provided a fixed tension field. The model accurately approximates the mean surface of a wrinkled membrane with a reduced number of degrees of freedom in comparison to a shell solution.
\end{abstract}

\begin{keyword}
isogeometric analysis \sep Kirchhoff-Love shell \sep Tension field theory \sep Dynamic Relaxation method \sep Hyperelasticity
\end{keyword}

\end{frontmatter}


\section{Introduction}
\input{Sections_Introduction.tex}

\section{Isogeometric membrane formulation}\label{sec:membrane_element}
\input{Sections_Membrane.tex}

\section{Linear elastic wrinkling model}\label{sec:elastic}
\input{Sections_Linear.tex}

\section{Hyperelastic wrinkling model}\label{sec:hyperelastic}
\input{Sections_Hyperelastic.tex}

\section{Numerical Solution Strategies}\label{sec:solution}
\input{Sections_SolutionStrategies.tex}

\section{Benchmarks}\label{sec:results}
\input{Sections_Benchmarks.tex}

\section{Conclusions}\label{sec:conclusions}
\input{Sections_Conclusions.tex}

\section*{Acknowledgments}
This research has been conducted as a part of the project FlexFloat with project number 19002 of the research programme Open Technology which is (partly) financed by the Dutch Research Council (NWO). In addition, HMV is grateful for the financial support from the Italian Ministry of University and Research (MUR) through the PRIN projects COSMIC (no. 2022A79M75) and ASTICE (no. 202292JW3F), with the contribution of the European union -- Next Generation EU.

\bibliography{references}

\end{document}

%% file: TUDColors.tex
\definecolor{TUcol1}{cmyk}{1,0,0,0}
\definecolor{TUcol2}{cmyk}{0.54,0,0.32,0}
\definecolor{TUcol3}{cmyk}{1,0.15,0.4,0}
\definecolor{TUcol4}{cmyk}{1,0.66,0,0.4}
\definecolor{TUcol5}{cmyk}{0.98,1,0,0.35}
\definecolor{TUcol6}{cmyk}{0.82,0,0.21,0.08}
\definecolor{TUcol7}{cmyk}{0.45,0,0.06,0.06}
\definecolor{TUcol8}{cmyk}{0.45,0.2,0,0.07}
\definecolor{TUcol9}{cmyk}{0.02,0.56,0.84,0}
\definecolor{TUcol10}{cmyk}{0.58,1,0,0.02}
\definecolor{TUcol11}{cmyk}{0.19,1,0,0.19}
\definecolor{TUcol12}{cmyk}{0.36,0,1,0}
\definecolor{TUcol13}{cmyk}{0.02,0,0.54,0}

\definecolor{col1}{HTML}{0480B0}
\definecolor{col2}{HTML}{FFB600}
\definecolor{col3}{HTML}{FF2C00}
\definecolor{col4}{HTML}{AA6F39}

\definecolor{chop1}{RGB}{92,76,161}
\definecolor{chop2}{RGB}{87,172,188}
\definecolor{chop3}{RGB}{161,217,163}
\definecolor{chop4}{RGB}{235,245,156}
\definecolor{chop5}{RGB}{250,230,154}
\definecolor{chop6}{RGB}{253,163,93}
\definecolor{chop7}{RGB}{227,81,74}
\definecolor{chop8}{RGB}{156,0,63}

\pgfplotscreateplotcyclelist{MyCyclelist}{%
thick,solid,col1,mark=*,mark size=1,every mark/.append style={solid,line width = 0.5}\\
thick,densely dashed,col2,mark=square*,mark size=1,every mark/.append style={solid,line width = 0.5}\\
thick,densely dotted,col3,mark=triangle*,mark size=1,every mark/.append style={solid,line width = 0.5}\\
thick,dotted,col4,mark=o*,mark size=1,every mark/.append style={solid,line width = 0.5}\\
}

\pgfplotsset{
    layers/legend behind plots/.define layer set={
            axis background,axis grid,axis ticks,axis lines,axis tick labels,main,axis descriptions,axis foreground
    }{
        grid style= {/pgfplots/on layer=axis grid},
        tick style= {/pgfplots/on layer=axis ticks},
        axis line style= {/pgfplots/on layer=axis lines},
        label style= {/pgfplots/on layer=axis descriptions},
        legend style= {/pgfplots/on layer=axis tick labels}, 
        title style= {/pgfplots/on layer=axis descriptions},
        colorbar style= {/pgfplots/on layer=axis descriptions},
        ticklabel style= {/pgfplots/on layer=axis tick labels},
        axis background@ style={/pgfplots/on layer=axis background},
        3d box foreground style={/pgfplots/on layer=axis foreground},
    },
    MyStyle1/.style={
        MyStyleA/.style={red},
    },
    MyStyle2/.style={
        MyStyleA/.style={blue},
    },
}

%% file: Triangle.tex
\newcommand{\logLogSlopeTriangleRev}[5]
{
\small

    \pgfplotsextra
    {
        \pgfkeysgetvalue{/pgfplots/xmin}{\xmin}
        \pgfkeysgetvalue{/pgfplots/xmax}{\xmax}
        \pgfkeysgetvalue{/pgfplots/ymin}{\ymin}
        \pgfkeysgetvalue{/pgfplots/ymax}{\ymax}

        \pgfmathsetmacro{\xArel}{#1-#2}
        \pgfmathsetmacro{\yArel}{#3}
        \pgfmathsetmacro{\xBrel}{#1}
        \pgfmathsetmacro{\yBrel}{\yArel}
        \pgfmathsetmacro{\xCrel}{\xArel}

        \pgfmathsetmacro{\lnxB}{\xmin*(1-(#1-#2))+\xmax*(#1-#2)} 
        \pgfmathsetmacro{\lnxA}{\xmin*(1-#1)+\xmax*#1} 
        \pgfmathsetmacro{\lnyA}{\ymin*(1-#3)+\ymax*#3} 
        \pgfmathsetmacro{\lnyC}{\lnyA-#4*(\lnxA-\lnxB)}
        \pgfmathsetmacro{\yCrel}{\lnyC-\ymin)/(\ymax-\ymin)} 

        \coordinate (A) at (rel axis cs:\xArel,\yArel);
        \coordinate (B) at (rel axis cs:\xBrel,\yBrel);
        \coordinate (C) at (rel axis cs:\xCrel,\yCrel);

        \draw[#5]   (A)-- node[pos=0.5,anchor=south] {1}
                    (B)--
                    (C)-- node[pos=0.5,anchor=east] {#4}
                    cycle;
    }
}

%% file: Sections_Introduction.tex
Wrinkling is a phenomenon that is omnipresent around us and appears at different scales: it influences the thermo-conductivity of graphene on the nanoscale \cite{Deng2016} and the reflectivity of solar sails on the metre scale \cite{Sakamoto2007}. Wrinkling shapes floating leaves \cite{Xu2020d} and plays an important role in cosmetics, for example in wound healing processes \cite{Swain2015}. In engineering, the wrinkling simulation is relevant in the design of airbags \cite{Graczykowski2016,Meduri2022}, sails \cite{Renzsch2018a}, parachutes \cite{Takizawa2016}, and floating solar platforms \cite{Verhelst2019}. Membrane wrinkling is influenced by balancing potential energy stored in bending deformations, membrane deformations, or in fluid or solid foundations, as discussed in the seminal works of \fullcite{Cerda2002} and \fullcite{Pocivavsek2008}. For further reference on the study of wrinkling from a physics perspective, the reader is referred to the review papers by \fullcite{Wang2022} on tension-induced wrinkles, by \fullcite{Tan2020} on wrinkles in curved surfaces, by \fullcite{Ma2020} on wrinkles in membranes with elasticity gradients, by \fullcite{Paulsen2019} on wrinkles in membranes with low bending stiffness and high membrane stiffness, and by \fullcite{Li2012} for a complete but less recent review on wrinkling.\\

For engineering applications where wrinkling is an relevant factor, accurate and efficient numerical modelling of wrinkling patterns becomes of great importance. In general, the numerical modelling can be done in different ways: by explicitly modelling the wrinkling amplitudes and wave lengths using shell models, by using reduced-order models based on the Föppl-Von Kármán plate equations, or by implicitly embedding the effects of wrinkling in element formulations.\\

Firstly, the modelling of membrane wrinkling can be done using mathematical models that account for the physics of thin films, including both membrane and bending effects. Plate and shell models are particularly suited to simulate wrinkling patterns, whereas membrane models lacking bending stiffness cannot capture wrinkle formations. Many studies on the physics of wrinkling employed commercial finite element methods to numerically investigate wrinkling under different conditions
\cite{Wong2006c,Nayyar2011,Vella2011,Panaitescu2019,Nazzal2022,Yan2014b,Li2017,Stoop2015,Yang2018a,Taffetani2017,Luo2020,Faghfouri2020,Matoz-Fernandez2020}. In addition, dedicated numerical models for the modelling of wrinkling patterns have been developed. For initially flat geometries, the Föppl-Von Kármán (FvK) model incorporating out-of-plane displacements and linear bending strains has been applied \cite{Kim2012,Healey2013,Sipos2016,Li2016,Wang2022,Fu2019,Liu2019,Yang2020}. This FvK model has been extended for hyperelastic materials, orthotropy, and general anisotropic properties \cite{Fu2021,Fu2022,Wang2020,Wang2021,Wang2022c,Wang2022g}. Another model based on Koiter's non-linear plate theory, proposed by \fullcite{Steigmann2013}, has been used to compute wrinkling cases involving shear, holes, annuli, graphene, and reinforced plates \cite{Taylor2014,Wong2006a,Taylor2015,Qin2014,Taylor2019,Taylor2020}. In order to find wrinkling patterns, different algorithms searching for equilibrium solutions have been used, including static methods like the dynamic relaxation method and Newton--Raphson, as demonstrated by Taylor \cite{Taylor2014}, or continuation methods such as the Arc-Length Method \cite{Riks1972,Crisfield1983} and the Asymptotic Numerical Method (ANM) \cite{Cochelin1994,Cochelin2007,Cochelin2013,Damil1990,Najah1998}. In general, the advantage of explicitly modelling wrinkling amplitudes and wave lengths is that they resemble actual physics. However, when wrinkling wave lengths decrease, the mesh size required to find the solutions typically decreases as well, implying increased computational costs.\\

Secondly, the so-called \emph{Fourier reduced model} is a technique where the Föppl-Von Kármán model is discretised using Fourier series expansions. Consequently, the model provides a multi-scale model where the large Fourier modes capture the macroscopic deformations and the high-frequency content captures the wrinkling patterns. This technique, introduced by \fullcitelist{Damil2014}{Damil2006,Damil2008,Damil2010}, is efficient as it can predict wrinkling patterns with few degrees of freedom, but it is inaccurate along boundaries. A remedy to that is to combine it with full shell models \cite{Huang2019a}. Detailed reviews of this approach are given by \fullcite{Potier-Ferry2016} and \fullcite{Huang2019a}. And recent works include the extension to cases with non-uniform wrinkling orientations \cite{Khalil2020}
and the combination with the ANM path following method \cite{Tian2021}. Although the Fourier reduced model provides a reduction in degrees of freedom, potentially independent of the wrinkling wave length, the applicability of the method has not been fully demonstrated yet.\\

Thirdly, when modelling wrinkles implicitly, the goal is not necessarily to establish the actual wrinkling pattern, but to estimate wrinkling sensitive parts of the structure instead. Such models are typically driven by a so-called \emph{tension field} \cite{Wagner1929,Wagner1931,Wagner1931,Reissner1938,Kondo1954,Mansfield1969,Mansfield1970}, describing the state of parts of a membrane by being either taut, slack, or wrinkled. Tension fields can be defined based on principal stresses, principal strains, or a combination of those, as discussed by \fullcitelist{Kang1997}{Kang1997,Kang1999} and \fullcite{Roddeman1987}, among others. Depending on the tension field, constitutive or kinematic equations can be modified to embed wrinkling effects into numerical methods; see the works of \fullcite{LeMeitour2021} and \fullcite{Miyazaki2006} for an overview. For example, \fullcitelist{Pipkin1986}{Pipkin1986,Pipkin1993} and \fullcite{Steigmann1990} proposed to modify the strain energy density function based on the tension field, for which variational methods \cite{Mosler2008,Mosler2008a,Mosler2009} and interior point models \cite{DeRooij2015} have been derived, and extensions for anisotropic and hyperelastic models have been proposed \cite{Massabo2006,Atai2012,Atai2014,Epstein2001}. In addition, material matrix modifications instead of strain energy density modifications have been developed \cite{Kang1997,Kang1999,Akita2007,Liu2001,Jarasjarungkiat2008,Jarasjarungkiat2009}. Alternatively, modification schemes based on deformation tensor modifications instead of constitutive relation modifications were proposed by \fullcitelist{Roddeman1987}{Roddeman1987,Roddeman1987a,Roddeman1991}, and applied to orthotropic materials \cite{Muttin1996,Raible2005}. Efficient implementations of this model were presented by \fullcite{Lu2001} and \fullcite{Nakashino2005}, the latter authors applying it to isogeometric membrane elements as well \cite{Nakashino2020}. Although the model based on modifications of the deformation tensor provides a more generic approach, it has not been applied to hyperelastic material modelling, to the best of the authors' knowledge. Although implicit models do not provide wrinkling amplitudes, post-processing methods for recovering wrinkling amplitudes have been proposed in the computer graphics community \cite{Skouras2014,Montes2020,Choi2002,Jin2017,Chen2021}. Finally, \fullcite{Zhang2024} recently presented a novel, variationally consistent wrinkling model for isogeometric membranes, employing a split of the strain tensor, restricted to linear elastic materials. The latter model does not require to evaluate a wrinkling criterion and to solve a root-finding problem at each quadrature point, as is the case for the model of \fullcite{Nakashino2020}.\\

In this paper, the wrinkling model of \fullcitelist{Nakashino2005}{Nakashino2005,Nakashino2020} is extended to hyperelastic materials. To this end, the model of \fullcite{Roddeman1987a} is used, where the deformation gradient is modified for a wrinkled material and, as a consequence, the stress and material tensors are modified based on the tension field in the membrane. The extension to hyperelastic materials introduces extra terms in the modified material tensor due to the dependency of the wrinkling stresses on the strain tensor. Although the formulations derived in the present work apply to finite element methods, they are applied in the context of isogeometric analysis \cite{Hughes2005}, as has been done by \fullcite{Nakashino2020}. Since finite bending stiffness is essential in the modelling of wrinkling amplitudes, reference solutions are computed using the isogeometric Kirchhoff-Love shell model \cite{Kiendl2009} with extensions to hyperelasticity \cite{Kiendl2015,Verhelst2021}.\\

The paper is outlined as follows: In \cref{sec:membrane_element}, the isogeometric membrane formulation is introduced. This element is equivalent to an isogeometric Kirchhoff-Love shell \cite{Kiendl2009} without bending contributions. \Cref{sec:elastic} recalls the original model by \fullcitelist{Nakashino2005}{Nakashino2005,Nakashino2020} for linear elastic materials to introduce the original concept. Thereafter, in \cref{sec:hyperelastic}, the model of \fullcite{Nakashino2020} is extended for hyperelastic materials. In \cref{sec:solution}, the numerical implementation of the model is discussed, and in \cref{sec:results}, numerical benchmark results are provided. Lastly, \cref{sec:conclusions} provides conclusions and recommendations for future work.

%% file: Sections_Membrane.tex
In this section, the isogeometric membrane formulation is derived. The primary purpose is to introduce the notations for geometric, kinematic, and constitutive quantities together with the variational formulation for a membrane. Since the membrane closely relates to parts of the Kirchhoff-Love shell, the notations are based on the ones used in \cite{Kiendl2015,Verhelst2021}. That is, basis vectors are denoted by lower-case bold and italic letters, e.g., $\avec$; second-order tensors are denoted by upper-case bold letters, e.g., $\TEN{A}$; discrete vectors or second-order tensors in Voigt notation are denoted by upper-case bold and italic letters, e.g., $\VEC{A}$; and matrices are denoted by upper-case letters, e.g., $A$. If needed, the notation $[A]$ is used to emphasise a matrix. Lastly, Greek sub- and superscripts take values of $1,2$, while Latin sub- and superscripts take values of $1,2,3$.\\

Consider surfaces $\Xvec(\theta^1,\theta^2)$ and $\xvec(\theta^1,\theta^2)=\Xvec(\theta^1,\theta^2)+\uvec(\theta^1,\theta^2)$ denoting points in the undeformed and deformed configurations of a membrane, respectively, with $\theta^\alpha$, $\alpha=1,2$ the parametric coordinates of the surface and $\uvec(\theta^1,\theta^2)$ the deformation vector field. Consequently, the covariant basis vectors of the deformed and undeformed configurations are defined by $\Avec_\alpha$ and $\avec_\alpha$, respectively, given by
\begin{equation}
	\Avec_{\alpha} = \pdv{\Xvec}{\theta^\alpha},
\end{equation}
and similar for the deformed configuration $\xvec$. In addition, $\mathring{a}_{\alpha\beta}=\Avec_\alpha\cdot\Avec_\beta$ are the coefficients of the covariant metric tensor. The vectors $\Avec^{\alpha}$ and $\avec^{\alpha}$ denote the contravariant basis vectors of the undeformed and deformed membrane surfaces, with identity $\Avec_\alpha\cdot\Avec^\beta=\delta_\alpha^\beta$ (similar for the basis vectors $\avec_\alpha$) with $\delta_\alpha^\beta=1$ if $\alpha=\beta$ and 0 otherwise. The contravariant basis vectors are constructed via the metric tensor:
\begin{equation}
	\Avec^{\alpha} = [\mathring{a}_{\alpha\beta}]^{-1}\Avec_\beta,
\end{equation}
where $[\mathring{a}]^{-1}$ denotes the inverse of the metric tensor coefficient matrix $[\mathring{a}]$. For the deformed configuration, the same relation holds between $\avec^{\alpha}$ and $\avec_{\alpha}$.

\begin{remark}\label{remark:thickness}
	Contrary to the coordinate system for Kirchhoff-Love shells \cite{Kiendl2009}, the through thickness coordinate $\theta^3$ is omitted for the isogeometric membrane, assuming small thickness. Therefore, it is assumed that deformations are constant through the thickness of the membrane.
\end{remark}

\subsection{Kinematic Equation}\label{subsec:kinematic}
The deformation gradient $\Ften$ or the deformation tensor $\Cten$ relate $\Xvec$ with $\xvec$ as follows:
\begin{equation}
	\Ften = \avec_\alpha\otimes\Avec^\beta,\quad\Cten = \Ften^\top\Ften = a_{\alpha\beta}\:\Avec^\alpha\otimes\Avec^\beta.
\end{equation}
Using these quantities, the Green-Lagrange strain tensor $\Eten$ is defined by
\begin{equation}\label{eq:Eten}
	\Eten = \frac{1}{2}\qty(\Ften^\top\Ften-\TEN{I}) = \frac{1}{2}\qty(\Cten-\TEN{I}) = \frac{1}{2}\qty(a_{\alpha\beta}-\mathring{a}_{\alpha\beta}) = E_{\alpha\beta}\:\Avec^\alpha\otimes\Avec^\beta,
\end{equation}
where $\TEN{I}=\mathring{a}_{\alpha\beta}\:\Avec^\alpha\otimes\Avec^\beta$ is the identity tensor. The second Piola-Kirchhoff stress tensor $\Sten$ is defined through the constitutive relation; see \cref{subsec:constitutive}.

\subsection{Constitutive Relation}\label{subsec:constitutive}
The relation between the Green-Lagrange strain tensor $\Eten$ to the second Piola-Kirchhoff stress tensor $\Sten$. In the following, the indices $i,j,k,l=1,\dots,3$ are used, representing the covariant and contravariant bases related to three parametric directions: parameters $\theta^1$ and $\theta^2$ represent the in-plane surface coordinates, and $\theta^3$ represents the through-thickness coordinate in the direction of the unit normal vector $\UVEC{a}_3$.\\

For linear elastic materials, the stress and strain tensors are simply related via the following relation:
\begin{equation}\label{eq:constitutive_linear}
	\Sten = \Ccten:\Eten = S^{ij}\:\Avec_i\otimes\Avec_j,
\end{equation}
where the coefficients of the stress tensor, $S^{ij}$, are given by $S^{ij} = \Ccoef^{ijkl}E_{kl}$. For a Saint Venant-Kirchhoff material with Lamé parameters $\lambda$ and $\mu$, the coefficients of the material tensor are given by:
\begin{equation}
	\Ccoef^{ijkl} = \frac{2\lambda\mu}{\lambda+2\mu}\mathring{a}^{ij}\mathring{a}^{kl} + \mu\qty(\mathring{a}^{ik}\mathring{a}^{jl} + \mathring{a}^{il}\mathring{a}^{jk}).
\end{equation}
For hyperelastic materials, the constitutive relationship is defined through the strain energy density function $\Psi(\Cten)$ or $\Psi(\Eten)$. In particular, the coefficients of the second Piola-Kirchhoff stress tensor are given by
\begin{equation}\label{eq:S}
	S^{ij} = 2\pdv{\Psi}{C_{ij}} = \pdv{\Psi}{E_{ij}}.
\end{equation}
The material tensor $\Ccoef^{ijkl}$ is not required in the derivation of the tensor $\Sten$; however, it plays a role in the definition of the variation of the stress tensor, $\delta \Sten$, as shown in \cref{eq:deltaS}. In terms of the strain energy density function, the coefficients of the material tensor are defined by
\begin{equation}\label{eq:C}
	C^{ijkl} = \pdv{S^{ij}}{C_{kl}} = 4\pdv[2]{\Psi}{C_{ij}}{C_{kl}}.
\end{equation}
As for the isogeometric Kirchhoff-Love shell, through thickness deformation is being neglected, meaning that $C_{33}=1$. This violates the plane stress condition ($S^{33}=0$) since $S^{33}=\pdv{\Psi}{C_{33}}\neq0$. To satisfy the plane stress condition, the normal deformation $C_{33}$ needs to be modified. As described for hyperelastic shells in \cite{Kiendl2015}, this can be done analytically for incompressible materials using the property that the Jacobian determinant, given by:
\begin{equation}\label{eq:J}
	J=\sqrt{\frac{\vert a_{\alpha\beta}\vert}{\vert \mathring{a}_{\alpha\beta}\vert}}\sqrt{C_{33}},
\end{equation}
is unity, i.e., $J=1$. For compressible materials, the plane stress condition is iteratively satisfied, as discussed in \cite{Kiendl2015}. Finally, static condensation of the material tensor results in the in-plane material tensor $\hat{\Ccten}^{\alpha\beta\gamma\delta}$\cite{Kiendl2015}:
\begin{equation}\label{eq:condensation}
	\hat{\Ccten}^{\alpha\beta\gamma\delta} = \Ccoef^{\alpha\beta\gamma\delta}-\frac{\Ccoef^{\alpha\beta33}\Ccoef^{33\gamma\delta}}{\Ccoef^{3333}}.
\end{equation}
For incompressible materials, it can be found that the coefficients of the statically condensed material tensor are \cite{Kiendl2015}:
\begin{multline}\label{eq:condensation_incompressible}
	\Ccoef^{\alpha\beta\gamma\delta} = 4\pdv[2]{\Psi}{C_{\alpha\beta}}{C_{\gamma\delta}} + 4\pdv[2]{\Psi}{C_{33}}J_0^{-4}a^{\alpha\beta}a^{\gamma\delta} - 4\pdv[2]{\Psi}{C_{33}}{C_{\alpha\beta}}J_0^{-2}a^{\gamma\delta}-4\pdv[2]{\Psi}{C_{33}}{C_{\gamma\delta}}J_0^{-2}a^{\alpha\beta} \\ +2\pdv{\Psi}{C_{33}}J_0^{-2}\qty(2a^{\alpha\beta}a^{\gamma\delta}+a^{\alpha\gamma}a^{\beta\delta}+a^{\alpha\delta}a^{\beta\gamma})
\end{multline}

Lastly, the tensor $\bm{\sigma}=J^{-1}\Ften^\top\Sten\Ften$ is the Cauchy stress tensor, which is used for stress recovery.

\begin{example}\label{example:NH_incompressible}
	For a Neo-Hookean material model with $\Psi(\Cten)=\frac{1}{2}\mu(I_1(\Cten)-3)$, where $\mu$ is Lamé's second parameter, \cref{eq:condensation_incompressible} simplifies to
	\begin{equation}\label{eq:NH_incompressible}
		\hat{\Ccten}^{\alpha\beta\gamma\delta} = \mu J_0^{-2}\qty(2a^{\alpha\beta}a^{\gamma\delta} + a^{\alpha\gamma}a^{\beta\delta} + a^{\alpha\delta}a^{\beta\gamma}).
	\end{equation}
\end{example}

\subsection{Variational Formulation}\label{subsec:variational}
The variational formulation for membranes follows from a variation of the internal and external energy contributions. For derivation of the variational formulation for isogeometric membranes, the reader is referred to \cite{Nakashino2020} or to works on the Kirchhoff-Love shell \cite{Kiendl2009,Kiendl2015} omitting the bending stiffness contributions. The elastic energy is given by
\begin{equation}\label{eq:Wint}
	\mathcal{W}_{\text{int}} = -\frac{1}{2}\int_{\Omega^\star} \Sten:\Eten \dd{\Omega^\star}.
\end{equation}
Here, $\Omega^\star=\Omega\times[-t/2,t/2]$ denotes the integration domain, with $t$ the thickness of the membrane and $\Omega$ the surface domain. Taking the Gateaux derivative with respect to the displacements $\uvec$, the variation of the internal elastic energy is:
\begin{equation}\label{eq:dWint}
	\delta \mathcal{W}_{\text{int}} = -\int_{\Omega^\star} \Sten:\delta\Eten \dd{\Omega^\star} = -\int_{\Omega^\star} \Nten:\delta\Eten \dd{\Omega^\star},
\end{equation}
where the tensor $\Nten$ denotes the membrane force tensor, obtained by integrating the stress tensor $\Sten$ through the thickness of the membrane. Since the thickness of the membrane is small (see \cref{remark:thickness}), the thickness integral becomes:
\begin{equation}\label{eq:Nten}
	\Nten(\theta^1,\theta^2)=\int_{[-t/2,t/2]} \Sten(\theta^1,\theta^2,\theta^3)\dd{\theta^3} \approx t\Sten(\theta^1,\theta^2,0).
\end{equation}
The external virtual work is provided by the following expression:
\begin{equation}\label{eq:dWext}
	\delta \mathcal{W}_{\text{ext}} = \int_{\Omega} \delta\uvec\cdot\VEC{f}(\uvec)\dd{\Omega} + \int_{\partial\Omega} \VEC{g}\cdot \delta\uvec\dd{\Gamma},
\end{equation}
where $\VEC{f}(\uvec)$ is a vector representing a follower-load acting on a point on the deformed surface $\xvec(\theta^1,\theta^2)$ and $\VEC{g}$ is a line load acting on the boundary $\partial\Omega$. When the sum of the internal and external virtual work is zero, i.e.
\begin{equation}\label{eq:dW}
	\delta \mathcal{W}(\uvec,\delta\uvec)=\delta \mathcal{W}_{\text{int}}-\delta \mathcal{W}_{\text{ext}}=0,
\end{equation}
equilibrium is found. Since $\delta \mathcal{W}(\uvec,\delta\uvec)$ is non-linear, solving the equation requires linearization. The second variation of the internal energy $\mathcal{W}_{\text{int}}$ in the system is given by:
\begin{equation}\label{eq:ddW}
	\delta^2 \mathcal{W}_{\text{int}}(\uvec,\delta\uvec,\Delta\uvec) = -\int_{\Omega} \delta\Nten:\delta\Eten + \Nten:\delta^2\Eten \dd{\Omega}.
\end{equation}
Here, $\delta^2\Eten$ denotes the second variation of the Green-Lagrange strain tensor. Since the external virtual work from \cref{eq:dWext} depends on the solution $\uvec$, its second variation is non-zero and given by:
\begin{equation}\label{eq:ddWext}
	\delta^2 \mathcal{W}_{\text{ext}}(\uvec,\delta\uvec,\Delta\uvec) = \int_{\Omega} \delta\uvec\cdot\VEC{f}^\prime(\uvec,\Delta\uvec) \dd{\Omega}.
\end{equation}
For a follower pressure, $\VEC{f}=p\UVEC{a}_3$, and its variation is $\VEC{f}^\prime=p\UVEC{a}_3^\prime(\uvec,\Delta\uvec)$. Here $\UVEC{a}_3(\uvec)$ and $\UVEC{a}_3^\prime(\uvec,\Delta\uvec)$ are the unit normal vector and its variation, respectively, which can be found in the derivation of the Kirchhoff--Love shell \cite{Kiendl2009}. Furthermore, the variation of the normal force tensor is:
\begin{equation}\label{eq:dNten}
	\delta\Nten = \int_{\Omega}\delta\Sten\dd{\Omega}.
\end{equation}
Here, the variation of the second Piola-Kirchhoff stress tensor is defined as:
\begin{equation}\label{eq:deltaS}
	\delta\Sten = \Ccten:\delta\Eten.
\end{equation}
\Cref{eq:ddW,eq:deltaS} show that given the variations $\delta \Eten$ and $\delta^2\Eten$, the first and second variations of the internal energy $\delta W_{\text{int}}$ and $\delta^2 \mathcal{W}_{\text{int}}$ from \cref{eq:dW,eq:ddW} can be found. Since the deformed configuration $\xvec$ is unknown, the problem needs to be discretised such that the variations of $\Eten$ can be defined.\\

\subsection{Discretisation}
Discretisation of the variational problem in \cref{eq:dW} is achieved by discretising the undeformed and deformed configurations $\Xvec$ and $\xvec$, respectively. In the context of isogeometric analysis, this is done by choosing splines as a basis for the geometry, i.e., by describing the geometries as a weighted sum of basis functions $\varphi_i(\theta^1,\theta^2)$ and control points $\Xvec^h_i$ and $\xvec^h_i$, respectively:
\begin{equation}\label{eq:xnurbs}
	\begin{aligned}
		\Xvec^h(\theta^1,\theta^2) = \sum_k \varphi_k(\theta^1,\theta^2) \Xvec^h_k,\\
		\xvec^h(\theta^1,\theta^2) = \sum_k \varphi_k(\theta^1,\theta^2) \xvec^h_k.
	\end{aligned}
\end{equation}
The superscript $h$ indicates a discretisation of the undeformed and deformed geometries, $\Xvec$ and $\xvec$, respectively. Using the same basis for $\Xvec^h$ and $\xvec^h$, the discrete displacement vector $\uvec^h$ is given as the difference between the two, i.e., $\uvec^h=\xvec^h-\Xvec^h$. Since the variational formulation is expressed in terms of the displacement field $\uvec$, the variations of the control points of the field $\uvec^h$ are the virtual displacements in the discrete system, hence the unknowns. In the sequel, the subscripts $r$ and $s$ denote the indices of the components of conveniently numbered degrees of freedom of $\uvec^h$ incorporating the spatial dimensions of the surface. Furthermore, the notation $(\cdot)_r=\partial (\cdot)\partial u_r$ is used for derivatives, and the superscript $h$ is omitted. Following from \cref{eq:xnurbs}, the variation of the deformed geometry is given by
\begin{equation}\label{eq:dx}
	\xvec_{,r} = \sum_k \qty(\Xvec_{k,r} + \uvec_{k,r}) = \sum_k \varphi_k\uvec_{k,r} = \uvec_{,r},
\end{equation}
where the last equality follows from the fact that the undeformed configuration is trivially independent of the deformation field $\uvec$. Similarly, the derivatives of the covariant basis vectors $\avec_\alpha$ of the discrete deformed configuration $\xvec_h$, see \cref{eq:dx}, are:
\begin{equation}
	\avec_{\alpha,r} = \qty(\pdv{\xvec_k}{\theta^{\alpha}})_{,r} = \sum_k \pdv{\varphi_k}{\theta^\alpha} \uvec_{k,r}.
\end{equation}
As a consequence, the variation of the surface metric tensor of the deformed configuration, $a_{\alpha\beta}$, becomes:
\begin{equation}
	a_{\alpha\beta,r} = \qty(\avec_\alpha\cdot\avec_\beta)_{r} = \avec_{\alpha,r}\cdot\avec_\beta + \avec_\alpha\cdot\avec_{\beta,r}.
\end{equation}
Since the undeformed configuration is invariant to the deformation field $\uvec$, the first variation of the membrane strain tensor $\epsten$ from \cref{eq:Eten} becomes
\begin{equation}\label{eq:Evar}
	E_{\alpha\beta,r} = \frac{1}{2}a_{\alpha\beta,r}.
\end{equation}
Similarly, the second variation of the deformed configuration, the deformed surface metric tensor, and the membrane strain can be derived. Starting with the first variation of the deformed configuration from \cref{eq:dx}, the second variation becomes
\begin{equation}
	\xvec_{,rs} = \sum_k \varphi_k \uvec_{k,rs} = \VEC{0}.
\end{equation}
The second variation of $\uvec_k$ is zero since the components of these nodal weights are linear in $u_r$. Similarly, $\avec_{\alpha,rs}=\VEC{0}$. As a consequence, the second variation of the surface metric tensor in the deformed configuration, $a_{\alpha\beta}$, becomes
\begin{equation}
	a_{\alpha\beta,rs} =
	\avec_{\alpha,rs}\cdot\avec_{\beta} + \avec_{\alpha,r}\cdot\avec_{\beta,s} + \avec_{\alpha,s}\cdot\avec_{\beta,s} + \avec_{\alpha}\cdot\avec_{\beta,rs},
	=
	\avec_{\alpha,r}\cdot\avec_{\beta,s} + \avec_{\alpha,s}\cdot\avec_{\beta,s}.
\end{equation}
Again, since the undeformed configuration is invariant to the deformation field $\uvec$, the second variation of the membrane strain tensor becomes
\begin{equation}\label{eq:Evar2}
	E_{\alpha\beta,rs} = \frac{1}{2}a_{\alpha\beta,rs}
\end{equation}

Since $E_{\alpha\beta,r}$ and $E_{\alpha\beta,rs}$ are the discrete variations of $\delta\Eten$ and $\delta^2\Eten$, respectively, the discrete residual and Jacobian can be derived. The discrete residual vector follows from \cref{eq:dW} with \cref{eq:dWext,eq:dWint,eq:Evar,eq:Nten}:
\begin{equation}\label{eq:residual}
	R_r(\uvec) = \int_\Omega \Nten(\uvec):\Eten_{,r}(\uvec)\dd{\Omega} - \int_{\Omega}\VEC{f}(\uvec)\cdot\uvec_{,r}\dd{\Omega}.
\end{equation}
Furthermore, the discrete Jacobian matrix follows from \cref{eq:ddW} with \cref{eq:ddWext,eq:dWint,eq:Evar,eq:Nten,eq:dNten,eq:Evar2}:
\begin{equation}\label{eq:jacobian}
	K_{rs} = \int_{\Omega} \Nten_{,s}(\uvec):\Eten_{,r}(\uvec) + \Nten(\uvec):\Eten_{,rs}(\uvec)\dd{\Omega} - \int_{\Omega} \VEC{f}(\uvec)_{,s}\cdot\uvec_{,r}\dd{\Omega}.
\end{equation}
Here, the product $\TEN{A}:\TEN{B}$ denotes an inner product of two second-order tensors. Furthermore, the contribution of the displacement-dependent load $\VEC{f}(\uvec)$ requires a derivative, defined as $\VEC{f}(\uvec)_{,s}=p\UVEC{a}_{3,s}$ for a follower pressure, where $\UVEC{a}_{3,s}(\uvec)$ is the discrete derivative of the surface normal vector $\UVEC{a}_3(\uvec)$, which can be found in the derivations of the Kirchhoff--Love shell \cite{Kiendl2009}. The contribution $\VEC{f}(\uvec)\cdot\uvec_{,rs}=0$ since $\uvec_{,rs}=\VEC{0}$ \cite{Kiendl2011}. Given the discretisation of $\xvec$ and $\Xvec$ using splines, see \cref{eq:xnurbs}, and given the residual and Jacobian from \cref{eq:residual,eq:jacobian}, respectively, only the definition of the stress tensor $\Sten$ and the material tensor $\Ccten$ are remaining undefined. In \cref{sec:elastic,sec:hyperelastic}, definitions for $\Sten$ and $\Ccten$ are provided, incorporating the wrinkling model of \fullcite{Nakashino2005} for linear elasticity and incorporating the extension for hyperelastic materials, which is the novelty of this paper.

\subsection{Implementation}
Since the tensors $\Eten$ and $\Sten$ are symmetric second-order tensors, they can be written in Voigt notation:
\begin{equation}\label{eq:voigt}
	\VEC{S} =
	\begin{bmatrix}
		S^{11} & S^{22} & S^{12}
	\end{bmatrix}^\top,
	\quad
	\VEC{E} =
	\begin{bmatrix}
		E^{11} & E^{22} & 2E^{12}
	\end{bmatrix}^\top.
\end{equation}
As a consequence, the material tensor is represented in Voigt notation as well, using \cref{eq:deltaS}:
\begin{equation}\label{eq:Cvoigt}
	\VEC{C} =
	\begin{bmatrix}
		\Ccoef^{1111} & \Ccoef^{1122} & \Ccoef^{1112}\\
		\Ccoef^{2211} & \Ccoef^{2222} & \Ccoef^{2212}\\
		\Ccoef^{1211} & \Ccoef^{1222} & \Ccoef^{1212}
	\end{bmatrix}.
\end{equation}
Using Voigt notation for the tensors $\Eten$, $\Sten$, and $\Ccten$, the second-order tensor inner products $:$ in \cref{eq:deltaS,eq:residual,eq:jacobian} can simply be evaluated as matrix-vector products and vector inner-products. In \cref{subsec:elastic_implementation,subsec:hyperelastic_implementation}, the elastic and hyperelastic constitutive laws for wrinkling are provided in Voigt notations for fast computer implementation.\\

%% file: Sections_Linear.tex
In the linear elastic wrinkling model, the taut, wrinkling, and slack conditions are governed by a modification of the deformation gradient and consequently the strain, stress, and material tensors, following the model originally proposed by \fullcitelist{Roddeman1987}{Roddeman1987,Roddeman1987a}. This model lays the foundation of the wrinkling model proposed by \fullcite{Nakashino2005}, which was later implemented for isogeometric membranes \cite{Nakashino2020}. This model is presented with the assumption that the material behaviour is linear; therefore, it is referred to herein as the linear elastic wrinkling model. In this section, an overview of the linear elastic wrinkling model proposed by \fullcite{Nakashino2005} is provided. In \cref{sec:hyperelastic}, the linear elastic wrinkling model will be extended to hyperelastic material models.\\

In general, tension field-based models rely on the definition of a tension field $\phi$ over the domain. A tension field classifies the stress state in a membrane as either slack, taut, or wrinkled, depending on the deformation tensor $\VEC{C}$. In general, three different definitions of the tension field are used in the literature: based on principal strains, principal stresses, or combinations of those (mixed), in this paper represented by $\phi^{E}$, $\phi^{S}$, and $\phi^{M}$, respectively. These tension fields are defined as:
\begin{multline}\label{eq:tensionfield}
	\phi^{E} =
	\begin{dcases}
		\text{Taut} & \text{if}\:E_{p,1} > 0\\
		\text{Slack} & \text{if}\:E_{p,2} \leq 0\\
		\text{Wrinkled} & \text{otherwise}\\
	\end{dcases},
	\qquad
	\phi^{S} =
	\begin{dcases}
		\text{Taut} & \text{if}\:S_{p,1} > 0\\
		\text{Slack} & \text{if}\:S_{p,2} \leq 0\\
		\text{Wrinkled} & \text{otherwise}\\
	\end{dcases},
	\qquad
	\phi^{M} =
	\begin{dcases}
		\text{Taut} & \text{if}\:S_{p,1} > 0\\
		\text{Slack} & \text{if}\:E_{p,2} \leq 0\\
		\text{Wrinkled} & \text{otherwise}\\
	\end{dcases},
\end{multline}
where $S_{p,1}$ and $S_{p,2}$ are the principal stresses such that $S_{p,1}\leq S_{p,2}$ and $E_{p,1}$ and $E_{p,2}$ are the principal strains such that $E_{p,1}\leq E_{p,2}$. Given the tension field, tension field models typically modify the stress tensor $\Sten$ and consequently the material tensor $\Ccten$ based on the tension field:
\begin{equation}\label{eq:S_TFT}
	\Sten =
	\begin{dcases}
		\VEC{0} & \text{if } \phi = \text{ Slack}\\
		\Sten & \text{if } \phi = \text{ Taut}\\
		\Sten^\prime & \text{if } \phi = \text{ Wrinkled}
	\end{dcases}
	,\quad
	\Ccten =
	\begin{dcases}
		\VEC{0} & \text{if } \phi = \text{ Slack}\\
		\Ccten & \text{if } \phi = \text{ Taut}\\
		\Ccten^\prime & \text{if } \phi = \text{ Wrinkled}
	\end{dcases},
\end{equation}
where $\Sten^\prime$ is a modified stress tensor. This modified stress tensor can be obtained in different ways, either by adjusting the constitutive or kinematic equations provided by the tension field. For example, the work of \cite{Alberini2021} modifies the Ogden constitutive relation based on the tension field $\phi^S$, whereas \cite{Nakashino2005} modifies the kinematic equation based on $\phi^M$. In this paper, the approach of \cite{Nakashino2005} is followed. As discussed in the work of \fullcite{Kang1997}, the definition of the tension field using $\phi^M$ has advantages over $\phi_E$ and $\phi_S$.\\

\subsection{Kinematic Equation}
Given the deformation gradient $\Ften$, the modified deformation gradient \cite{Roddeman1987,Roddeman1987a} is given by
\begin{equation}
	\Ften^\prime = (\TEN{I} + b\UVEC{w}\otimes\UVEC{w})\cdot \Ften.
\end{equation}
Here, $\TEN{I}$ is the second-order identity tensor, and $b$ is the measure of the amount of `wrinkliness' \cite{Nakashino2005,Roddeman1987,Lu2001}, by definition $b>0$. Furthermore, $\UVEC{w}$ is the unit vector transverse to the wrinkles. Using the modified deformation tensor, the modified strain tensor can be computed, given by:
\begin{equation}\label{eq:wrinklingstrain}
	\Eten^\prime = E_{\alpha\beta}^\prime \Avec^\alpha\otimes \Avec^\beta = \frac{1}{2}\qty(\Ften^{\prime\top}\cdot \Ften^\prime -\TEN{I}) = \Eten + \Eten_W,
\end{equation}
where $\Eten=\Ften^\top\Ften-\TEN{I}=\Cten-\TEN{I}$ is the Green-Lagrange strain with
\begin{equation}\label{eq:deformationtensor}
	\Cten=C_{\alpha\beta}\:\Avec^\alpha\otimes\Avec^\beta = a_{\alpha}a_{\beta}\:\Avec^\alpha\otimes\Avec^\beta,
\end{equation}
being the deformation tensor. $\Eten_W$ is the wrinkling strain, given by
\begin{equation}
	\Eten_W = \frac{1}{2}b\qty(b+2)\UVEC{w}\otimes\UVEC{w}.
\end{equation}
Here, $\UVEC{w}=\VEC{w}\cdot\Ften=w_a\Avec^\alpha$ is the projection of $\UVEC{w}$ onto the undeformed contravariant basis. Introducing the rotation $\vartheta$ and magnitude $a\neq0$ of the projected wrinkling direction $\UVEC{w}$ using $\UVEC{w}_1=a n_1$ and $\UVEC{w}_2=a n_2$ using $n_1=\cos\vartheta$ and $n_2=\sin\vartheta$, the coefficients of the wrinkling strain tensor can be written as
\begin{equation}\label{eq:wrinklingstrain2}
	E_{\alpha\beta}^\prime = E_{\alpha\beta} + \frac{1}{2}\qty(\Ften^{\prime\top} \Ften^\prime -\TEN{I}) = E_{\alpha\beta} + \gamma n_\alpha n_\beta,
\end{equation}
where $\gamma=\frac{1}{2}a^2b(b+1)$. In this definition of the strain tensor $\Eten^\prime$, the wrinkling strain amplitude $\gamma$ and the angle of the wrinkles $\vartheta$ are unknown. Through the uniaxial tension condition for wrinkled materials, these unknowns will be determined in the next subsection.

\subsection{Constitutive Relation}
Since the definition of the strain tensor $\Eten$ changes for the wrinkled state of the membrane, the constitutive relation from \cref{eq:constitutive_linear} also changes. For linear materials, the wrinkled stress tensor $\Sten^\prime$ simply becomes:
\begin{equation}
	\Sten^\prime = \Ccten:\Eten^\prime = \Ccten:\qty(\Eten + \Eten_W),
\end{equation}
where $\Eten_W$ depends on the unknowns $\gamma$ and $\vartheta$, as in \cref{eq:wrinklingstrain2}. This equation can also be written in terms of the component of the stress tensor $\Sten^\prime=S^{\prime\alpha\beta}\:\Avec_\alpha\otimes\Avec_\beta$:
\begin{equation}\label{eq:wrinklingstress}
	S^{\prime\alpha\beta} = \Ccoef^{\alpha\beta\gamma\delta}(E_{\gamma\delta} + \gamma n_\gamma n_\delta) = S^{\alpha\beta} + \gamma\Ccoef^{\alpha\beta\gamma\delta} n_\gamma n_\delta.
\end{equation}
For a wrinkled membrane, an uniaxial tension state is assumed \cite{Liu2001}, meaning the stress orthogonal to the wrinkles should vanish,
\begin{equation}\label{eq:uniaxialtension_cauchy}
	\sigma^\prime\cdot\VEC{w} = 0,
\end{equation}
and the stress parallel to the wrinkles should be positive
\begin{equation}\label{eq:uniaxialtension_cauchy_parallel}
	\VEC{t}\cdot\sigma^\prime\cdot\VEC{t} > 0.
\end{equation}
Here, $\sigma^\prime=\frac{1}{\det \Ften^\prime}\Ften^\prime \Sten^\prime \Ften^{\prime\top}$ is the modified Cauchy stress tensor. The uniaxial tension condition in \cref{eq:uniaxialtension_cauchy} can be written as
\begin{equation}\label{eq:uniaxialtension_cauchyS}
	\Sten^\prime\cdot\VEC{w} = 0,
\end{equation}
using the components $n_\alpha$ and $n_\beta$, the conditions \cref{eq:uniaxialtension_cauchy,eq:uniaxialtension_cauchy_parallel} can be written as \cite{Lu2001,Nakashino2005}:
\begin{equation}\label{eq:uniaxialtension_cauchy2}
	\begin{aligned}
		S^{\prime\alpha\beta} n_\alpha n_\beta = 0,\\
		S^{\prime\alpha\beta} m_\alpha n_\beta = 0,\\
		E^{\prime\alpha\beta} m_\alpha m_\beta > 0.\\
	\end{aligned}
\end{equation}
Here, $m_\alpha = \partial n_\alpha/\partial \vartheta$. From \cref{eq:wrinklingstress}, the uniaxial tension condition \cref{eq:uniaxialtension_cauchy2} becomes:
\begin{equation}\label{eq:uniaxialtension_cauchy3}
	\begin{aligned}
		S^{\alpha\beta} n_\alpha n_\beta + \gamma\Ccoef^{\alpha\beta\gamma\delta}n_\alpha n_\beta n_\gamma n_\delta = 0,\\
		S^{\alpha\beta} m_\alpha n_\beta + \gamma\Ccoef^{\alpha\beta\gamma\delta}m_\alpha n_\beta n_\gamma n_\delta = 0.
	\end{aligned}
\end{equation}
From the first line of \cref{eq:uniaxialtension_cauchy3}, the variable $\gamma$ can be found as:
\begin{equation}\label{eq:gamma}
	\gamma = - \frac{S^{\alpha\beta} n_\alpha n_\beta}{ \Ccoef^{\alpha\beta\gamma\delta}n_\gamma n_\delta n_\alpha n_\beta },
\end{equation}
and substituting $\gamma$ in the second line of \cref{eq:uniaxialtension_cauchy3}, the following equation is found:
\begin{equation}\label{eq:f}
	f(\vartheta) \equiv S^{\alpha\beta} m_\alpha n_\beta +\gamma \Ccoef^{\alpha\beta\gamma\delta} m_\alpha n_\beta n_\gamma n_\delta = 0.
\end{equation}
The only unknown in this equation is the angle $\vartheta$, thus the equation can be solved by root finding. The root $\vartheta$ of $f(\vartheta)=0$ must satisfy the uniaxial tension conditions \cref{eq:uniaxialtension_cauchy,eq:uniaxialtension_cauchy_parallel}. Since $f(\vartheta)$ follows from \cref{eq:uniaxialtension_cauchy}, the first uniaxial tension condition is satisfied when the root is found. The condition for positive stress along the wrinkling direction, i.e., \cref{eq:uniaxialtension_cauchy_parallel}, is satisfied by selecting the feasible root. In the work of \cite{Liu2001}, the bounds of an interval for the root that satisfies \cref{eq:uniaxialtension_cauchy_parallel} are derived, such that a bounded root finding algorithm, e.g., Brent's method \cite{Brent1971}, can be used. If the procedure to find the root $\vartheta$ fails on the prescribed interval, the domain $[0,2\pi]$ can be subdivided into sub-intervals $[\vartheta_A,\vartheta_B)\subset[0,2\pi]$, and Brent's method can be started for each sub-interval $[\vartheta_A,\vartheta_B)\subset[0,2\pi]$ that satisfies $f(\vartheta_A)f(\vartheta_B)<0$. The root $\vartheta=0$ is a root if and only if $\lim_{\vartheta\to 0}f(\vartheta)f(2\pi-\vartheta)<0$. As soon as the root $\vartheta$ is found, the wrinkling strain $\Eten^\prime$ from \cref{eq:wrinklingstrain2} and the wrinkling stress $\Sten^\prime$ from \cref{eq:wrinklingstress} can be computed.\\

\subsection{Variational Formulation}
In the variational equation from \cref{eq:dW}, the variation of the strain tensor $\Eten$ and stress tensor $\Sten$, respectively $\delta\Eten$ and $\delta\Sten$, are required. For the taut state, the variations remain unchanged, as seen in \cref{eq:S_TFT}. In the slack state, the stress tensor and its variation are equal to zero, hence the variation of the internal energy becomes zero. For the wrinkling state, the variations of $\Eten^\prime$ and $\Sten^\prime$ need to be found.\\

Firstly, the variation of $\Eten^\prime$ is independent of the constitutive law \cite{Nakashino2005}. Here, it is shown that the contribution of the virtual wrinkling strain, $\delta\Eten_W$, in the variational formulation is zero, since the product of the wrinkling stress tensors $\Sten^\prime$ and $\delta\Eten_W$ is zero. Physically, this means that the wrinkling strain corresponds to the rigid body movements to stretch the wrinkled membrane, hence not altering the strain energy.\\

Secondly, the variation of $\Sten^\prime$ needs to be found. Taking the variation of \cref{eq:wrinklingstress}, it follows that:
\begin{equation}\label{eq:dS}
	\begin{aligned}
		\delta S^{\prime\alpha\beta}(\Eten) = \dv{S^{\prime\alpha\beta}}{E_{\sigma\tau}} \delta E_{\sigma\tau} &= \Ccoef^{\prime\alpha\beta\sigma\tau} \delta E_{\sigma\tau}\\
		&= \qty(\pdv{S^{\alpha\beta}}{E_{\sigma\tau}} + \dv{\gamma}{E_{\sigma\tau}}\Ccoef^{\alpha\beta\gamma\delta} n_\gamma n_\delta + \gamma\Ccoef^{\alpha\beta\gamma\delta} \dv{(n_\gamma n_\delta)}{E_{\sigma\tau}})\delta E_{\sigma\tau}.
	\end{aligned}
\end{equation}
The full derivative of $\gamma$ with respect to $E_{\sigma\tau}$ can be found using the definition of $\gamma$ from \cref{eq:gamma}:
\begin{equation}
	\dv{\gamma}{E_{\sigma\tau}} = \pdv{\gamma}{E_{\sigma\tau}} + \pdv{\gamma}{\vartheta}\pdv{\vartheta}{E_{\sigma\tau}}.
\end{equation}

The derivative of $(n_\gamma n_\delta)$ to $E_{\sigma\tau}$ directly follows from the definitions of $n_\alpha$, $m_\alpha$ and the chain rule,
\begin{equation}\label{eq:dnndT}
	\dv{(n_\gamma n_\delta)}{E_{\sigma\tau}} = \pdv{(n_\gamma n_\delta)}{\vartheta} \pdv{\vartheta}{E_{\sigma\tau}} = \qty(m_\gamma n_\delta + n_\gamma m_\delta)\pdv{\vartheta}{E_{\sigma\tau}}.
\end{equation}
The derivative of $\gamma$ with respect to $\vartheta$ follows from \cref{eq:gamma}:
\begin{equation}\label{eq:dgammadT}
	\begin{aligned}
		\pdv{\gamma}{\vartheta} &= -\frac{\qty(\Ccoef^{\alpha\beta\gamma\delta}n_\gamma n_\delta n_\alpha n_\beta) \qty(S^{\alpha\beta} \pdv{(n_\alpha n_\beta)}{\vartheta}) - \qty(S^{\alpha\beta} n_\alpha n_\beta) \qty(\Ccoef^{\alpha\beta\gamma\delta}\qty(\pdv{(n_\gamma n_\delta)}{\vartheta} n_\alpha n_\beta + n_\gamma n_\delta \pdv{(n_\alpha n_\beta)}{\vartheta}))}{\qty(\Ccoef^{\alpha\beta\gamma\delta}n_\gamma n_\delta n_\alpha n_\beta)^2}\\
		&= -\frac{\gamma \Ccoef^{\alpha\beta\gamma\delta}n_\gamma n_\delta \pdv{(n_\alpha n_\beta)}{\vartheta}}{\Ccoef^{\alpha\beta\gamma\delta}n_\gamma n_\delta n_\alpha n_\beta},
	\end{aligned}
\end{equation}
where \cref{eq:dnndT,eq:uniaxialtension_cauchy3,eq:gamma} are used in the second equality. Furthermore, the derivative $\pdv{\gamma}{E_{\alpha\beta}}$ follows directly from \cref{eq:gamma}
\begin{equation}\label{eq:dgammadE}
	\pdv{\gamma}{E_{\sigma\tau}} = - \frac{\Ccoef^{\alpha\beta\sigma\tau} n_\alpha n_\beta}{ \Ccoef^{\alpha\beta\gamma\delta}n_\gamma n_\delta n_\alpha n_\beta }.
\end{equation}
Lastly, the derivative of the angle $\vartheta$ with respect to the strain tensor component $E_{\sigma\tau}$ can be found by using \cref{eq:f}
\begin{equation}\label{eq:dfdE}
	\pdv{f}{E_{\sigma\tau}} + \pdv{f}{\vartheta}\pdv{\vartheta}{E_{\sigma\tau}} = 0,
\end{equation}
giving:
\begin{equation}\label{eq:dTdE}
	\pdv{\vartheta}{E_{\sigma\tau}} = -\pdv{f}{E_{\sigma\tau}}\qty(\pdv{f}{\vartheta})^{-1}.
\end{equation}
The derivative of $f$ with respect to $\Eten$ directly follows from \cref{eq:f}:
\begin{equation}\label{eq:dFdE}
	\pdv{f}{E_{\sigma\tau}} = \Ccoef^{\alpha\beta\sigma\tau}m_\alpha n_\beta - \pdv{\gamma}{E_{\sigma\tau}}\Ccoef^{\alpha\beta\gamma\delta}m_\gamma n_\delta n_\alpha n_\beta,
\end{equation}
and the derivative of $f$ with respect to $\vartheta$ follows from \cref{eq:f} as well:
\begin{equation}\label{eq:dFdT}
	\pdv{f}{\vartheta} = S^{\alpha\beta} \pdv{(m_\alpha n_\beta)}{\vartheta} + \pdv{\gamma}{\vartheta} \Ccoef^{\alpha\beta\gamma\delta} m_\alpha n_\beta n_\gamma n_\delta + \gamma \Ccoef^{\alpha\beta\gamma\delta} \qty(\pdv{(m_\alpha n_\beta)}{\vartheta} n_\gamma n_\delta + m_\alpha n_\beta \pdv{(n_\gamma n_\delta)}{\vartheta}),
\end{equation}
where the derivative of $m_\alpha n_\beta$ with respect to $\vartheta$ can be easily obtained from the definitions of $m_\alpha$ and $n_\alpha$:
\begin{equation}\label{eq:dmndT}
	\pdv{(m_\alpha n_\beta)}{\vartheta} = -n_\alpha n_\beta + m_\alpha m_\beta.
\end{equation}

Using \cref{eq:dgammadE,eq:dgammadT,eq:dfdE,eq:dTdE,eq:dFdE,eq:dFdT}, the variation of $\Sten$ from \cref{eq:dS} is found, as well as the definition of the modified material tensor $\Ccten^\prime$.

\subsection{Implementation}\label{subsec:elastic_implementation}
As discussed in the work of \fullcite{Nakashino2005}, the modifications of the stress and material tensor for the wrinkling model, see \cref{eq:wrinklingstress,eq:dS}, can be expressed in terms of matrix-vector multiplications when employing Voigt notation for the stress and strain tensors, $\Sten$ and $\Eten$, respectively, see \cref{eq:voigt,eq:Cvoigt}. Following the notation of \citeauthor{Nakashino2005} \cite{Nakashino2005}, the terms $n_\alpha$ and $m_\alpha$ are collected in the following vectors:
\begin{equation}\label{eq:n}
	\begin{aligned}
		\VEC{n}_1 &=
		\begin{bmatrix}
			n_1n_1 & n_2n_2 & 2n_1n_2
		\end{bmatrix}^\top, \\
		\VEC{n}_2 &=
		\begin{bmatrix}
			n_1m_1 & n_2m_2 & n_1m_2 + m_1n_2
		\end{bmatrix}^\top = \frac{1}{2}\pdv{\VEC{n}_1}{\vartheta}\\
		\VEC{n}_3 &=
		\begin{bmatrix}
			m_1m_1 - n_1n_1 & m_2m_2-n_2n_2 & 2\qty(m_1m_2 -n_1n_2)
		\end{bmatrix}^\top = \pdv{\VEC{n}_2}{\vartheta},\\
		\VEC{n}_4 &=
		\begin{bmatrix}
			m_1m_1 & m_2m_2 & 2m_1m_2
		\end{bmatrix}^\top = \VEC{n}_3+\VEC{n}_1.
	\end{aligned}
\end{equation}
Using \cref{eq:n}, the wrinkled stress and strain tensors from \cref{eq:wrinklingstrain2,eq:wrinklingstress,eq:wrinklingstrain2} are written as
\begin{align}
	\VEC{S}^\prime &= \VEC{S} + \gamma\VEC{C}\cdot\VEC{n}_1,\\
	\VEC{E}^\prime &= \VEC{E} + \gamma \VEC{n}_1.
\end{align}
Furthermore, using the stress and strain tensor coefficients (\cref{eq:voigt}) together with \cref{eq:n}, the formulations in \cref{eq:uniaxialtension_cauchy2} become:
\begin{equation}
	\begin{aligned}
		\VEC{S}^\prime \cdot \VEC{n}_1 &= 0,\\
		\VEC{S}^\prime \cdot \VEC{n}_2 &= 0,\\
		\VEC{E}^\prime \cdot \VEC{n}_4 &> 0.\\
	\end{aligned}
\end{equation}
Following from these relations, the equations \cref{eq:gamma} and \cref{eq:f} are written as
\begin{align}
	\gamma &= -\frac{\VEC{S}\cdot\VEC{n}_1}{\VEC{n}_1^\top \cdot \VEC{C}\cdot\VEC{n}_1},\\
	f(\vartheta) &= \VEC{S}\cdot\VEC{n}_1 - \gamma \VEC{n}_2^\top \cdot \VEC{C}\cdot\VEC{n}_1.
\end{align}
To compute the wrinkling material tensor $\VEC{C}^\prime$, the derivatives from \cref{eq:dgammadT,eq:dgammadE,eq:dfdE,eq:dfdE} need to be expressed in terms of $\VEC{E}$, $\VEC{S}$, $\VEC{C}$, and $\VEC{n}_k$, $k=1,\dots,4$; see \cref{eq:voigt,eq:Cvoigt,eq:n}. The derivative of $\gamma$ with respect to $\vartheta$ from \cref{eq:dgammadT} becomes:
\begin{equation}\label{eq:dgammadT_voigt}
	\pdv{\gamma}{\vartheta} = -2\gamma\frac{\VEC{n}_2^\top\cdot\VEC{C}\cdot\VEC{n}_1}{\VEC{n}_1^\top\cdot\VEC{C}\cdot\VEC{n}_1}.
\end{equation}
Furthermore, the derivative of $\gamma$ with respect to $\VEC{E}$ becomes (see \cref{eq:dgammadE}):
\begin{equation}\label{eq:dgammadE_voigt}
	\qty[\pdv{\gamma}{\VEC{E}}] = -\frac{\VEC{C}\cdot\VEC{n}_1}{\VEC{n}_1^\top\cdot\VEC{C}\cdot\VEC{n}_1},
\end{equation}
where the bracket $[\cdot]$ is used to stress that $\pdv{\gamma}{\VEC{E}}$ is a vector with the derivatives of $\gamma$ with respect to $E_{\alpha\beta}$. To obtain the derivative of $\vartheta$ with respect to $\VEC{E}$, the derivatives of $f$ with respect to $\vartheta$ and $\VEC{E}$ are used, as in \cref{eq:dTdE}. The derivative of $f$ with respect to $\VEC{E}$ is a vector as well. Following from \cref{eq:dfdE}, it is given by:
\begin{equation}\label{eq:dfdE_voigt}
	\qty[\pdv{f}{\VEC{E}}] = \VEC{C}\cdot\VEC{n}_2 + \qty(\VEC{n}_2^\top\cdot\VEC{C}\cdot\VEC{n}_1)\qty[\pdv{\gamma}{\VEC{E}}],
\end{equation}
and the derivative of $f$ with respect to $\vartheta$ is obtained from \cref{eq:dFdT}:
\begin{equation}\label{eq:dfdT_voigt}
	\pdv{f}{\vartheta} = (\VEC{n}_4^\top\cdot\Sten) + \pdv{\gamma}{\vartheta}\qty(\VEC{n}_2^\top\cdot\VEC{C}\cdot\VEC{n}_1) + \gamma \qty(\VEC{n}_4^\top\cdot\VEC{C}\cdot\VEC{n}_1 + 2\VEC{n}_2^\top\cdot\VEC{C}\cdot\VEC{n}_2).
\end{equation}
Using \cref{eq:dfdE_voigt,eq:dfdT_voigt}, the derivative of $\vartheta$ with respect to $\VEC{E}$ is simply obtained by scalar division, using \cref{eq:dTdE}:
\begin{equation}\label{eq:dTdE_voigt}
	\qty[\pdv{\vartheta}{\VEC{E}}] = \qty[\pdv{f}{\VEC{E}}] / \pdv{f}{\vartheta}.
\end{equation}

Using the definition of $\Ccoef^\prime$ from \cref{eq:dS}, the matrix $\VEC{C}^\prime$ can be expressed in linear algebra operations using the scalars and vectors defined in \cref{eq:dgammadT_voigt,eq:dgammadE_voigt,eq:dTdE_voigt}:
\begin{equation}\label{eq:C_voigt}
	\VEC{C}^\prime = \VEC{C} \qty( \TEN{I} + \VEC{n}_1\otimes\qty[\pdv{\gamma}{\VEC{E}}] + \pdv{\gamma}{\vartheta}\VEC{n_1}\otimes\qty[\pdv{\vartheta}{\VEC{E}}] + 2\gamma\VEC{n}_2\otimes\qty[\pdv{\vartheta}{\VEC{E}}]).
\end{equation}

By using the definition of the tension field $\phi$ from \cref{eq:tensionfield}, the stress and material tensors are defined using \cref{eq:S_TFT}. Therefore, the assembly of the residual and Jacobian from \cref{eq:residual,eq:jacobian}, respectively, involves computing the tension field and inserting the corresponding option from \cref{eq:S_TFT} on each integration point.

\begin{remark}\label{remark:consistency}
	The modification scheme presented in this section, originally proposed by \fullcite{Nakashino2005}, provides an analytical derivation of the derivative $\Ccten^\prime$ of the wrinkling stress tensor $\Sten^\prime$. Consequently, the definitions in \cref{eq:S_TFT} are consistent and should provide optimal convergence in Newton--Raphson iterations. However, the definition of $\Sten^\prime$ in \cref{eq:S_TFT} depends on the tension field $\phi$ from \cref{eq:tensionfield}. The dependency of $\Sten^\prime$ on $\phi$ is not included in its derivative, and therefore the convergence behaviour of Newton--Raphson iterations can be suboptimal or diverging.
\end{remark}

%% file: Sections_Hyperelastic.tex
In this section, the novelty of this paper is presented. The theory from \fullcite{Nakashino2005}, recalled in the previous section, is extended to hyperelastic materials. The outline of this section is similar to that of \cref{sec:elastic}, but since the assumption for hyperelastic materials only affects the constitutive relation, the kinematic equations are as before, hence not included in this section. The present section primarily presents the differences with the elastic theory from \cref{sec:elastic}, first for the constitutive relation for the wrinkled membrane and then for the variational formulation.

\subsection{Constitutive Relation}
The derivation in \cref{sec:elastic} assumes a linear elastic constitutive model in \cref{eq:wrinklingstress}. In the case of non-linear hyperelastic material models, the constitutive relation is defined by a strain energy density function $\Psi(\Cten)$ (see \cref{eq:S}), where $\Cten=\Ften^\top\Ften=a_{\alpha\beta}\:\Avec^\alpha\otimes \Avec^\beta$ is the deformation tensor (see \cref{eq:deformationtensor}). In order to derive the wrinkling stress for hyperelastic materials, denoted by $\Sten^\prime = \Sten(\Eten^\prime)$, the stress tensor is linearized around the Green-Lagrange strain $\Eten$, i.e.
\begin{equation}\label{eq:wrinklingstress_HE}
	\begin{aligned}
		\Sten^\prime = \Sten(\Eten^\prime) = \Sten(\Eten+\Eten_W) = \Sten(\Eten) + \pdv{\Sten}{\Eten}:\Eten_W + \mathcal{O}(\Eten_W^2),
		&= \Sten(\Eten) + \Ccten(\Eten):\Eten_W + \mathcal{O}(\Eten_W^2)\\ &\approx \Sten(\Eten) + \Ccten(\Eten):\Eten_W.
	\end{aligned}
\end{equation}
In the second last equality, the definition of the material tensor $\Ccten(\Eten)$ from \cref{eq:C} is used. The definition in \cref{eq:wrinklingstress_HE} is similar to \cref{eq:wrinklingstress} when the Taylor expansion is truncated and the contribution of $\mathcal{O}(\Eten_W^2)$ is neglected under the assumption that $\Eten_W^2$ is small. Therefore, equivalent to \cref{eq:wrinklingstress}, the coefficients of the wrinkling stress tensor $\Sten^\prime$ are given by:
\begin{equation}\label{eq:wrinklingstress_HE2}
	S^{\prime\alpha\beta}(\Eten) = S^{\alpha\beta}(\Eten) + \gamma\Ccoef(\Eten) ^{\alpha\beta\gamma\delta}n_{\gamma}n_{\delta}.
\end{equation}
Following from the hyperelastic modified wrinkling strain tensor from \cref{eq:wrinklingstress_HE}, the hyperelastic counterpart of \cref{eq:uniaxialtension_cauchy3} is derived:
\begin{equation}\label{eq:uniaxialtension_cauchy3_HE}
	\begin{aligned}
		S^{\alpha\beta}(\Eten) n_\alpha n_\beta + \gamma\Ccten(\Eten)^{\alpha\beta\gamma\delta}n_\alpha n_\beta n_\gamma n_\delta = 0,\\
		S^{\alpha\beta}(\Eten) m_\alpha n_\beta + \gamma\Ccten(\Eten)^{\alpha\beta\gamma\delta}m_\alpha n_\beta n_\gamma n_\delta = 0.
	\end{aligned}
\end{equation}
From these equations, the derivation of the hyperelastic counterpart of \cref{eq:gamma} is straightforward:
\begin{equation}\label{eq:gamma_HE}
	\gamma(\Eten) = - \frac{S^{\alpha\beta}(\Eten) n_\alpha n_\beta}{ \Ccoef^{\alpha\beta\gamma\delta}(\Eten)n_\gamma n_\delta n_\alpha n_\beta },
\end{equation}
as well as the equation for $f(\vartheta)$:
\begin{equation}\label{eq:f_HE}
	f(\vartheta,\Eten) \equiv S^{\alpha\beta}(\Eten) m_\alpha n_\beta + \gamma(\Eten) \Ccoef^{\alpha\beta\gamma\delta}(\Eten) m_\alpha n_\beta n_\gamma n_\delta = 0.
\end{equation}
Like in the linear expression in \cref{eq:f}, the result in \cref{eq:f_HE} is also dependent on $\Eten$, which is fixed while finding the root. In this way, the same root-finding procedure can be applied as discussed for the linear elastic model.

\subsection{Variational Formulation}
For the hyperelastic model, the variation of the wrinkling stress tensor $\Sten^\prime$ is derived from \cref{eq:wrinklingstress_HE2}:
\begin{equation}\label{eq:dS_HE}
	\begin{aligned}
		\delta S^{\prime\alpha\beta}(\Eten) = \pdv{S^{\prime\alpha\beta}}{E_{\sigma\tau}} \delta E_{\sigma\tau} &= \Ccoef^{\prime\alpha\beta\sigma\tau}\delta E_{\sigma\tau}\\ 
		&= \qty(\pdv{S^{\alpha\beta}}{E_{\sigma\tau}} + \pdv{\gamma}{E_{\sigma\tau}}\Ccoef^{\alpha\beta\gamma\delta} n_\gamma n_\delta + \gamma\pdv{\Ccoef^{\alpha\beta\gamma\delta}}{E_{\sigma\tau}}n_\gamma n_\delta + \gamma\Ccoef^{\alpha\beta\gamma\delta} \pdv{(n_\gamma n_\delta)}{E_{\sigma\tau}})\delta E_{\sigma\tau}.
	\end{aligned}
\end{equation}
Compared to \cref{eq:dS}, the expression in \cref{eq:dS_HE} contains an extra term with the derivative of the material tensor $\Ccten$ with respect to the strains $\Eten$. This derivative can be written in terms of the deformation tensor $\Cten$ following \cref{eq:Eten}
\begin{equation}
	\pdv{\Ccoef^{\alpha\beta\gamma\delta}}{E_{\sigma\tau}} = 2\pdv{\Ccoef^{\alpha\beta\gamma\delta}}{C_{\sigma\tau}}.
\end{equation}
Considering \cref{eq:dgammadE,eq:dnndT}, the derivatives of $\gamma$ and $f$ with respect to $\Eten$, see \cref{eq:dgammadE,eq:dgammadT,eq:dfdE,eq:dFdT,eq:dTdE}, need to be re-defined due to the dependency of $\Sten$ and $\Ccten$ on $\Eten$. Starting with $\gamma$, the derivative with respect to the strain tensor $\Eten$ becomes:

\begin{equation}\label{eq:dgammadE_HE}
	\begin{aligned}
		\pdv{\gamma}{E^{\sigma\tau}} &= - \frac{(\Ccoef^{\alpha\beta\gamma\delta}n_\gamma n_\delta n_\alpha n_\beta)(\Ccoef^{\alpha\beta\sigma\tau}n_\alpha n_\beta ) - \qty(S^{\alpha\beta} n_\alpha n_\beta)\qty(\pdv{\Ccoef^{\alpha\beta\gamma\delta}}{E_{\sigma\tau}}n_\gamma n_\delta n_\alpha n_\beta)}{\qty(\Ccoef^{\alpha\beta\gamma\delta}n_\gamma n_\delta n_\alpha n_\beta)^2}\\ 
		&= - \frac{\Ccoef^{\alpha\beta\sigma\tau}n_\alpha n_\beta + \gamma\pdv{\Ccoef^{\alpha\beta\gamma\delta}}{E_{\sigma\tau}}n_\gamma n_\delta n_\alpha n_\beta}{\Ccoef^{\alpha\beta\gamma\delta}n_\gamma n_\delta n_\alpha n_\beta}.
	\end{aligned}
\end{equation}
Furthermore, the derivative of the root \cref{eq:f_HE} with respect to the strain becomes:
\begin{equation}\label{eq:dFdE_HE}
	\pdv{f}{E_{\sigma\tau}} = \qty(\Ccoef^{\alpha\beta\sigma\tau}m_\alpha n_\beta + \pdv{\gamma}{E_{\sigma\tau}}\Ccoef^{\alpha\beta\gamma\delta}m_\gamma n_\delta n_\alpha n_\beta + \gamma \pdv{\Ccoef^{\alpha\beta\sigma\tau}}{E_{\sigma\tau}}m_\gamma n_\delta n_\alpha n_\beta).
\end{equation}
Since the other derivatives with respect to $\vartheta$ (see \cref{eq:dgammadT,eq:dFdT}) do not change for hyperelastic materials, \cref{eq:dnndT,eq:dTdE} can be evaluated using \cref{eq:dgammadT,eq:dgammadE_HE,eq:dFdT,eq:dFdE_HE}, given the derivative of $\Ccten$ with respect to $\Eten$. Since $\Ccten$ depends on the strain energy density function $\Psi(\Cten)$ (see \cref{eq:C}), its derivative can be computed analytically. However, when static condensation to satisfy the plane-stress criterion is performed numerically, which could be the case for compressible materials \cite{Kiendl2015,Verhelst2021}, it is not straightforward to compute the derivative of $\Ccten$. Alternatively, the analytical expression for the statically condensed material tensor (see \cref{eq:condensation}) can become lengthy, making the analytical derivation of its derivative a tedious exercise. As an alternative to analytical derivation of the material tensor, it can therefore be opted to use finite differences or automatic differentiation to obtain the derivative of $\Ccten$ with respect to the strain. In \cref{example:NH_incompressible_dC}, the analytical derivative of $\Ccten$ with respect to $\Eten$ for a statically condensed incompressible Neo-Hookean material model is provided.

\begin{example}\label{example:NH_incompressible_dC}
	From \cref{example:NH_incompressible}, it follows that for an incompressible Neo-Hookean material model,
	\begin{equation}\label{eq:dNHC}
		\begin{aligned}
			\pdv{\Ccoef^{\alpha\beta\gamma\delta}}{E_{\sigma\tau}} &= 2\pdv{\Ccoef^{\alpha\beta\gamma\delta}}{C_{\sigma\tau}}\\
			&= \mu\pdv{(J_0^{-2})}{C_{\sigma\tau}} \qty(2a^{\alpha\beta}a^{\gamma\delta} + a^{\alpha\xi}a^{\beta\eta} + a^{\alpha\eta}a^{\beta\xi}) + \mu J_0^{-2}\pdv{}{C_{\sigma\tau}} \qty(2a^{\alpha\beta}a^{\gamma\delta} + a^{\alpha\xi}a^{\beta\eta} + a^{\alpha\eta}a^{\beta\xi}).
		\end{aligned}
	\end{equation}
	Since $J_0^2=\vert a_{\alpha\beta}\vert/\vert \mathring{a}_{\alpha\beta}\vert$ and since $\pdv{J_0}{\sigma\tau}=\frac{J}{2} a^{\sigma\tau}$ \cite{Holzapfel2000}, the derivative of $J_0^{-2}$ is
	\begin{equation}\label{eq:dJ0}
		\pdv{(J_0^{-2})}{C_{\sigma\tau}} = -J_0^{-2}a^{\sigma\tau}.
	\end{equation}
	Furthermore, the derivative of the contravariant metric tensor $a^{\alpha\beta}$ is given by
	\begin{equation}\label{eq:dgcon}
		\pdv{a^{\alpha\beta}}{C_{\sigma\tau}} = -\frac{1}{2}\qty(a^{\alpha\sigma}a^{\beta\tau} + a^{\alpha\tau}a^{\beta\sigma}).
	\end{equation}
	Using \cref{eq:dJ0,eq:dgcon}, the derivatives in \cref{eq:dNHC} can be evaluated, and the analytical expression for $\partial\Ccoef^{\alpha\beta\gamma\delta}/\partial E_{\sigma\tau}$ for the incompressible Neo-Hookean material is found.
\end{example}


\subsection{Implementation}\label{subsec:hyperelastic_implementation}
As for the linear elastic model, the modified wrinkling model for hyperelastic materials can be expressed in terms of linear algebra operations using the Voigt notation of the strain, stress, and material tensor as a basis (see \cref{eq:voigt,eq:Cvoigt}). First of all, the derivative of $\gamma$ with respect to the strain tensor $\Eten$ for hyperelastic materials becomes, following from \cref{eq:dgammadE_HE},
\begin{equation}\label{eq:dgammadE_HE_voigt}
	\qty[\pdv{\gamma}{\VEC{E}}] = -\frac{\VEC{C}\cdot\VEC{n}_1 + \gamma\VEC{n}_1^\top\cdot\qty[\pdv{\Ccten}{\Eten}]\cdot \VEC{n}_1}{\VEC{n}_1^\top\cdot\VEC{C}\cdot\VEC{n}_1}.
\end{equation}
Moreover, using \cref{eq:dFdE_HE}, the derivative of $f$ with respect to $\Eten$ for hyperelastic materials becomes:
\begin{equation}\label{eq:dfdE_HE_voigt}
	\qty[\pdv{f}{\VEC{E}}] = \VEC{C}\cdot\VEC{n}_2 + \qty(\VEC{n}_2^\top\cdot\VEC{C}\cdot\VEC{n}_1)\qty[\pdv{\gamma}{\VEC{E}}] + \gamma \VEC{n}_2^\top \cdot\qty[\pdv{\Ccten}{\Eten}]\cdot \VEC{n}_1.
\end{equation}
Both expressions in \cref{eq:dgammadE_HE_voigt,eq:dfdE_HE_voigt}, as well as the extra contribution of $\qty[\partial\Ccten/\partial\Eten]$ in \cref{eq:dS_HE}, contain the inner-product of $\partial\Ccten/\partial\Eten$ with $\VEC{n}_1$. Instead of storing $\qty[\partial\Ccten/\partial\Eten]$ and multiplying it in both expressions with $\VEC{n}_1$, the product $\qty[\partial\Ccten/\partial\Eten]\cdot\VEC{n}_1$ can also be stored as a matrix:
\begin{equation}\label{eq:dCdE_voigt}
	\qty[\pdv{\Ccten}{\Eten}\cdot \VEC{n}_1] =
	\begin{bmatrix}
		\pdv{\Ccten}{E_{11}}\cdot\VEC{n}_1 & \pdv{\Ccten}{E_{22}}\cdot\VEC{n}_1 & \pdv{\Ccten}{E_{12}}\cdot\VEC{n}_1\\
	\end{bmatrix}.
\end{equation}
Using \cref{eq:dgammadE_HE_voigt,eq:dfdE_HE_voigt,eq:dCdE_voigt} together with \cref{eq:dgammadT_voigt,eq:dfdT_voigt,eq:dTdE_voigt}, the modified material tensor for hyperelastic wrinkled materials employing Voigt notations, $\VEC{C}^\prime$, can be computed, based on \cref{eq:dS_HE}:
\begin{equation}\label{eq:C_HE_voigt}
	\VEC{C}^\prime = \VEC{C} \qty( \TEN{I} + \VEC{n}_1\otimes\qty[\pdv{\gamma}{\VEC{E}}] + \pdv{\gamma}{\vartheta}\VEC{n_1}\otimes\qty[\pdv{\vartheta}{\VEC{E}}] + 2\gamma\VEC{n}_2\otimes\qty[\pdv{\vartheta}{\VEC{E}}]) + \gamma\qty[\pdv{\Ccten}{\Eten}\cdot \VEC{n}_1].
\end{equation}

Similar to the linear elastic wrinkling modification scheme from \cref{sec:elastic}, the modified tensors $\VEC{S}^\prime$ and $\VEC{C}^\prime$ can be used in the definition in \cref{eq:S_TFT} depending on the tension field, evaluated per quadrature point. In addition, \cref{remark:consistency} regarding the variation of the stress tensor with respect to the tension field also applies for the hyperelastic model.

%% file: Sections_SolutionStrategies.tex
To obtain the solution $\uvec$ to the variational formulation \cref{eq:dW}, discretised by the residual from \cref{eq:residual}, different solution strategies can be adopted. In this section, a brief overview of the Newton--Raphson and the Dynamic Relaxation methods employed in the numerical examples is given.\\

\subsection{Newton-Raphson Method}
Firstly, the Newton--Raphson method solves the system of equations
\begin{equation}
	\MAT{K}(\uvec_k)\Delta\uvec_{k+1} = -\VEC{R}(\uvec_k),
\end{equation}
for $\Delta \uvec_{k+1}$ given $\uvec_{k}$ and updating $\uvec_{k+1}=\uvec_{k}+\Delta \uvec_{k+1}$ in iteration $k$. Here, $\MAT{K}$ is the Jacobian of the system from \cref{eq:jacobian} and $\VEC{R}$ is the residual vector from \cref{eq:residual}. Typically, the iterations are terminated if the update norm is $\Vert\Delta\uvec_k\Vert/\Vert\Delta\uvec_0\Vert<\epsilon_{\Delta\uvec}$ or if the relative residual norm is $\Vert\VEC{R}(\uvec_k)\Vert/\Vert\VEC{R}(\uvec_0)\Vert<\epsilon_{\VEC{R}}$. Although the Newton--Raphson iterations provide second-order convergence towards the solution $\uvec$, the convergence region is bounded, meaning that the method is guaranteed to converge only for an initial guess $\uvec_0$ sufficiently close to the final solution $\uvec$. Furthermore, if the Jacobian matrix $\MAT{K}$ is not exact, the speed of convergence can be decreased. In the case of the methods presented by \cite{Nakashino2005} and in this paper, the variation of the stress tensor with respect to the tension field is not included in the Jacobian (see \cref{remark:consistency}), possibly deteriorating the convergence behaviour.\\

\subsection{Dynamic Relaxation Method}
A commonly used alternative for Newton--Raphson iterations for solving problems involving wrinkling stabilities is the Dynamic Relaxation (DR) method \cite{Otter1960}. In this method, a dynamic system with artificial stiffness and damping is solved. The advantages of the dynamic relaxation method are that it is robust given a sufficiently small step size and that it only requires the residual vector $\VEC{R}$. However, its convergence is very slow.
The dynamic relaxation method is based on the solution of the structural dynamics equation:
\begin{equation}\label{eq:preliminaries_DR_equation}
	M\ddot{\uvec}(t)+C\dot{\uvec}(t)-\VEC{R}(\uvec(t))=\VEC{0},
\end{equation}
where $M$ is the mass matrix, $\ddot{\uvec}$ is the vector of discrete accelerations, $C$ is the damping matrix, $\dot{\uvec}$ is the vector of discrete velocities, and $\VEC{R}(\uvec)$ is the residual vector. Using central finite differences, the acceleration vector $\ddot{\uvec}$ can be expressed in terms of the velocity vector $\dot{\uvec}$ and a time step $\Delta t$:
\begin{equation}\label{eq:preliminaries_DR_velocity}
	\ddot{\uvec}_{t} = \frac{\dot{\uvec}_{t+\Delta t / 2} - \dot{\uvec}_{t-\Delta t / 2}}{\Delta t},
\end{equation}
where the notation $u_{t}=u(t)$ is adopted for the sake of clarity. A common assumption in dynamic relaxation methods is to define the damping matrix proportional to the mass matrix, i.e., $C = cM$ or by other appropriate scaling techniques \cite{Rodriguez2011,Rezaiee-Pajand2017}. Since these approaches involve an extra parameter in the solver, an alternative approach is to use the so-called \emph{kinetic damping} approach \cite{Cundall1976}. In this approach, the kinetic energy in the system is traced, and the nodal velocities $\dot{\uvec}$ are set to zero when a peak in kinetic energy is detected. The advantage of this method is that no parameter for damping is required and that it provides robustness \cite{Barnes1988,Barnes1999,Shugar1990}. Firstly, the kinetic energy in the system is defined by:
\begin{equation}
	E^K_{t} = \frac{1}{2}\dot{\uvec}^\top M\dot{\uvec}.
\end{equation}
Hence, a peak is detected if $E^K_{t}> E^K_{K,t+\Delta t}$ for $\Delta t > 0$. It is assumed that a peak occurs in the middle of the interval $[t-\Delta t,t]$, hence at $t-\Delta t/2$ if $E^K_{t-3\Delta t/2}<E^K_{t-\Delta t /2 }$ and $E^K_{t-\Delta t /2} > E^K_{t+\Delta t/2}$ \cite{Topping1994}. In that situation, the displacement vector $\uvec_{t+\Delta t}$ and the velocity vector $\dot{\uvec}_{t+\Delta t /2}$ are known. Using these solutions, the displacements at the peak can be computed by: 
\begin{equation}\label{eq:preliminaries_DR_displacement}
	\uvec_{t^\star} = \uvec_{t + \Delta t} - \frac{3}{2}\dot{\uvec}_{t+\Delta t/2} + \frac{\Delta t}{2}M^{-1}\VEC{R}(\uvec_{t}),
\end{equation}
Where the peak time is denoted by $t^\star$. Using the displacement vector $\uvec_{t^\star}$, the method is re-initiated using $\uvec_{t^\star}$.
Since the velocities are fully damped after a kinetic energy peak, they are set to zero upon re-initialization. Hence, to compute the next step after the restart at a peak on $\uvec_{t^\star}$, the velocity vector for $\uvec_{t^\star+\Delta t/2}$ becomes:
\begin{equation}\label{eq:preliminaries_DR_peak_velo}
	\dot{\uvec}_{\uvec_{t^\star+\Delta t / 2}} = \frac{\Delta t}{2}M^{-1}\VEC{R}(\uvec_{t^\star})
\end{equation}
Using \cref{eq:preliminaries_DR_peak_velo}, the displacement vector $\uvec_{t^\star+\Delta t}$ can be found using \cref{eq:preliminaries_DR_displacement}. This kinematic damping procedure is successfully applied in \cite{Barnes1988,Barnes1999,Taylor2014,Rezaiee-Pajand2017}, among others, showing the robustness of the method while eliminating the need to determine the damping coefficient $c$. \\

The equations above can be solved if and only if the mass matrix $M$ is invertible. However, having a full mass matrix drastically increases the computational costs of the method, since it involves solving a linear system. Therefore, alternatives include to select $M=\rho\:\text{diag}(K)$ with linear stiffness matrix $K$ \cite{Papadrakakis1981} to use a diagonal lumped mass matrix \cite{Joldes2010,Joldes2011,Barnes1988,Barnes1999,Topping1994,Alic2016}
or to use a column-sum of the stiffness matrix \cite{Underwood1983}.\\ 

In this paper, the dynamic relaxation method is equipped with the kinetic damping approach and a diagonal lumped mass matrix. The latter is assembled based on a unit-density and scaled with a parameter $\alpha$ as in, tuning the speed of convergence of the DR iterations. In the work of \cite{Lee2011}, a brief study on the tuning of the parameter $\alpha$ is provided. In brief, a too low value of $\alpha$ leads to divergence of the DR method, whereas a too high value leads to slow convergence. Typically, the DR iterations are terminated if the relative residual norm is below a tolerance $\Vert\VEC{R}(\uvec_k)\Vert/\Vert\VEC{R}(\uvec_0)\Vert<\epsilon_{\VEC{R}}$, or if the relative kinetic energy is below a certain tolerance $E^K_k/E^K_0<\epsilon_{E^K}$.\\


%% file: Sections_Benchmarks.tex
In this section, four benchmarks are presented for verification of the model presented in this paper. The benchmarks are selected from previous works on hyperelastic wrinkling simulations. Firstly, an uniaxial tension test is performed to verify the implementation of the model, based on the examples in the works of \cite{Kiendl2015} and \cite{Verhelst2021}. In this case, the full domain is in wrinkling condition, hence the tension field does not change during the iterations. Secondly, the inflation of a square membrane is modelled in \cref{subsec:benchmarks_pillow}, inspired by \cite{Diaby2006}. Thirdly, \cref{subsec:benchmarks_annulus} provides an example of a planar annular sheet in which the inner boundary is pulled out of the plane and twisted. This benchmark is inspired by the example given by \cite{Taylor2014} for linear elastic materials modified using a hyperelastic material model. Lastly, \cref{subsec:benchmarks_cylinder} models a cylindrical surface subject to large axial strain and a radial twist to demonstrate the capabilities of the present model on conic surfaces under large strains. For all benchmarks, a wrinkling simulation resolving wrinkling amplitudes is provided as a reference, along with results from the literature if available. The former are generated using isogeometric Kirchhoff-Love shells with hyperelastic constitutive models \cite{Kiendl2009,Kiendl2015,Verhelst2021}.\\

In the sequel, different hyperelastic material models are used. The compressible and incompressible Neo-Hookean (NH) and Mooney-Rivlin (MR) material models are given by the following strain energy density functions:
\begin{align}
	\Psi(\Cten) &= \frac{\mu}{2}\qty(J^{-\frac{2}{3}}I_1-3) + \Psi_{\text{vol}}(J)&& \text{NH Compressible}\label{eq:NH_comp},\\
	\Psi(\Cten) &= \frac{\mu}{2}\qty(I_1-3), && \text{NH Incompressible}\label{eq:NH_incomp},\\
	\Psi(\Cten) &= \frac{c_1}{2}\qty(J^{-\frac{2}{3}}I_1-3) + \frac{c_2}{2}\qty(J^{-\frac{4}{3}}I_2-3) + \Psi_{\text{vol}}(J), && \text{MR Compressible}\label{eq:MR_comp},\\
	\Psi(\Cten) &= \frac{c_1}{2}\qty(I_1-3) + \frac{c_2}{2}\qty(I_2-3), && \text{MR Incompressible} \label{eq:MR_incomp}.
\end{align}
Here, $\mu$ is the second Lamé parameter, defined as $\mu=E / (2 (1 + \nu))$ and $c_1$ and $c_2$ are the coefficients controlling the Mooney-Rivlin model via $\mu=c_1+c_2$. Furthermore, $\Psi_{\text{vol}}$ the volumetric strain energy density function using bulk modulus $K$ and parameter $\beta=-2$:
\begin{equation}
	\Psi_{\text{vol}} = K\mathcal{G}(J) = K\beta^{-2}\qty( \beta\log(J) + J^{-\beta}-1).
\end{equation}
For the incompressible Neo-Hookean material model, the analytical derivative of the material tensor is implemented, whereas for the other material models a finite-difference technique is used to obtain this term.

\subsection{Square subject to Tension}\label{subsec:benchmarks_UAT}
As a first example, a uniaxial tension test is performed on a square membrane, see \cref{fig:benchmarks_UAT_setup} for the model parameters inspired by \cite{Kiendl2015} and \cite{Verhelst2021}. The uniaxial tension test is performed both with an increasing load and with the value given in \cref{fig:benchmarks_UAT_setup}. The load stepping simulation provides a load displacement curve, which will be used against the implementation in our previous work \cite{Verhelst2021}, replicating analytical solutions. The simulation with the fixed value of the line load $p$ is used to assess the convergence in Newton--Raphson iterations. This is used to confirm the correct derivation and implementation of the modified stress and material tensors $\VEC{S}^\prime$ and $\Ccten^\prime$. The tests are performed for compressible and incompressible Neo-Hookean (NH) and Mooney-Rivlin (MR) material models; see \cref{eq:NH_comp,eq:NH_incomp,eq:MR_comp,eq:MR_incomp}.\\

\Cref{fig:benchmarks_UAT} presents the results of the uniaxial tension test. In the top plots, the stretch $\lambda$ is plotted against the applied line load $p$ for the compressible and incompressible Neo-Hookean and Mooney-Rivlin material models. In the bottom plots, the convergence is presented using the current and previous relative residual norms, $\Vert \VEC{R}_{i+1} \Vert / \Vert \VEC{R}_0 \Vert$ and $\Vert \VEC{R}_{i} \Vert / \Vert \VEC{R}_0 \Vert$, respectively. The load-displacement curves (top) show that the tension field theory modification scheme accurately predicts the constitutive behaviour of the membrane compared with the original model. Furthermore, the convergence plots (bottom) show that for compressible and incompressible Neo-Hookean and Mooney-Rivlin models, the convergence rate is optimal, i.e., quadratic, for Newton--Raphson iterations, but flattens out due to machine precision.

\begin{figure}
	\begin{minipage}{0.6\linewidth}
		\centering
		\resizebox{\linewidth}{!}
		{
			\begin{tikzpicture}[scale=0.8]
				\def\B{4}
				\def\L{4}
				\def\Bdef{2.828427}
				\def\Ldef{8}

				\filldraw[opacity=1, bottom color=black!10, top color=black!05, draw=black] (0,0) -- (\Ldef,0) -- (\Ldef,\Bdef) -- (0,\Bdef) --cycle;

				\draw[draw=black] (0,0) -- (\L,0) node[midway,below]{$\Gamma_1$} -- (\L,\B) node[midway,left]{$\Gamma_2$} -- (0,\B) node[midway,below]{$\Gamma_3$} --cycle node[midway,right]{$\Gamma_4$};

				\draw[latex-latex] (0,\B+0.1*\B) -- node[midway,above] {$L$} (\L,\B+0.1*\B);
				\draw[latex-latex] (-0.1*\B,0) -- node[midway,left] {$W$} (-0.1*\B,\B);

				\draw[-latex] (0.1,0.1) -- node[above right] {$x$}(1,0.1);
				\draw[-latex] (0.1,0.1) -- node[above right] {$y$}(0.1,1);

				\foreach \k in {0.5,1.0,...,9.5}
				{
					\pgfmathsetmacro\y{\k*\B / 10}
					\draw[-latex] (\L,\y) -- (\L + \L/10,\y);
				}
				\node[right] at (\L + \L/10,\B / 2) {$p$};

			\end{tikzpicture}
		}
	\end{minipage}
	\hspace{0.05\linewidth}
	\begin{minipage}{0.3\linewidth}
		\centering
		\footnotesize
		\begin{tabular}{llr}
			\toprule
			\multicolumn{3}{c}{Geometry}\\
			\midrule
			$L$, $W$ & 1.0 & $[\text{m}]$\\
			$t$ & 1.0 & $[\text{mm}]$\\
			\midrule
			\multicolumn{3}{c}{Material}\\
			\multicolumn{3}{c}{\textit{Neo-Hookean or Mooney-Rivlin}}\\
			\midrule
			$\mu$ & $1.5$ & $[\text{MPa}]$\\
			$\nu$ & 0.45 (comp.) & $[-]$\\
			& 0.50 (incomp.) & $[-]$\\
			$E$ & $2\mu(1+\nu)$ & $[\text{MPa}]$\\
			$c_1/c_2$ & $7.0$ & $[-]$\\
			\midrule
			\multicolumn{3}{c}{Boundary Conditions}\\
			\midrule
			$\Gamma_1$ & $u_y=0$ & \\
			$\Gamma_2$,$\Gamma_3$ & Free & \\
			$\Gamma_4$ & $u_x=0$ & \\
			\midrule
			\multicolumn{3}{c}{Loads}\\
			\midrule
			$p$ & $1.0$ & $[\text{MPa}]$\\
			\bottomrule
		\end{tabular}
	\end{minipage}
	\caption{Set-up for the uniaxial tension benchmark problem. In the figure on the left, the filled geometry represents the deformed configuration, and the dashed line indicates the undeformed geometry. The load $p$ indicates a line load acting on the undeformed geometry. The table on the right provides the parameter values for the specific benchmark problem for Neo-Hookean (NH) and Mooney-Rivlin (MR) materials.}
	\label{fig:benchmarks_UAT_setup}
\end{figure}

\begin{figure}[]
	\centering
	\begin{tikzpicture}
		\begin{groupplot}
			[
			height=0.2\textheight,
			width=0.38\linewidth,
			enlarge x limits = true,
			enlarge y limits = true,
			scale only axis,
			xmin = 0,
			xmax = 12,
			ymin = 0,
			ymax = 1.8e7,
			xlabel = {Stretch $\lambda$},
			ylabel = {Applied Load $p$},
			group style={
				group name = top,
				group size=2 by 1,
				ylabels at=edge left,
				y descriptions at=edge left,
				xlabels at=edge bottom,
				x descriptions at=edge bottom,
				vertical sep=0.05\textheight,
				horizontal sep=0.04\linewidth,
			},
			legend pos = north west,
			]
			\nextgroupplot[title={Incompressible}]
			\addlegendimage{style=style1,only marks,opacity=0.5}\addlegendentry{NH};
			\addlegendimage{style=style1,no markers}\addlegendentry{NH - TFT};
			\addlegendimage{style=style2,only marks,opacity=0.5}\addlegendentry{MR};
			\addlegendimage{style=style2,no markers}\addlegendentry{MR - TFT};

			\addplot+[style=style1,mark repeat=5,only marks,opacity=0.5] table[col sep=comma,header=true,x expr=\thisrowno{0}/1+1,y expr=\thisrowno{3}/1e-3/1] {Results_UAT_NH_Incomp_data.csv};
			\addplot+[style=style1,no markers] table[col sep=comma,header=true,x expr=\thisrowno{0}/1+1,y expr=\thisrowno{3}/1e-3/1] {Results_UAT_NH_Incomp_TFT_data.csv};
			\addplot+[style=style2,mark repeat=5,only marks,opacity=0.5] table[col sep=comma,header=true,x expr=\thisrowno{0}/1+1,y expr=\thisrowno{3}/1e-3/1] {Results_UAT_MR_Incomp_data.csv};
			\addplot+[style=style2,no markers] table[col sep=comma,header=true,x expr=\thisrowno{0}/1+1,y expr=\thisrowno{3}/1e-3/1] {Results_UAT_MR_Incomp_TFT_data.csv};

			\nextgroupplot[title={Compressible}]
			\addplot+[style=style1,mark repeat=5,only marks,opacity=0.5] table[col sep=comma,header=true,x expr=\thisrowno{0}/1+1,y expr=\thisrowno{3}/1e-3/1] {Results_UAT_NH_Comp_data.csv};
			\addplot+[style=style1,no markers] table[col sep=comma,header=true,x expr=\thisrowno{0}/1+1,y expr=\thisrowno{3}/1e-3/1] {Results_UAT_NH_Comp_TFT_data.csv};
			\addplot+[style=style2,mark repeat=5,only marks,opacity=0.5] table[col sep=comma,header=true,x expr=\thisrowno{0}/1+1,y expr=\thisrowno{3}/1e-3/1] {Results_UAT_MR_Comp_data.csv};
			\addplot+[style=style2,no markers] table[col sep=comma,header=true,x expr=\thisrowno{0}/1+1,y expr=\thisrowno{3}/1e-3/1] {Results_UAT_MR_Comp_TFT_data.csv };

		\end{groupplot}

		\begin{groupplot}
			[
			height=0.1\textheight,
			width=0.17\linewidth,
			enlarge x limits = true,
			enlarge y limits = true,
			scale only axis,
			xmode = log,
			ymode = log,
			xmin = 1e-5,
			xmax = 1e0,
			ymin = 1e-7,
			ymax = 1e0,
			xlabel = {$\Vert R_{i+1} \Vert / \Vert R_0 \Vert$},
			ylabel = {$\Vert R_{i} \Vert / \Vert R_0 \Vert$},
			legend pos = north west,
			group style={
				group size=4 by 1,
				ylabels at=edge left,
				y descriptions at=edge left,
				xlabels at=edge bottom,
				x descriptions at=edge bottom,
				vertical sep=0.05\textheight,
				horizontal sep=0.04\linewidth,
			},
			]

			\nextgroupplot[title={NH},anchor=north west, at={($(top c1r1.south west) + (0,-0.075\textheight)$)}]
			\addplot+[style=style1,only marks,opacity=0.5] table[col sep=comma,header=true,x expr=exp(\thisrowno{0}),y expr=exp(\thisrowno{1})] {Results_UAT_IncompressibleNH3.csv};
			\addplot+[style=style1,no markers] table[col sep=comma,header=true,x expr=exp(\thisrowno{0}),y expr=exp(\thisrowno{1})] {Results_UAT_IncompressibleNH3_TFT.csv};
			\logLogSlopeTriangleRev{0.8}{0.2}{0.835}{2}{black}
			\nextgroupplot[title={MR}]
			\addplot+[style=style2,only marks,opacity=0.5] table[col sep=comma,header=true,x expr=exp(\thisrowno{0}),y expr=exp(\thisrowno{1})] {Results_UAT_IncompressibleMR3.csv};
			\addplot+[style=style2,no markers] table[col sep=comma,header=true,x expr=exp(\thisrowno{0}),y expr=exp(\thisrowno{1})] {Results_UAT_IncompressibleMR3_TFT.csv};
			\logLogSlopeTriangleRev{0.8}{0.2}{0.835}{2}{black}

			\nextgroupplot[title={NH}]
			\addplot+[style=style1,only marks,opacity=0.5] table[col sep=comma,header=true,x expr=exp(\thisrowno{0}),y expr=exp(\thisrowno{1})] {Results_UAT_CompressibleNH3.csv};
			\addplot+[style=style1,no markers] table[col sep=comma,header=true,x expr=exp(\thisrowno{0}),y expr=exp(\thisrowno{1})] {Results_UAT_CompressibleNH3_TFT.csv};
			\logLogSlopeTriangleRev{0.8}{0.2}{0.835}{2}{black}

			\nextgroupplot[title={MR}]
			\addplot+[style=style2,only marks,opacity=0.5] table[col sep=comma,header=true,x expr=exp(\thisrowno{0}),y expr=exp(\thisrowno{1})] {Results_UAT_CompressibleMR3.csv};
			\addplot+[style=style2,no markers] table[col sep=comma,header=true,x expr=exp(\thisrowno{0}),y expr=exp(\thisrowno{1})] {Results_UAT_CompressibleMR3_TFT.csv};
			\logLogSlopeTriangleRev{0.8}{0.2}{0.835}{2}{black}
		\end{groupplot}
	\end{tikzpicture}
	\caption{Results for the uniaxial tension benchmark problem from \cref{fig:benchmarks_UAT_setup}. The top plots depict the stretch $\lambda$ versus the applied load $p$. The bottom plots present the convergence of the relative residual norm $\Vert \VEC{R}_i \Vert / \Vert \VEC{R}_0 \Vert$ for a load step with $p=1.0\:[\text{MPa}]$, with the triangle indicating second-order convergence. The results are depicted for Neo-Hookean (NH) and Mooney-Rivlin (MR) materials with the parameters from \cref{fig:benchmarks_UAT_setup}. The lines indicate the results without a tension field theory (TFT) modification, and the markers indicate the results including the modification proposed in this paper. }
	\label{fig:benchmarks_UAT}
\end{figure}

\subsection{Square subject to Inflation}\label{subsec:benchmarks_pillow}
In the next example, the inflation of a square membrane is modelled (see \cref{fig:benchmarks_pillow_setup}). This example is inspired by the works of \cite{Diaby2006}, among others. In most previous works, the inflated square membrane was modelled using linear elastic models \cite{Nakashino2020,Jarasjarungkiat2009,Kang1999,Contri1988a,Ziegler2003}, but in the work of \cite{Diaby2006} the case was used with a hyperelastic Neo-Hookean material model, which is adopted in this paper as well. Furthermore, the inflation is modelled using a surface loading $\VEC{f}(\uvec)=p\hat{\VEC{n}}$, see \cref{eq:dWext}, with the pressure $p=5000\:[\text{Pa}]$. The simulation is solved in two stages, inspired by \cite{Diaby2006}, to ensure bounded equilibrium iterations for the membrane. First, an in-plane boundary load of $p=5000\:[\text{Pa}]$ orthogonal to the boundary is applied to $\Gamma_1$ and $\Gamma_4$ to pre-stretch the membrane. Thereafter, the pressure is applied to the pre-stretched domain, and the boundary load is removed. The equilibrium iterations are performed using a Dynamic Relaxation method with $\epsilon_{\VEC{R}} = 10^{-6}$. In addition, the solutions are computed on meshes with $8\times8$, $16\times16$ and $32\times32$ elements. For the membrane model, elements with degree $p\in\{1,2,3\}$ are used, whereas the Kirchhoff--Love shell model uses $p\in\{2,3\}$ due to the $C^1$ continuity requirement of the shell model.\\

\begin{figure}
	\centering
	\begin{minipage}{0.6\linewidth}
		\centering
		\resizebox{\linewidth}{!}
		{
			\tdplotsetmaincoords{60}{-40}
			\begin{tikzpicture}[
				scale=1.5,
				tdplot_main_coords,bullet/.style={circle,inner
					sep=1pt,fill=black,fill opacity=1}]
				\begin{scope}[canvas is xy plane at z=0]
					\fill[top color=black!40,bottom color=black!20](0,0) rectangle (5,5);
					\fill[top color=black!60,bottom color=black!40](0,0) rectangle (2.5,2.5);
					\draw[thick](2.5,0)--(2.5,2.5) node[midway,above right,inner sep=4pt]{$\Gamma_1$};
					\draw[thick](0,0)--(2.5,0) node[midway,below right,inner sep=4pt]{$\Gamma_2$};
					\draw[thick](0,0)--(0,2.5) node[midway,below left,inner sep=4pt]{$\Gamma_3$};
					\draw[thick](0,2.5)--(2.5,2.5) node[midway,above left,inner sep=4pt]{$\Gamma_4$};
					\draw[latex-latex] (0,5.2) -- (5,5.2) node[midway,above left]{$L$};
					\draw[latex-latex] (5.2,0) -- (5.2,5) node[midway,above right]{$L$};
					\filldraw[draw=black,fill=white] (2.5,2.5) circle[radius=0.05];
					\draw[latex-] ($(2.5,2.5)+(0.2,0.1)$) to[in=80,out=20] (3,2.5) node[right]{$M$};
					\filldraw[draw=black,fill=white] (0.0,0.0) circle[radius=0.05];
					\draw[latex-] ($(0.0,0.0)+(-0.0,-0.1)$) to[in=180,out=-110] (0.0,-0.5) node[right]{$A$};
				\end{scope}
				\foreach \x in {0.2,0.4,...,2.3}
				{
					\foreach \y in {0.2,0.4,...,2.3}
					{
						\draw[-latex] (\x,\y,0)--(\x,\y,0.2);
					}
				}
				\node[above] at (2.3,0.0,0.2){$p$};
				\draw[-latex](3.0,3.0,0)--(3.5,3.0,0) node[right]{$x$};
				\draw[-latex](3.0,3.0,0)--(3.0,3.5,0) node[left]{$y$};
				\draw[-latex](3.0,3.0,0)--(3.0,3.0,0.5) node[above]{$z$};
			\end{tikzpicture}
		}
	\end{minipage}
	\hfill
	\begin{minipage}{0.3\linewidth}
		\centering
		\footnotesize
		\begin{tabular}{llr}
			\toprule
			\multicolumn{3}{c}{Geometry}\\
			\midrule
			$L$ & $\sqrt{1.2^2/2}$ & $[\text{m}]$\\
			$t$ & 0.1 & $[\text{mm}]$\\
			\midrule
			\multicolumn{3}{c}{Material}\\
			\multicolumn{3}{c}{\textit{Compressible Neo-Hookean}}\\
			\midrule
			$E$ & 588 & $[\text{MPa}]$\\
			$\nu$ & 0.4 & $[-]$\\
			\midrule
			\multicolumn{3}{c}{Boundary Conditions}\\
			\midrule
			$\Gamma_1$ & $u_x=0$ & \\
			$\Gamma_2$,$\Gamma_3$ & $u_z=0$ & \\
			$\Gamma_4$ & $u_y=0$ & \\
			\midrule
			\multicolumn{3}{c}{Loads}\\
			\midrule
			$p$ & $5000$ & $[\text{Pa}]$\\
			\bottomrule
		\end{tabular}
	\end{minipage}
	\caption{Problem definition for a square membrane with diagonal length $1.2\:[\text{m}]$ subject to pressure $p$. The membrane is modelled using in-plane symmetry boundary conditions on $\Gamma_1$ and $\Gamma_4$. Furthermore, the sides $\Gamma_2$ and $\Gamma_3$ have restricted $z$-displacement. The square membrane has a Neo-Hookean material model with the parameters provided in the table on the right.}
	\label{fig:benchmarks_pillow_setup}
\end{figure}

\begin{figure}
	\centering
	\begin{subfigure}[t]{0.45\linewidth}
		\centering
		\includegraphics[width=0.9\linewidth]{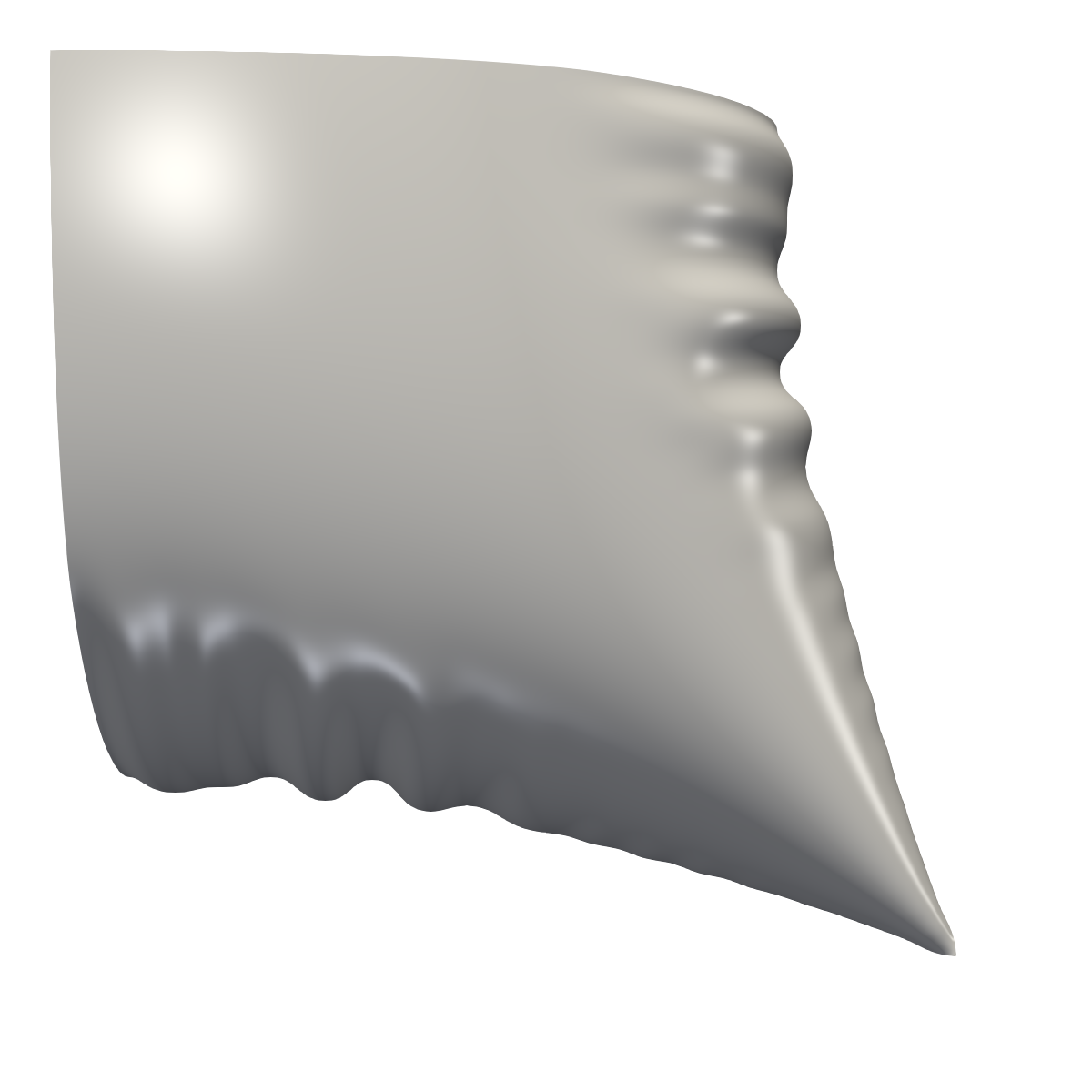}
		\caption{Top view of a quarter of the inflated square membrane for $32\times32$ elements with degree $p=3$ using the KL shell.}
		\label{fig:benchmarks_pillow_deformed}
	\end{subfigure}
	\hfill
	\begin{subfigure}[t]{0.45\linewidth}
		\centering
		\includegraphics[width=0.9\linewidth]{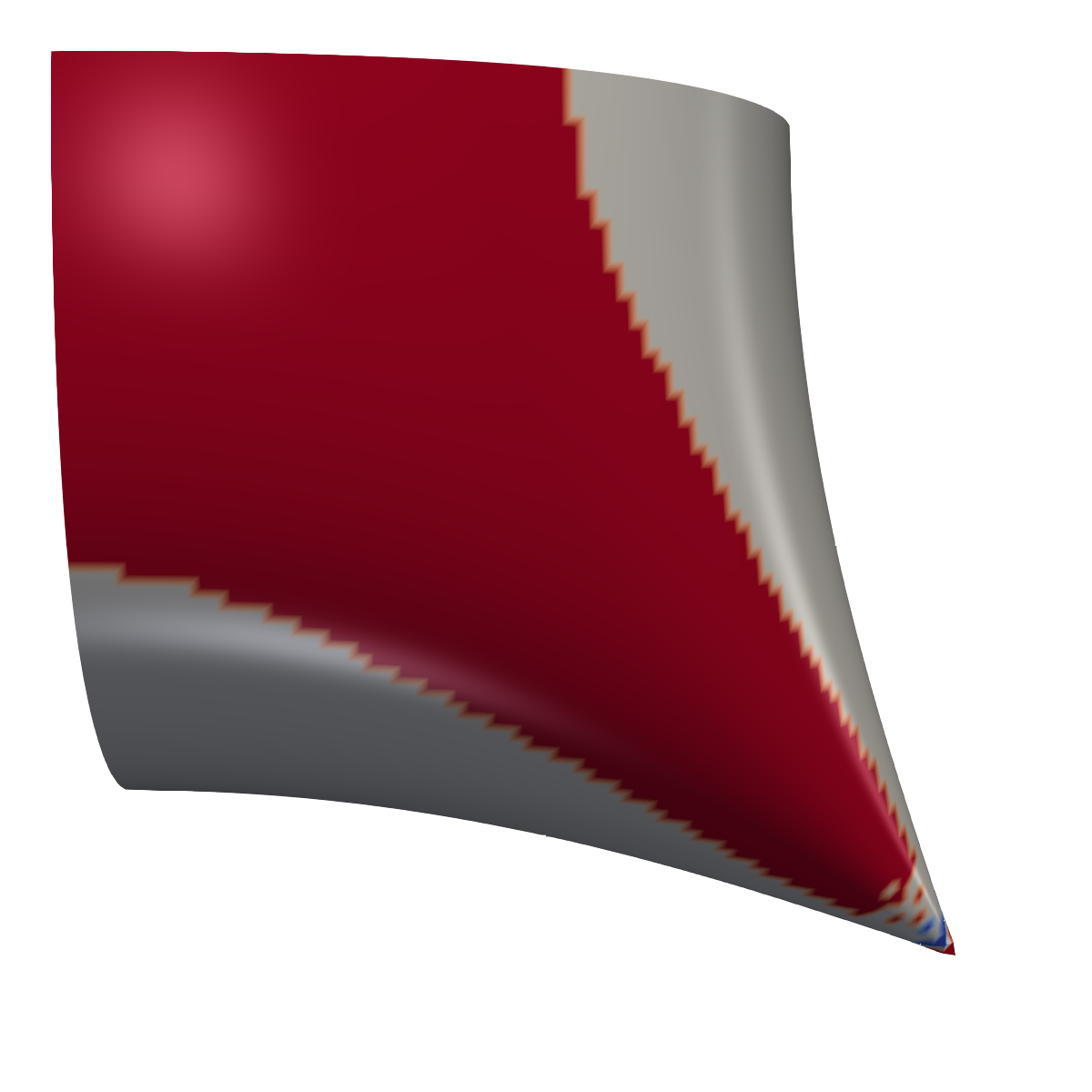}
		\caption{Top view of a quarter of the inflated square membrane for $32\times32$ elements with degree $p=3$ using the membrane model with tension field theory modifications. The red colour represents a taut region, the grey colour represents a wrinkled region and the blue color (bottom right) represents a slack region.}
		\label{fig:benchmarks_pillow_tensionfield}
	\end{subfigure}


	\centering
	\begin{subfigure}{\linewidth}
		\begin{tikzpicture}
			\begin{groupplot}
				[
				height=0.3\linewidth,
				width=0.3\linewidth,
				enlarge x limits=false,
				enlarge y limits=false,
				scale only axis,
				hide axis,
				group style={
					group size=3 by 1,
					ylabels at=edge left,
					y descriptions at=edge left,
					xlabels at=edge bottom,
					x descriptions at=edge bottom,
					vertical sep=0.05\textheight,
					horizontal sep=0.05\linewidth,
				},
				no markers,
				]
				\nextgroupplot[title={$p=1$},legend pos = north west]
				\addlegendimage{style=style2,solid}\addlegendentry{$16\times16$}
				\addlegendimage{style=style3,solid}\addlegendentry{$32\times32$}
				\addlegendimage{style=style0,solid,black}\addlegendentry{KL Shell}
				\addlegendimage{style=style0,dashed,black}\addlegendentry{TFT Membrane}
				\addplot+[gray,opacity=0.25,solid] table[col sep=comma,header=true,x index={0},y index={1}] {Pillow_r4_e1-M1-c1_line_u=0.000000_dir=1.csv};
				\addplot+[gray,opacity=0.25,solid] table[col sep=comma,header=true,x index={0},y index={1}] {Pillow_r4_e2-M1-c1_line_u=0.000000_dir=1.csv};
				\addplot+[gray,opacity=0.25,solid] table[col sep=comma,header=true,x index={0},y index={1}] {Pillow_r5_e1-M1-c1_line_u=0.000000_dir=1.csv};
				\addplot+[gray,opacity=0.25,solid] table[col sep=comma,header=true,x index={0},y index={1}] {Pillow_r5_e2-M1-c1_line_u=0.000000_dir=1.csv};

				\addplot+[gray,opacity=0.25,solid] table[col sep=comma,header=true,x index={0},y index={1}] {Pillow_r4_e1-M1-c1_line_u=1.000000_dir=0.csv};
				\addplot+[gray,opacity=0.25,solid] table[col sep=comma,header=true,x index={0},y index={1}] {Pillow_r4_e2-M1-c1_line_u=1.000000_dir=0.csv};
				\addplot+[gray,opacity=0.25,solid] table[col sep=comma,header=true,x index={0},y index={1}] {Pillow_r5_e1-M1-c1_line_u=1.000000_dir=0.csv};
				\addplot+[gray,opacity=0.25,solid] table[col sep=comma,header=true,x index={0},y index={1}] {Pillow_r5_e2-M1-c1_line_u=1.000000_dir=0.csv};

				\addplot+[gray,opacity=0.25,dashed] table[col sep=comma,header=true,x index={0},y index={1}] {Pillow_r4_e1-M1-c1_TFT_line_u=0.000000_dir=1.csv};
				\addplot+[gray,opacity=0.25,dashed] table[col sep=comma,header=true,x index={0},y index={1}] {Pillow_r4_e2-M1-c1_TFT_line_u=0.000000_dir=1.csv};
				\addplot+[gray,opacity=0.25,dashed] table[col sep=comma,header=true,x index={0},y index={1}] {Pillow_r5_e0-M1-c1_TFT_line_u=0.000000_dir=1.csv};
				\addplot+[gray,opacity=0.25,dashed] table[col sep=comma,header=true,x index={0},y index={1}] {Pillow_r5_e1-M1-c1_TFT_line_u=0.000000_dir=1.csv};
				\addplot+[gray,opacity=0.25,dashed] table[col sep=comma,header=true,x index={0},y index={1}] {Pillow_r5_e2-M1-c1_TFT_line_u=0.000000_dir=1.csv};

				\addplot+[gray,opacity=0.25,dashed] table[col sep=comma,header=true,x index={0},y index={1}] {Pillow_r4_e0-M1-c1_TFT_line_u=1.000000_dir=0.csv};
				\addplot+[gray,opacity=0.25,dashed] table[col sep=comma,header=true,x index={0},y index={1}] {Pillow_r4_e2-M1-c1_TFT_line_u=1.000000_dir=0.csv};
				\addplot+[gray,opacity=0.25,dashed] table[col sep=comma,header=true,x index={0},y index={1}] {Pillow_r5_e0-M1-c1_TFT_line_u=1.000000_dir=0.csv};
				\addplot+[gray,opacity=0.25,dashed] table[col sep=comma,header=true,x index={0},y index={1}] {Pillow_r5_e1-M1-c1_TFT_line_u=1.000000_dir=0.csv};
				\addplot+[gray,opacity=0.25,dashed] table[col sep=comma,header=true,x index={0},y index={1}] {Pillow_r5_e2-M1-c1_TFT_line_u=1.000000_dir=0.csv};

				\addplot+[style=style2,dashed] table[col sep=comma,header=true,x index={0},y index={1}] {Pillow_r4_e0-M1-c1_TFT_line_u=0.000000_dir=1.csv};
				\addplot+[style=style3,dashed] table[col sep=comma,header=true,x index={0},y index={1}] {Pillow_r4_e1-M1-c1_TFT_line_u=1.000000_dir=0.csv};

				\nextgroupplot[title={$p=2$},legend pos = north west]
				\addplot+[gray,opacity=0.25,solid] table[col sep=comma,header=true,x index={0},y index={1}] {Pillow_r4_e2-M1-c1_line_u=0.000000_dir=1.csv};
				\addplot+[gray,opacity=0.25,solid] table[col sep=comma,header=true,x index={0},y index={1}] {Pillow_r5_e1-M1-c1_line_u=0.000000_dir=1.csv};
				\addplot+[gray,opacity=0.25,solid] table[col sep=comma,header=true,x index={0},y index={1}] {Pillow_r5_e2-M1-c1_line_u=0.000000_dir=1.csv};

				\addplot+[gray,opacity=0.25,solid] table[col sep=comma,header=true,x index={0},y index={1}] {Pillow_r4_e1-M1-c1_line_u=1.000000_dir=0.csv};
				\addplot+[gray,opacity=0.25,solid] table[col sep=comma,header=true,x index={0},y index={1}] {Pillow_r4_e2-M1-c1_line_u=1.000000_dir=0.csv};
				\addplot+[gray,opacity=0.25,solid] table[col sep=comma,header=true,x index={0},y index={1}] {Pillow_r5_e2-M1-c1_line_u=1.000000_dir=0.csv};

				\addplot+[gray,opacity=0.25,dashed] table[col sep=comma,header=true,x index={0},y index={1}] {Pillow_r4_e0-M1-c1_TFT_line_u=0.000000_dir=1.csv};
				\addplot+[gray,opacity=0.25,dashed] table[col sep=comma,header=true,x index={0},y index={1}] {Pillow_r4_e2-M1-c1_TFT_line_u=0.000000_dir=1.csv};
				\addplot+[gray,opacity=0.25,dashed] table[col sep=comma,header=true,x index={0},y index={1}] {Pillow_r5_e0-M1-c1_TFT_line_u=0.000000_dir=1.csv};
				\addplot+[gray,opacity=0.25,dashed] table[col sep=comma,header=true,x index={0},y index={1}] {Pillow_r5_e1-M1-c1_TFT_line_u=0.000000_dir=1.csv};
				\addplot+[gray,opacity=0.25,dashed] table[col sep=comma,header=true,x index={0},y index={1}] {Pillow_r5_e2-M1-c1_TFT_line_u=0.000000_dir=1.csv};

				\addplot+[gray,opacity=0.25,dashed] table[col sep=comma,header=true,x index={0},y index={1}] {Pillow_r4_e0-M1-c1_TFT_line_u=1.000000_dir=0.csv};
				\addplot+[gray,opacity=0.25,dashed] table[col sep=comma,header=true,x index={0},y index={1}] {Pillow_r4_e1-M1-c1_TFT_line_u=1.000000_dir=0.csv};
				\addplot+[gray,opacity=0.25,dashed] table[col sep=comma,header=true,x index={0},y index={1}] {Pillow_r5_e0-M1-c1_TFT_line_u=1.000000_dir=0.csv};
				\addplot+[gray,opacity=0.25,dashed] table[col sep=comma,header=true,x index={0},y index={1}] {Pillow_r5_e1-M1-c1_TFT_line_u=1.000000_dir=0.csv};
				\addplot+[gray,opacity=0.25,dashed] table[col sep=comma,header=true,x index={0},y index={1}] {Pillow_r5_e2-M1-c1_TFT_line_u=1.000000_dir=0.csv};

				\addplot+[style=style2,solid] table[col sep=comma,header=true,x index={0},y index={1}] {Pillow_r4_e1-M1-c1_line_u=0.000000_dir=1.csv};
				\addplot+[style=style3,solid] table[col sep=comma,header=true,x index={0},y index={1}] {Pillow_r4_e2-M1-c1_line_u=1.000000_dir=0.csv};
				\addplot+[style=style2,dashed] table[col sep=comma,header=true,x index={0},y index={1}] {Pillow_r4_e1-M1-c1_TFT_line_u=0.000000_dir=1.csv};
				\addplot+[style=style3,dashed] table[col sep=comma,header=true,x index={0},y index={1}] {Pillow_r4_e2-M1-c1_TFT_line_u=1.000000_dir=0.csv};

				\nextgroupplot[title={$p=3$},legend pos = north west]
				\addplot+[gray,opacity=0.25,solid] table[col sep=comma,header=true,x index={0},y index={1}] {Pillow_r4_e1-M1-c1_line_u=0.000000_dir=1.csv};
				\addplot+[gray,opacity=0.25,solid] table[col sep=comma,header=true,x index={0},y index={1}] {Pillow_r4_e2-M1-c1_line_u=0.000000_dir=1.csv};
				\addplot+[gray,opacity=0.25,solid] table[col sep=comma,header=true,x index={0},y index={1}] {Pillow_r5_e2-M1-c1_line_u=0.000000_dir=1.csv};

				\addplot+[gray,opacity=0.25,solid] table[col sep=comma,header=true,x index={0},y index={1}] {Pillow_r4_e1-M1-c1_line_u=1.000000_dir=0.csv};
				\addplot+[gray,opacity=0.25,solid] table[col sep=comma,header=true,x index={0},y index={1}] {Pillow_r4_e2-M1-c1_line_u=1.000000_dir=0.csv};
				\addplot+[gray,opacity=0.25,solid] table[col sep=comma,header=true,x index={0},y index={1}] {Pillow_r5_e1-M1-c1_line_u=1.000000_dir=0.csv};

				\addplot+[gray,opacity=0.25,dashed] table[col sep=comma,header=true,x index={0},y index={1}] {Pillow_r4_e1-M1-c1_TFT_line_u=0.000000_dir=1.csv};
				\addplot+[gray,opacity=0.25,dashed] table[col sep=comma,header=true,x index={0},y index={1}] {Pillow_r4_e2-M1-c1_TFT_line_u=0.000000_dir=1.csv};
				\addplot+[gray,opacity=0.25,dashed] table[col sep=comma,header=true,x index={0},y index={1}] {Pillow_r5_e1-M1-c1_TFT_line_u=0.000000_dir=1.csv};

				\addplot+[gray,opacity=0.25,dashed] table[col sep=comma,header=true,x index={0},y index={1}] {Pillow_r4_e1-M1-c1_TFT_line_u=1.000000_dir=0.csv};
				\addplot+[gray,opacity=0.25,dashed] table[col sep=comma,header=true,x index={0},y index={1}] {Pillow_r4_e2-M1-c1_TFT_line_u=1.000000_dir=0.csv};
				\addplot+[gray,opacity=0.25,dashed] table[col sep=comma,header=true,x index={0},y index={1}] {Pillow_r5_e1-M1-c1_TFT_line_u=1.000000_dir=0.csv};

				\addplot+[style=style2,solid] table[col sep=comma,header=true,x index={0},y index={1}] {Pillow_r5_e1-M1-c1_line_u=0.000000_dir=1.csv};
				\addplot+[style=style3,solid] table[col sep=comma,header=true,x index={0},y index={1}] {Pillow_r5_e2-M1-c1_line_u=1.000000_dir=0.csv};
				\addplot+[style=style2,dashed] table[col sep=comma,header=true,x index={0},y index={1}] {Pillow_r5_e2-M1-c1_TFT_line_u=0.000000_dir=1.csv};
				\addplot+[style=style3,dashed] table[col sep=comma,header=true,x index={0},y index={1}] {Pillow_r5_e2-M1-c1_TFT_line_u=1.000000_dir=0.csv};
			\end{groupplot}
		\end{tikzpicture}
		\caption{Contour lines at $z=0$ for the inflated square membrane for a mesh with degree $p=1$ (left), $p=2$ (middle) and $p=3$ (right). The results of the shell simulation are depicted as solid lines, and the membrane tension field theory results are represented by dashed lines. The coloured lines highlight a certain number of elements, using symmetry for the horizontal and vertical directions, for the purpose of compactness. The grey lines indicate the other non-highlighted lines.}
		\label{fig:benchmarks_pillow_contours}
	\end{subfigure}
	\caption{Results for the square membrane subject to a pressure load from \cref{fig:benchmarks_pillow_setup}. (\subref{fig:benchmarks_pillow_deformed}) represents the deformed shape from a Kirchhoff--Love shell simulation, providing wrinkles. (\subref{fig:benchmarks_pillow_tensionfield}) provides the deformed shape from a membrane simulation with the proposed tension field theory modification scheme, together with the tension field. (\subref{fig:benchmarks_pillow_contours}) provides the contours of the deformation for both models for different numbers of elements and degrees.}
	\label{fig:benchmarks_pillow}
\end{figure}

\Cref{fig:benchmarks_pillow} depicts the results of the inflated square membrane. As can be seen in \cref{fig:benchmarks_pillow_deformed}, the membrane inflates into a pillow-shape with wrinkles along the free boundaries. Indeed, \cref{fig:benchmarks_pillow_tensionfield} shows that the tension field in the boundary regions indicates wrinkling, whereas the other parts are in a taut state, which is in line with similar observations as in \cite{Nakashino2020}. The contour lines on the $z=0$-plane, see \cref{fig:benchmarks_pillow_contours}, show large differences in the wrinkling pattern, varying the number of elements and the degree of the basis for the Kirchhoff--Love shell model. However, the membrane model using the tension field theory modification for hyperelastic materials as presented accurately captures the wrinkling mid-plane and is consistent across mesh refinement and degree elevation. The numerical results in \cref{tab:benchmarks_pillow} show that the mid-point and corner-point displacements of the membrane rapidly converge for the tension field membrane simulation, compared to the shell simulation. The final values of both simulations show small differences between the results obtained by the KL shell model and the membrane TFT model and slightly bigger differences with the results obtained by the shell model of \cite{Diaby2006}, possibly because of the lower degree finite element method and the lower number of elements used there.\\

\begin{table}
	\centering
	\caption{Results for the inflated square membrane; see \cref{fig:benchmarks_pillow_setup} for the problem parameters; for the Kirchhoff--Love (KL) shell resolving actual wrinkles; and for the membrane model using the tension field theory (TFT) modifications proposed in this paper. The reference results are provided as the vertical displacement in the mid-point $M$, $u_{z,M}$, and as the horizontal displacement along the $x$-axis in the point $A$, $u_{x,A}$. The results are compared with the work of \cite{Diaby2006}.}
	\label{tab:benchmarks_pillow}
	\footnotesize
	\begin{tabular}{llllllllll}
		\toprule
		&       & \multicolumn{3}{c}{KLShell}   & \multicolumn{3}{c}{Membrane (TFT)} & \multicolumn{2}{c}{\cite{Diaby2006}}\\
		\cmidrule(lr){3-5}\cmidrule(lr){6-8}\cmidrule(lr){9-10}
		\multicolumn{2}{l}{\# elements}                 & $8\times8$            & $16\times16$           & $32\times32$          & $8\times8$            & $16\times16$           & $32\times32$      & $5\times5$        & $25\times25$\\
		\midrule
		\multirow{3}{*}{$u_{z,M}$}  & $p=1$ &   &   &   & 0.2153	    & 0.2181	    & 0.2188    & \multirow{2}{*}{0.2144}       & \multirow{2}{*}{0.2245}   \\
		& $p=2$ & 0.2094	    & 0.2161	    & 0.2182		& 0.2198	    & 0.2192	    & 0.2191    &
		&                           \\
		& $p=3$ & 0.2133	    & 0.2173	    & 0.2205		& 0.2189	    & 0.2190	    & 0.2190    &                               &                           \\
		\midrule
		\multirow{3}{*}{$u_{x,A}$}  & $p=1$ &   &   &   & 0.0357        & 0.0333        & 0.0315    & \multirow{2}{*}{0.0265}       & \multirow{2}{*}{0.0307}   \\
		& $p=2$ & 0.0236	    & 0.0268	    & 0.0282        & 0.0309        & 0.0302      & 0.0298    &
		&                           \\
		& $p=3$ & 0.0251	    & 0.0270	    & 0.0293        & 0.0302        & 0.0299        & 0.0296    &                               &                           \\
		\bottomrule
	\end{tabular}
\end{table}

\subsection{Annulus Subject to Tension and Twist}\label{subsec:benchmarks_annulus}
As a next example, an annular planar surface subject to an out-of-plane translation and a twist of $90^\circ$ of the inner-boundary $\Gamma_i$ is modelled, see \cref{fig:benchmarks_annulus_setup}. The example is inspired by the work of \cite{Taylor2014} using the same geometric dimensions but different material parameters. An incompressible Neo-Hookean material model with the strain energy density function from \cref{eq:NH_incomp} is used as constitutive model, with the parameters from \cref{fig:benchmarks_annulus_setup}. The geometry is modelled using 4 patches representing a quarter annulus of degree $p=2$ with maximum regularity. The interfaces of the patches are smoothened with $C^1$ continuity over the interfaces using basic interface smoothing, i.e. by replacing every interface basis function with a weighted linear combination of its neighbours orthogonal to the patch interface. We refer to \cite{Toshniwal2017}, where this approach is used for regular interface functions, for more details. Every patch consists of one element. For the tension field theory membrane simulation, quadratic or cubic splines with $8\times8$, $16\times16$ or $32\times32$ elements per patch are used. The simulations are compared to Kirchhoff--Love on the same number of elements, but only with cubic degree to improve accuracy. \\

The example is solved using displacement controlled stepping for the rotation of the inner-boundary. In every displacement step, the example is solved using a Dynamic Relaxation method with tolerance $\epsilon_{\VEC{R}}=10^{-4}$ and a Newton--Raphson method with tolerance $\epsilon_{\VEC{R}}=10^{-6}$. In case of divergence or slow convergence (more than 50 iterations) of the Newton--Raphson iterations, the displacement step is bisected and re-started. After four subsequent failures of a displacement step, the solution obtained from the Dynamic Relaxation stage is accepted and the non-converged solution from the Newton--Raphson method is discarded. \\

\begin{figure}
	\centering
	\begin{minipage}{0.6\linewidth}
		\centering
		\resizebox{\linewidth}{!}
		{
			\tdplotsetmaincoords{60}{-30}
			\begin{tikzpicture}[
				scale=2,
				tdplot_main_coords,bullet/.style={circle,inner
					sep=1pt,fill=black,fill opacity=1}]
				\begin{scope}[canvas is xy plane at z=0]
					\filldraw[thick,shade,left color=gray,bottom color=white] (0,0) circle [radius=2];
					\filldraw[thick,fill=white] (0,0) circle [radius=1];
					\node[below] at(0,-1) {$\Gamma_i$};
					\node[right] at(2,0) {$\Gamma_o$};
					\draw[-latex] (0.9,0) arc [start angle=0,end angle=90,radius=0.9] node[midway,below left]{$\varTheta$};
					\draw[-latex] (0,0) -- (-2,0) node[midway,below]{$R_o$};
					\draw[-latex] (0,0) -- (-0.707,-0.707) node[midway,left]{$R_i$};
				\end{scope}
				\foreach \t in {0,30,...,360}
				{
					\pgfmathsetmacro\Xold{cos(\t)}
					\pgfmathsetmacro\Yold{sin(\t)}
					\pgfmathsetmacro\Xnew{cos(\t)}
					\pgfmathsetmacro\Ynew{sin(\t)}
					\draw[-latex](\Xold,\Yold,0) -- (\Xnew,\Ynew,0.5);
				}
				\begin{scope}[canvas is xy plane at z=0.5]
					\draw[thick] (0,0) circle [radius=1];
					\node[right] at(1,0) {$u_z$};
				\end{scope}
			\end{tikzpicture}
		}
	\end{minipage}
	\hfill
	\begin{minipage}{0.3\linewidth}
		\centering
		\footnotesize
		\begin{tabular}{llr}
			\toprule
			\multicolumn{3}{c}{Geometry}\\
			\midrule
			$R_i$ & 62.5  & $[\text{mm}]$\\
			$R_o$ & 250 & $[\text{mm}]$\\
			$t$ & 0.05 & $[\text{mm}]$\\
			\midrule
			\multicolumn{3}{c}{Material}\\
			\multicolumn{3}{c}{\textit{Incompressible Neo-Hookean}}\\
			\midrule
			$E$ & 1.0 & $[\text{GPa}]$\\
			$\nu$ & 0.5 & $[-]$\\
			\midrule
			\multicolumn{3}{c}{Boundary Conditions}\\
			\midrule
			$\varTheta$ & $\pi/2$ & $[\text{rad}]$\\
			$u_z$ & $125$ & $[\text{mm}]$\\
			\bottomrule
		\end{tabular}
	\end{minipage}
	\caption{Problem definition for an annulus with inner radius $R_i$ and outer radius $R_o$ subject to a vertical translation $u_z$ and a rotation $\varTheta$ on the inner boundary $\Gamma_i$ while being fixed on the outer boundary $\Gamma_o$. The annulus has a Neo-Hookean material model with the parameters provided in the table on the right.}
	\label{fig:benchmarks_annulus_setup}
\end{figure}

In \cref{fig:benchmarks_annulus}, the displacement fields are provided for the benchmark results obtained using the Kirchhoff--Love shell model, along with the displacement and tension fields obtained using the proposed membrane model. The results are provided for meshes with $64\times64$ elements  per patch and degree $p=3$. In addition, \cref{fig:benchmarks_annulus} provides selected contour lines of the displaced geometry with solid lines indicating the Kirchhoff--Love shell results and dashed lines indicating the membrane model results for different meshes and degrees. All results in \cref{fig:benchmarks_annulus} and \cref{fig:benchmarks_annulus_contours} are provided for applied angles $\{0.35,0.40,0.45,0.50\}\:[\pi\:\text{rad}]$.\\

From the results in \cref{fig:benchmarks_annulus,fig:benchmarks_annulus_contours}, multiple observations can be made. Firstly, \cref{fig:benchmarks_annulus_side} and \cref{fig:benchmarks_annulus_side_tensionfield} show a qualitatively good comparison between the wrinkled region given by the Kirchhoff--Love shell model and the tension field from the membrane model. Indeed, as shown by the contours in \cref{fig:benchmarks_annulus_contours}, the tension field theory membrane model accurately predicts the average of the wrinkles given by the wrinkling simulation using the shell model, which is as expected from the construction of the model. Furthermore, the performance of the membrane model is consistent across refinements and degree elevations, starting from coarse meshes compared to the Kirchhoff--Love shell model.\\

\begin{figure}[tb!]
	\centering
	\begin{subfigure}[t]{0.45\linewidth}
		\centering
		\includegraphics[width=\linewidth]{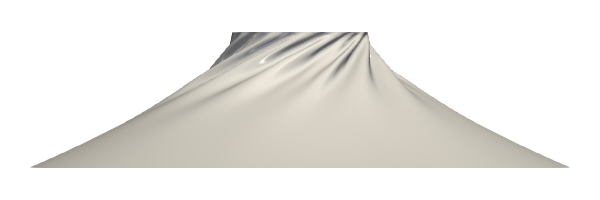}
		\includegraphics[width=\linewidth]{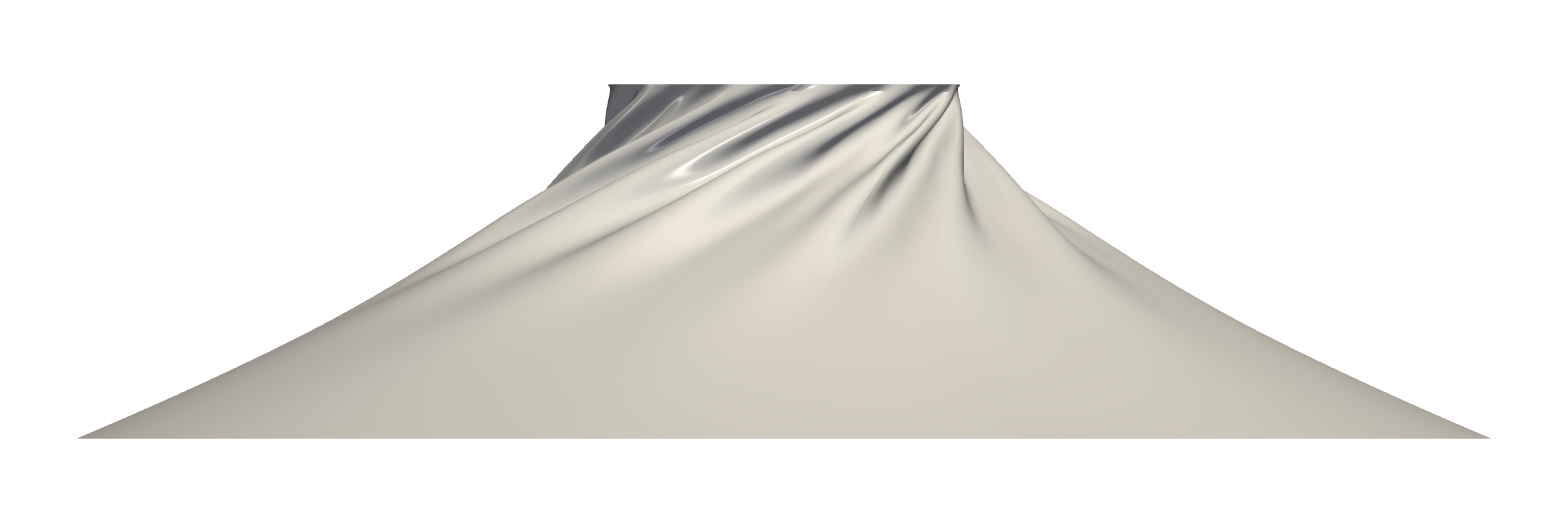}
		\includegraphics[width=\linewidth]{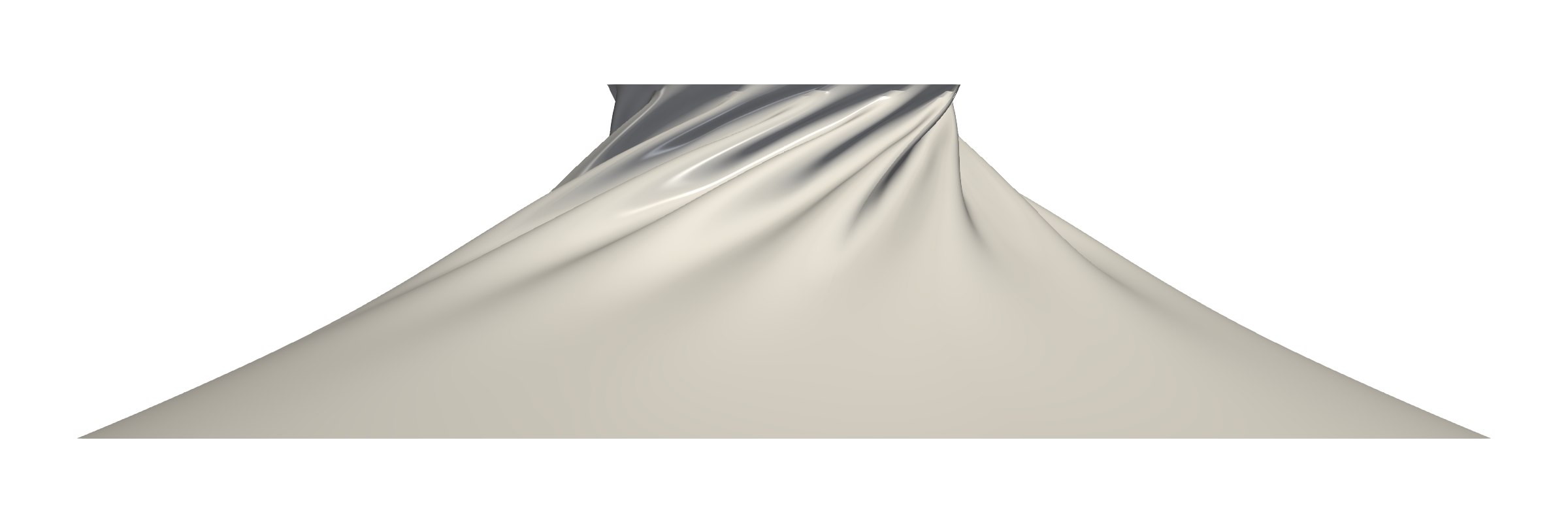}		\includegraphics[width=\linewidth]{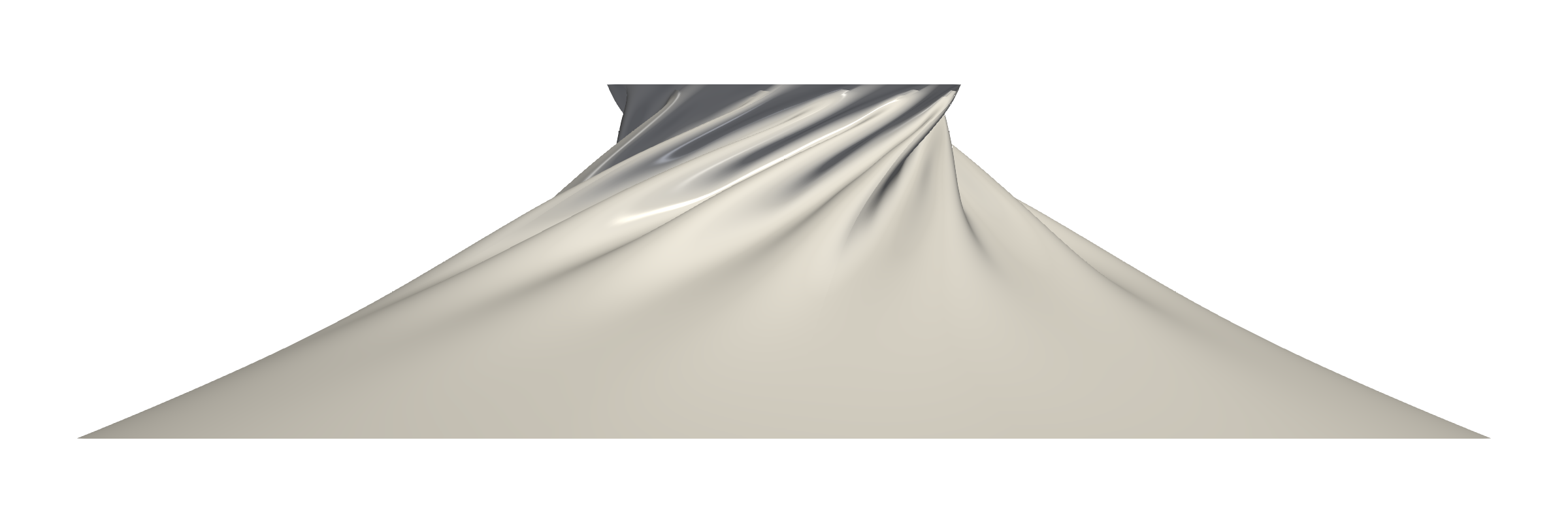}
		\caption{Side view of the deformed annulus using the Kirchhoff--Love shell model with a rotation of $\{0.35,0.40,0.45,0.50\}\:[\pi\:\text{rad}]$ from top to bottom.}
		\label{fig:benchmarks_annulus_side}
	\end{subfigure}
	\hfill
	\begin{subfigure}[t]{0.45\linewidth}
		\centering
		\includegraphics[width=\linewidth]{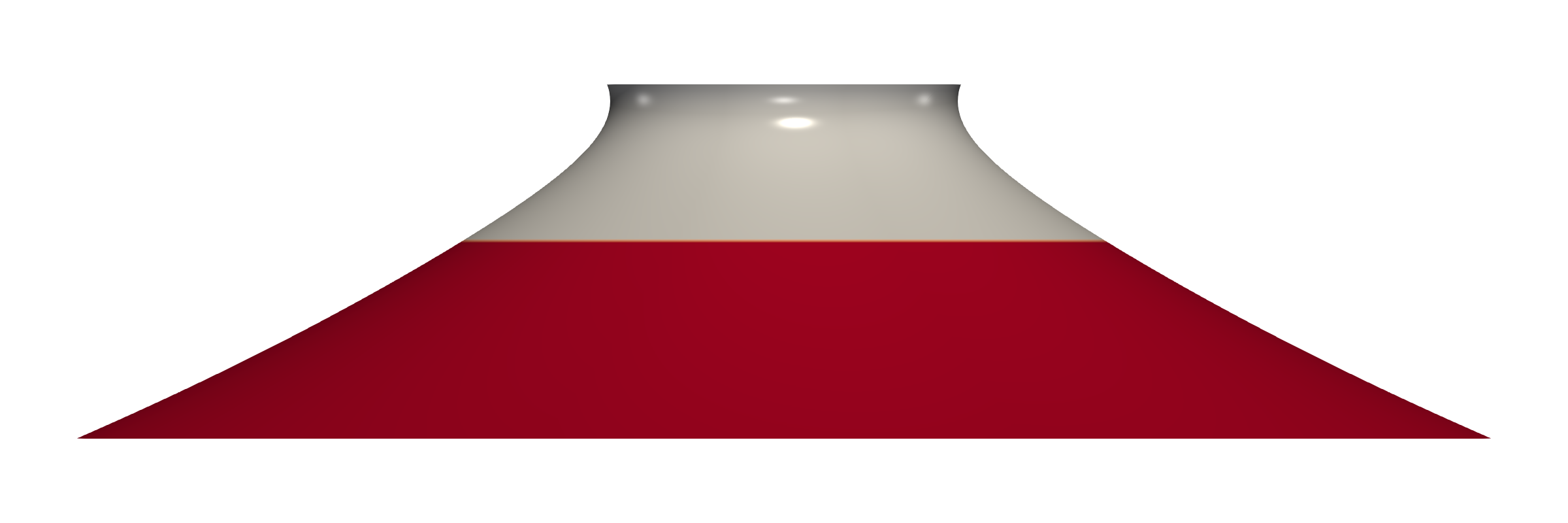}
		\includegraphics[width=\linewidth]{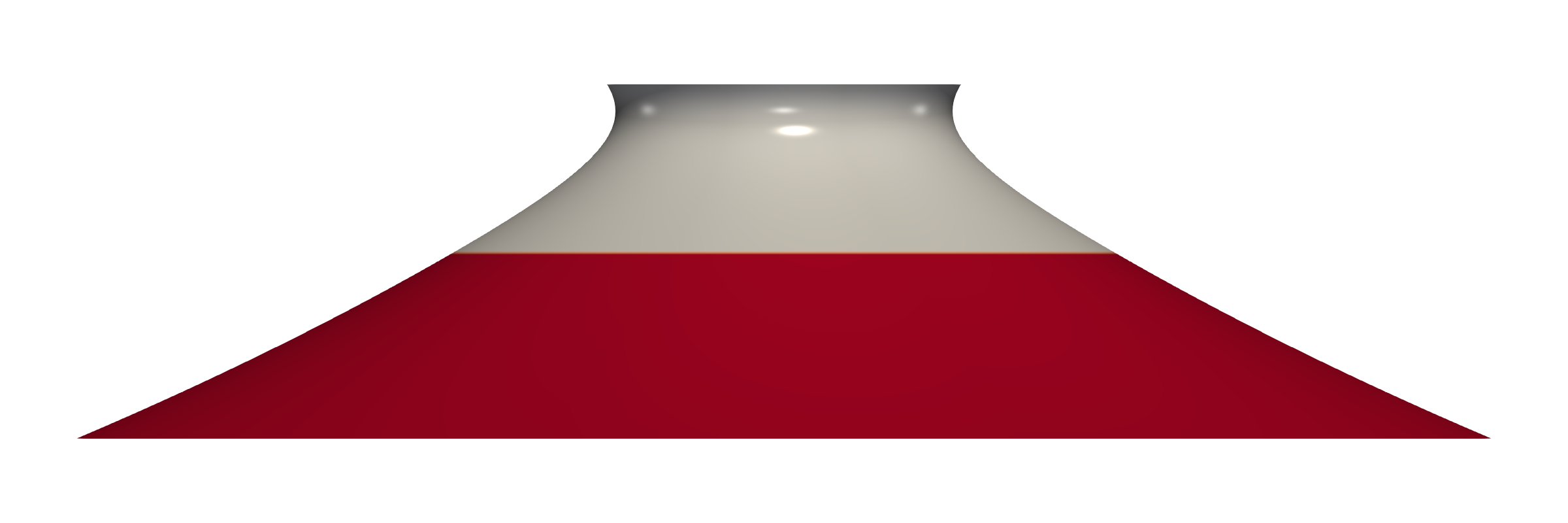}
		\includegraphics[width=\linewidth]{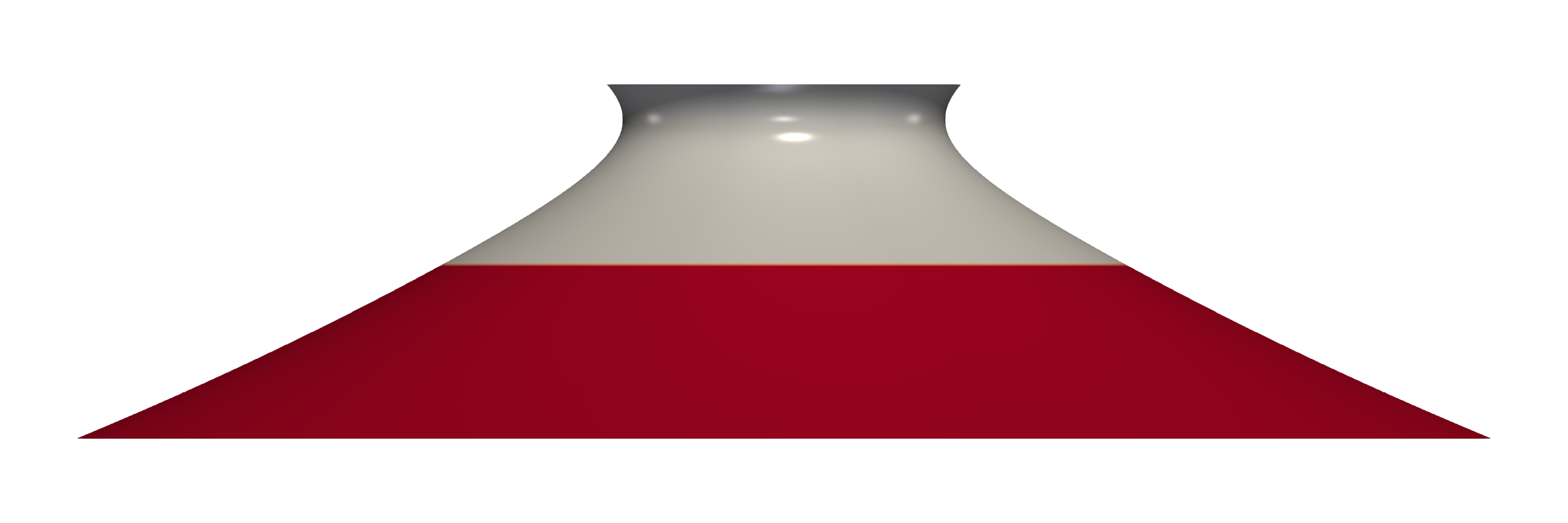}
		\includegraphics[width=\linewidth]{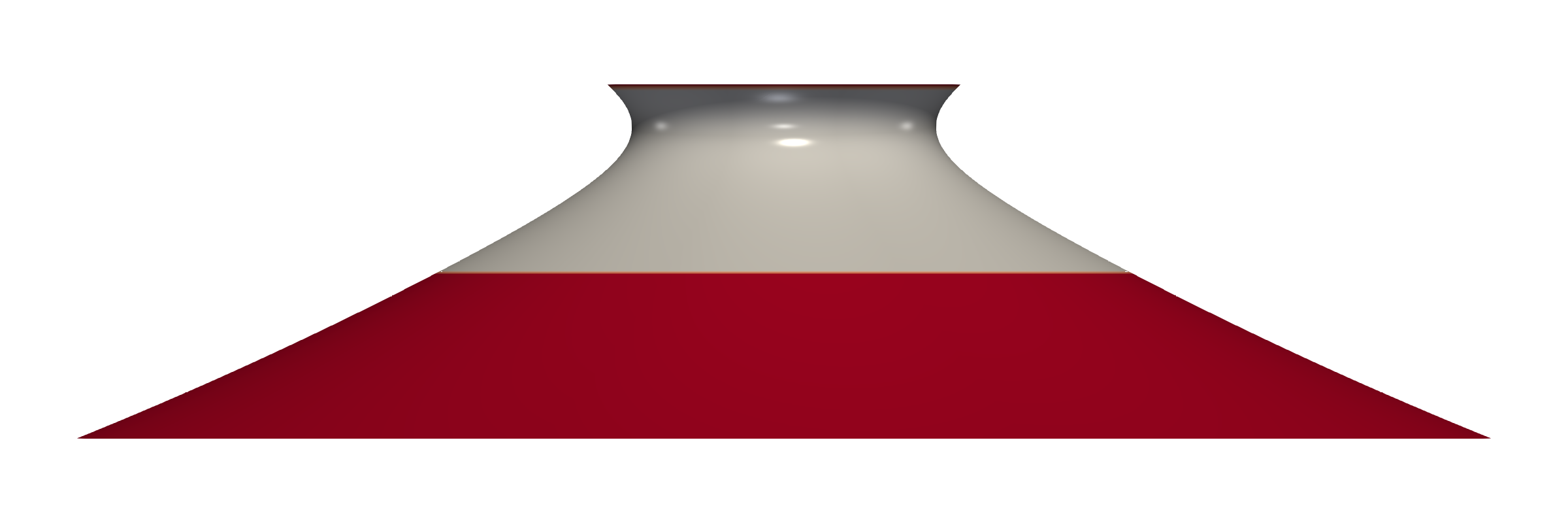}
		\caption{Side view of the deformed annulus using the tension field theory membrane model with a rotation of $\{0.35,0.40,0.45,0.50\}\:[\pi\:\text{rad}]$ from top to bottom. The red region denotes a taut, and the gray region denotes a wrinkled region.}
		\label{fig:benchmarks_annulus_side_tensionfield}
	\end{subfigure}
	 \caption{(Caption on next page)}
	 \end{figure}

 \begin{figure}
	 \ContinuedFloat
	\centering
	\begin{subfigure}[t]{0.45\linewidth}
		\centering
		\includegraphics[width=\linewidth]{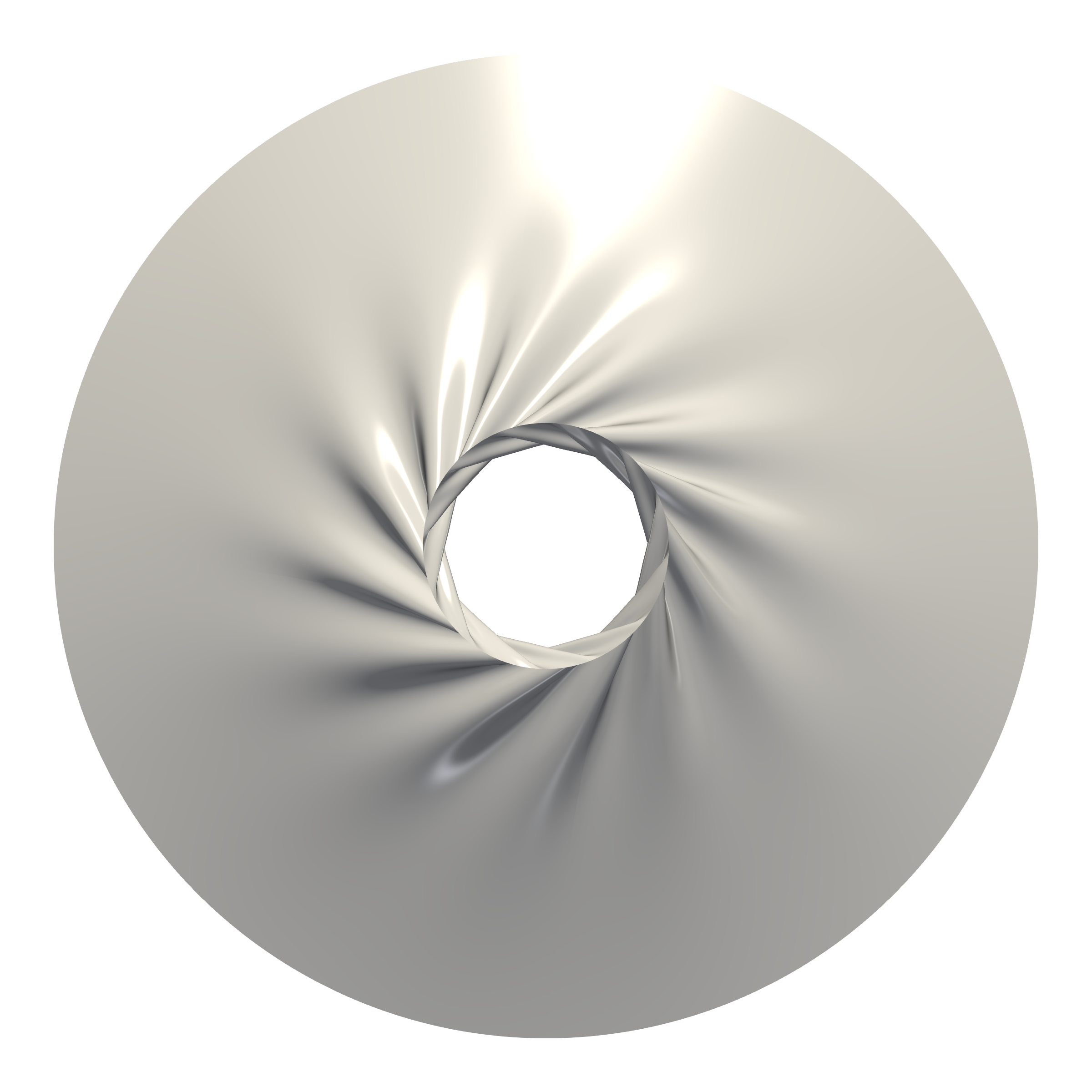}
		\caption{Top view of the deformed annulus using the Kirchhoff--Love shell model.}
		\label{fig:benchmarks_annulus_top}
	\end{subfigure}
	\hfill
	\begin{subfigure}[t]{0.45\linewidth}
		\centering
		\includegraphics[width=\linewidth]{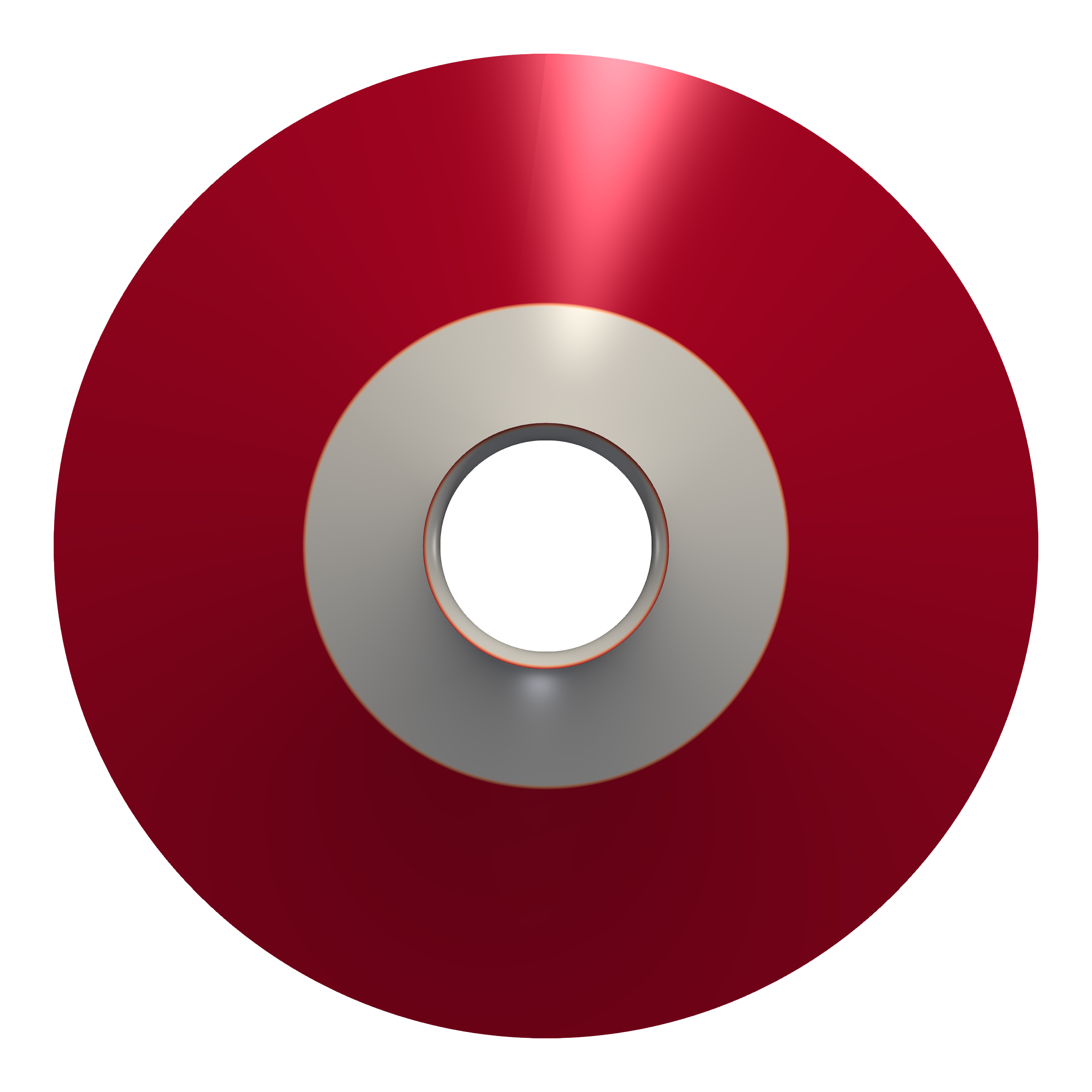}
		\caption{Top view of the deformed annulus using the tension field theory membrane model. The red region denotes a taut, and the gray region denotes a wrinkled region.}
		\label{fig:benchmarks_annulus_top_tensionfield}
	\end{subfigure}
	\caption{Results of the example with an annulus with fixed outer boundary and with an inner boundary subject to a translation and a rotation, see \cref{fig:benchmarks_annulus_setup}. A side view (\subref{fig:benchmarks_annulus_side}) and a top view (\subref{fig:benchmarks_annulus_top}) of the wrinkled membrane using the Kirchhoff--Love shell simulation are provided, as well as the deformed geometry from the tension field theory membrane simulation with the tension field for colouring (\subref{fig:benchmarks_annulus_side_tensionfield}). Furthermore, contour lines of the deformation for different parametric coordinates are provided; see (\subref{fig:benchmarks_annulus_top_tensionfield}).}
	\label{fig:benchmarks_annulus}
\end{figure}

\begin{figure}
\centering
	\begin{subfigure}[t]{0.225\linewidth}
		\centering
		\begin{tikzpicture}
			\begin{axis}
				[
				height=\linewidth,
				width=\linewidth,
				scale only axis,
				hide axis,
				]
				\addplot+[style=style1,solid,opacity=0.5,no markers] table[col sep=comma,header=true,x index={0},y index={1}] {annulus_e1_r2_line_D=0.700000_u=0.9_dir=1.csv};
				\addplot+[style=style2,solid,opacity=0.5,no markers] table[col sep=comma,header=true,x index={0},y index={1}] {annulus_e1_r3_line_D=0.700000_u=0.9_dir=1.csv};
				\addplot+[style=style3,solid,opacity=0.5,no markers] table[col sep=comma,header=true,x index={0},y index={1}] {annulus_e1_r4_line_D=0.700000_u=0.9_dir=1.csv};

				\addplot+[style=style1,dotted,no markers] table[col sep=comma,header=true,x index={0},y index={1}] {annulus_e0_r2_TFT_line_D=0.700000_u=0.9_dir=1.csv};
				\addplot+[style=style2,dotted,no markers] table[col sep=comma,header=true,x index={0},y index={1}] {annulus_e0_r3_TFT_line_D=0.700000_u=0.9_dir=1.csv};
				\addplot+[style=style3,dotted,no markers] table[col sep=comma,header=true,x index={0},y index={1}] {annulus_e0_r4_TFT_line_D=0.700000_u=0.9_dir=1.csv};

				\addplot+[style=style1,dashed,no markers] table[col sep=comma,header=true,x index={0},y index={1}] {annulus_e1_r2_TFT_line_D=0.700000_u=0.9_dir=1.csv};
				\addplot+[style=style2,dashed,no markers] table[col sep=comma,header=true,x index={0},y index={1}] {annulus_e1_r3_TFT_line_D=0.700000_u=0.9_dir=1.csv};
				\addplot+[style=style3,dashed,no markers] table[col sep=comma,header=true,x index={0},y index={1}] {annulus_e1_r4_TFT_line_D=0.700000_u=0.9_dir=1.csv};
			\end{axis}
		\end{tikzpicture}
		\caption*{$0.35\:[\pi\:\text{rad}]$}
	\end{subfigure}
	\hfill
	\begin{subfigure}[t]{0.225\linewidth}
		\centering
		\begin{tikzpicture}
			\begin{axis}
				[
				height=\linewidth,
				width=\linewidth,
				scale only axis,
				hide axis,
				]
				\addplot+[style=style1,solid,opacity=0.5,no markers] table[col sep=comma,header=true,x index={0},y index={1}] {annulus_e1_r2_line_D=0.800000_u=0.9_dir=1.csv};
				\addplot+[style=style2,solid,opacity=0.5,no markers] table[col sep=comma,header=true,x index={0},y index={1}] {annulus_e1_r3_line_D=0.800000_u=0.9_dir=1.csv};
				\addplot+[style=style3,solid,opacity=0.5,no markers] table[col sep=comma,header=true,x index={0},y index={1}] {annulus_e1_r4_line_D=0.800000_u=0.9_dir=1.csv};

				\addplot+[style=style1,dotted,no markers] table[col sep=comma,header=true,x index={0},y index={1}] {annulus_e0_r2_TFT_line_D=0.800000_u=0.9_dir=1.csv};
				\addplot+[style=style2,dotted,no markers] table[col sep=comma,header=true,x index={0},y index={1}] {annulus_e0_r3_TFT_line_D=0.800000_u=0.9_dir=1.csv};
				\addplot+[style=style3,dotted,no markers] table[col sep=comma,header=true,x index={0},y index={1}] {annulus_e0_r4_TFT_line_D=0.800000_u=0.9_dir=1.csv};

				\addplot+[style=style1,dashed,no markers] table[col sep=comma,header=true,x index={0},y index={1}] {annulus_e1_r2_TFT_line_D=0.800000_u=0.9_dir=1.csv};
				\addplot+[style=style2,dashed,no markers] table[col sep=comma,header=true,x index={0},y index={1}] {annulus_e1_r3_TFT_line_D=0.800000_u=0.9_dir=1.csv};
				\addplot+[style=style3,dashed,no markers] table[col sep=comma,header=true,x index={0},y index={1}] {annulus_e1_r4_TFT_line_D=0.800000_u=0.9_dir=1.csv};
			\end{axis}
		\end{tikzpicture}
		\caption*{$0.40\:[\pi\:\text{rad}]$}
	\end{subfigure}
	\hfill
	\begin{subfigure}[t]{0.225\linewidth}
		\centering
		\begin{tikzpicture}
			\begin{axis}
				[
				height=\linewidth,
				width=\linewidth,
				scale only axis,
				hide axis,
				]
				\addplot+[style=style1,solid,opacity=0.5,no markers] table[col sep=comma,header=true,x index={0},y index={1}] {annulus_e1_r2_line_D=0.900000_u=0.9_dir=1.csv};
				\addplot+[style=style2,solid,opacity=0.5,no markers] table[col sep=comma,header=true,x index={0},y index={1}] {annulus_e1_r3_line_D=0.900000_u=0.9_dir=1.csv};
				\addplot+[style=style3,solid,opacity=0.5,no markers] table[col sep=comma,header=true,x index={0},y index={1}] {annulus_e1_r4_line_D=0.900000_u=0.9_dir=1.csv};

				\addplot+[style=style1,dotted,no markers] table[col sep=comma,header=true,x index={0},y index={1}] {annulus_e0_r2_TFT_line_D=0.900000_u=0.9_dir=1.csv};
				\addplot+[style=style2,dotted,no markers] table[col sep=comma,header=true,x index={0},y index={1}] {annulus_e0_r3_TFT_line_D=0.900000_u=0.9_dir=1.csv};
				\addplot+[style=style3,dotted,no markers] table[col sep=comma,header=true,x index={0},y index={1}] {annulus_e0_r4_TFT_line_D=0.900000_u=0.9_dir=1.csv};

				\addplot+[style=style1,dashed,no markers] table[col sep=comma,header=true,x index={0},y index={1}] {annulus_e1_r2_TFT_line_D=0.900000_u=0.9_dir=1.csv};
				\addplot+[style=style2,dashed,no markers] table[col sep=comma,header=true,x index={0},y index={1}] {annulus_e1_r3_TFT_line_D=0.900000_u=0.9_dir=1.csv};
				\addplot+[style=style3,dashed,no markers] table[col sep=comma,header=true,x index={0},y index={1}] {annulus_e1_r4_TFT_line_D=0.900000_u=0.9_dir=1.csv};
			\end{axis}
		\end{tikzpicture}
		\caption*{$0.45\:[\pi\:\text{rad}]$}
	\end{subfigure}
	\hfill
	\begin{subfigure}[t]{0.225\linewidth}
		\centering
		\begin{tikzpicture}
			\begin{axis}
				[
				height=\linewidth,
				width=\linewidth,
				scale only axis,
				hide axis,
				]
				\addplot+[style=style1,solid,opacity=0.5,no markers] table[col sep=comma,header=true,x index={0},y index={1}] {annulus_e1_r2_line_D=1.000000_u=0.9_dir=1.csv};
				\addplot+[style=style2,solid,opacity=0.5,no markers] table[col sep=comma,header=true,x index={0},y index={1}] {annulus_e1_r3_line_D=1.000000_u=0.9_dir=1.csv};
				\addplot+[style=style3,solid,opacity=0.5,no markers] table[col sep=comma,header=true,x index={0},y index={1}] {annulus_e1_r4_line_D=1.000000_u=0.9_dir=1.csv};

				\addplot+[style=style1,dotted,no markers] table[col sep=comma,header=true,x index={0},y index={1}] {annulus_e0_r2_TFT_line_D=1.000000_u=0.9_dir=1.csv};
				\addplot+[style=style2,dotted,no markers] table[col sep=comma,header=true,x index={0},y index={1}] {annulus_e0_r3_TFT_line_D=1.000000_u=0.9_dir=1.csv};
				\addplot+[style=style3,dotted,no markers] table[col sep=comma,header=true,x index={0},y index={1}] {annulus_e0_r4_TFT_line_D=1.000000_u=0.9_dir=1.csv};

				\addplot+[style=style1,dashed,no markers] table[col sep=comma,header=true,x index={0},y index={1}] {annulus_e1_r2_TFT_line_D=1.000000_u=0.9_dir=1.csv};
				\addplot+[style=style2,dashed,no markers] table[col sep=comma,header=true,x index={0},y index={1}] {annulus_e1_r3_TFT_line_D=1.000000_u=0.9_dir=1.csv};
				\addplot+[style=style3,dashed,no markers] table[col sep=comma,header=true,x index={0},y index={1}] {annulus_e1_r4_TFT_line_D=1.000000_u=0.9_dir=1.csv};
			\end{axis}
		\end{tikzpicture}
		\caption*{$0.50\:[\pi\:\text{rad}]$}
	\end{subfigure}
	\caption{Contour line for parametric coordinate $\xi_r=0.9$ such that $\xi_r=0$ corresponds to $R_o$ and $\xi_r=1$ with $R_i$. The solid lines represent the Kirchhoff-Love shell results, and the dashed lines represent the tension field theory membrane results. Note that the coarsest mesh for the shell simulation (blue solid line) results in a self-intersecting solution for $\pi/2\:[\text{rad}]$. See \cref{fig:benchmarks_annulus_contours} for the legend. }
	\label{fig:benchmarks_annulus_contours}
\end{figure}

\begin{figure}
	\centering
	\begin{tikzpicture}
		\begin{groupplot}
			[
			height=0.2\textheight,
			width=0.38\linewidth,
			enlarge x limits = true,
			enlarge y limits = true,
			scale only axis,
			group style={
				group name = top,
				group size=2 by 1,
				xlabels at=edge bottom,
				x descriptions at=edge bottom,
				vertical sep=0.05\textheight,
				horizontal sep=0.04\linewidth,
			},
			legend pos = south east,
			]
			\nextgroupplot[xlabel = {Rotation $[\pi\:\text{rad}]$},
			ylabel = {Applied torque $M$\:[Nm]}]

			\addplot+[style=style1,no markers,solid] table[col sep=comma,header=true,x expr=\thisrowno{0}/2, y index={1}] {annulus_e1_r2_out.csv};
			\addplot+[style=style2,no markers,solid] table[col sep=comma,header=true,x expr=\thisrowno{0}/2, y index={1}] {annulus_e1_r3_out.csv};
			\addplot+[style=style3,no markers,solid] table[col sep=comma,header=true,x expr=\thisrowno{0}/2, y index={1}] {annulus_e1_r4_out.csv};

			\addplot+[style=style1,no markers,dotted] table[col sep=comma,header=true,x expr=\thisrowno{0}/2, y index={1}] {annulus_e0_r2_TFT_out.csv};
			\addplot+[style=style2,no markers,dotted] table[col sep=comma,header=true,x expr=\thisrowno{0}/2, y index={1}] {annulus_e0_r3_TFT_out.csv};
			\addplot+[style=style3,no markers,dotted] table[col sep=comma,header=true,x expr=\thisrowno{0}/2, y index={1}] {annulus_e0_r4_TFT_out.csv};

			\addplot+[style=style1,no markers,dashed] table[col sep=comma,header=true,x expr=\thisrowno{0}/2, y index={1}] {annulus_e1_r2_TFT_out.csv};
			\addplot+[style=style2,no markers,dashed] table[col sep=comma,header=true,x expr=\thisrowno{0}/2, y index={1}] {annulus_e1_r3_TFT_out.csv};
			\addplot+[style=style3,no markers,dashed] table[col sep=comma,header=true,x expr=\thisrowno{0}/2, y index={1}] {annulus_e1_r4_TFT_out.csv};

			\draw[gray,dashed] (axis cs:0.2,280) rectangle (axis cs: 0.5,480);
			\node[anchor=south east,gray] at (axis cs:0.5,280) {Inset};

			\nextgroupplot[xlabel = {Rotation $[\pi\:\text{rad}]$},
			xmin = 0.2, xmax = 0.5,
			ymin = 280, ymax = 480,
			enlargelimits=false,
			draw=gray,
			xlabel style={color=gray},
			x tick label style={color = gray},
			ylabel style={color=gray},
			y tick label style={color = gray},
			ylabel near ticks, yticklabel pos=right]

			\addlegendimage{style=style1,solid,no markers}\addlegendentry{$16\times16$};
			\addlegendimage{style=style2,solid,no markers}\addlegendentry{$32\times32$};
			\addlegendimage{style=style3,solid,no markers}\addlegendentry{$64\times64$};
			\addlegendimage{black,solid}\addlegendentry{KL-shell, $p=3$};
			\addlegendimage{black,dotted}\addlegendentry{TFT-membrane, $p=2$};
			\addlegendimage{black,dashed}\addlegendentry{TFT-membrane, $p=3$};

			\addplot+[style=style1,no markers,solid] table[col sep=comma,header=true,x expr=\thisrowno{0}/2, y index={1}] {annulus_e1_r2_out.csv};
			\addplot+[style=style2,no markers,solid] table[col sep=comma,header=true,x expr=\thisrowno{0}/2, y index={1}] {annulus_e1_r3_out.csv};
			\addplot+[style=style3,no markers,solid] table[col sep=comma,header=true,x expr=\thisrowno{0}/2, y index={1}] {annulus_e1_r4_out.csv};

			\addplot+[style=style1,no markers,dotted] table[col sep=comma,header=true,x expr=\thisrowno{0}/2, y index={1}] {annulus_e0_r2_TFT_out.csv};
			\addplot+[style=style2,no markers,dotted] table[col sep=comma,header=true,x expr=\thisrowno{0}/2, y index={1}] {annulus_e0_r3_TFT_out.csv};
			\addplot+[style=style3,no markers,dotted] table[col sep=comma,header=true,x expr=\thisrowno{0}/2, y index={1}] {annulus_e0_r4_TFT_out.csv};

			\addplot+[style=style1,no markers,dashed] table[col sep=comma,header=true,x expr=\thisrowno{0}/2, y index={1}] {annulus_e1_r2_TFT_out.csv};
			\addplot+[style=style2,no markers,dashed] table[col sep=comma,header=true,x expr=\thisrowno{0}/2, y index={1}] {annulus_e1_r3_TFT_out.csv};
			\addplot+[style=style3,no markers,dashed] table[col sep=comma,header=true,x expr=\thisrowno{0}/2, y index={1}] {annulus_e1_r4_TFT_out.csv};

		\end{groupplot}
	\end{tikzpicture}
	\caption{Load-displacement curves of the torque $M$ applied on the inner boundary (vertical axis) versus the rotation of the inner boundary (horizontal axis). The full diagram up to a rotation of $\pi/2$ radians is given in the left figure, while an inset of the final end of the curve is given in the right figure. The results are provided for different mesh sizes for the Kirchhoff--Love shell (KL-shell) and the tension field theory membrane (TFT-membrane) models fir the considered degrees.}
	\label{fig:benchmarks_annulus_torque}
\end{figure}

In addition, the resulting torque around the $z$-axis at the top boundary is given in \cref{fig:benchmarks_annulus_torque}. This quantity is computed by applying the variational formulation from \cref{eq:dW} on the relevant boundary, with the obtained solution $\uvec$ and the distance to the center of rotation as variation. As can be seen in the figure, the tension field membrane results provide a good approximation of the results obtained by the finest mesh of the Kirchhoff--Love shell simulation, even for the coarsest mesh. However, the resulting torque for the angle $\pi\:[\text{rad}]$ at the end of the simulation seems to be underestimated by the membrane model. An explanation for the substantial differences is unclear yet, but most likely the modelling assumption of a fully vanishing stress orthogonal to the wrinkle direction is invalid, e.g. because of remaining bending stiffness.\\

\subsection{Cylinder Subject to Tension and Twist}\label{subsec:benchmarks_cylinder}
Similar to the previous example regarding the pulled hyperelastic annulus (see \cref{subsec:benchmarks_annulus}), the next benchmark models a cylinder subject to a translation along its length and a rotation around its centre axis (see \cref{fig:benchmarks_cylinder_setup}). As for the annulus in \cref{subsec:benchmarks_annulus}, an incompressible Neo-Hookean material model is used. As in the previous benchmark, the geometry is modelled using quadratic patches with a smoothed basis over the interfaces. The degrees and mesh sizes of the bases used for the tension field membrane element and the Kirchhoff--Love shell are the same as in the previous benchmark. In addition, the solver setting is the same as in the previous benchmark, and the displacement of the top boundary is controlled. Although there are similarities in the problem set-up between the annulus from \cref{subsec:benchmarks_annulus} and the cylinder benchmark, it should be noted that the case of the cylinder involves larger strains, hence this example is more suitable for hyperelastic material models.\\

\begin{figure}[tb!]
	\centering
	\begin{minipage}{0.6\linewidth}
		\centering
		\begin{tikzpicture}[scale=0.7]
			\def\R{2}
			\begin{scope}[rotate=00]
				\begin{axis}[
					axis equal image,
					hide axis,
					width=\linewidth,
					view = {20}{10},
					scale = 2,
					xmin = 0.0,
					xmax = 8.0,
					ymin = -\R-0.5,
					ymax = \R+0.5,
					zmin = -\R-0.5,
					zmax = \R+0.5,
					line cap=round
					]\

					\addplot3[domain=-0.5*pi:0.5*pi, samples=100, samples y=0, no marks, smooth, thick,gray](
					{0},
					{\R*cos(deg(\x)},
					{\R*sin(deg(\x))}
					);

					\addplot3[
					surf,
					samples = 20,
					fill opacity=0.5,
					samples y = 2,
					domain = 0.0*pi:2*pi,
					domain y = 0:3,
					draw=none,
					z buffer = sort,
					no marks,
					mesh/interior colormap={blueblack}{color=(black) color=(white)},
					colormap ={blueblack}{color=(black) color=(white)},
					shader=interp,
					opacity=0.5,
					point meta=x+3*y*y-0.25*z,
					](
					{2*\y},
					{(\R)*cos(deg(\x)},
					{(\R)*sin(deg(\x))}
					);

					\addplot3[domain=-0.5*pi:0.5*pi, samples=100, samples y=0, no marks, smooth, thick,black](
					{0},
					{-\R*cos(deg(\x)},
					{-\R*sin(deg(\x))}
					);

					\addplot3[domain=0:2*pi, samples=100, samples y=0, no marks, smooth, thick,black](
					{6},
					{\R*cos(deg(\x)},
					{\R*sin(deg(\x))}
					);

					\addplot3[domain=0.0*pi:0.5*pi, samples=100, samples y=0, no marks, smooth, thick,black,latex-](
					{7.5},
					{\R*cos(deg(-\x)},
					{\R*sin(deg(-\x))}
					);

					\foreach \t in {0.0,0.1,...,2} {

						\edef\x{\R*cos(deg(\t*3.1415} 
						\edef\y{\R*sin(deg(\t*3.1415} 
						\edef
						\temp{
							\noexpand
							\draw [-latex] (axis cs:6,\x,\y) -- (axis cs:7,\x,\y);
						}
						\temp
					}

					\addplot3[domain=0.0*pi:2*pi, samples=100, samples y=0, no marks, smooth, thick](
					{7},
					{\R*cos(deg(\x)},
					{\R*sin(deg(\x))}
					);



					\node[right] at (axis cs: 7,0.707*\R,0.707*\R) {$u_x$};
					\node[right] at (axis cs: 7.5,0.707*\R,-0.707*\R) {$\varTheta$};

					\draw[thick] (axis cs: 0,0,\R) -- (axis cs: 6,0,\R);
					\draw[thick] (axis cs: 0,0,-\R) -- (axis cs: 6,0,-\R);

					\draw[latex-latex] (axis cs: 0,0,0) -- (axis cs: 0,-0.707*\R,-0.707*\R) node[midway, right]{$R$};

					\draw[latex-latex] (axis cs: 0,0,-2.2) -- (axis cs: 6,0,-2.2) node[midway,below]{$H$};
					\node at (axis cs: 0,0,0) {$\times$};

					\draw[-latex] (axis cs: 2,0,0) -- (axis cs: 2,0,1.0) node[above]{$z$};
					\draw[-latex] (axis cs: 2,0,0) -- (axis cs: 2,2.0,0) node[above right]{$y$};
					\draw[-latex] (axis cs: 2,0,0) -- (axis cs: 3,0,0) node[right]{$x$};

				\end{axis}
			\end{scope}
		\end{tikzpicture}
	\end{minipage}
	\hfill
	\begin{minipage}{0.3\linewidth}
		\centering
		\footnotesize
		\begin{tabular}{llr}
			\toprule
			\multicolumn{3}{c}{Geometry}\\
			\midrule
			$R$ & 250  & $[\text{mm}]$\\
			$L$ & 1.0 & $[\text{m}]$\\
			$t$ & 0.05 & $[\text{mm}]$\\
			\midrule
			\multicolumn{3}{c}{Material}\\
			\multicolumn{3}{c}{\textit{Incompressible Neo-Hookean}}\\
			\midrule
			$E$ & 1.0 & $[\text{GPa}]$\\
			$\nu$ & 0.5 & $[-]$\\
			\midrule
			\multicolumn{3}{c}{Boundary Conditions}\\
			\midrule
			$\varTheta$ & $\pi/2$ & $[\text{rad}]$\\
			$u_x$ & $1.0$ & $[\text{m}]$\\
			\bottomrule
		\end{tabular}
	\end{minipage}

	\caption{Problem definition for an cylinder with inner radius $R$ and length $L$ subject to an elongation $u_x$ and a rotation $\varTheta$ on the right boundary $\Gamma_r$ while being fixed on the left boundary $\Gamma_l$. The cylinder has a Neo-Hookean material model with the parameters provided in the table on the right.}
	\label{fig:benchmarks_cylinder_setup}
\end{figure}

\Cref{fig:benchmarks_cylinder} provides the results obtained using the Kirchhoff--Love shell model, along with the displacement and tension fields obtained using the proposed membrane model and \cref{fig:benchmarks_cylinder_contours} provides the contour line of the deformed geometries halfway the height at different load steps for both models. The results in \cref{fig:benchmarks_cylinder} are obtained for meshes with $64\times64$ elements with degree $p=3$ and the results in \cref{fig:benchmarks_cylinder_contours} are obtained for meshes with $16\times16$, $32\times32$ and $64\times64$ elements of degree $p=3$ for the Kirchhoff--Love shell and degrees $p=2$ and $p=3$ for the proposed membrane model. All results are provided for applied angles $\{0.35,0.40,0.45,0.50\}\:[\pi\:\text{rad}]$.\\


\begin{figure}
	\centering
	\begin{subfigure}[t]{0.45\linewidth}
		\centering
		\includegraphics[width=\linewidth]{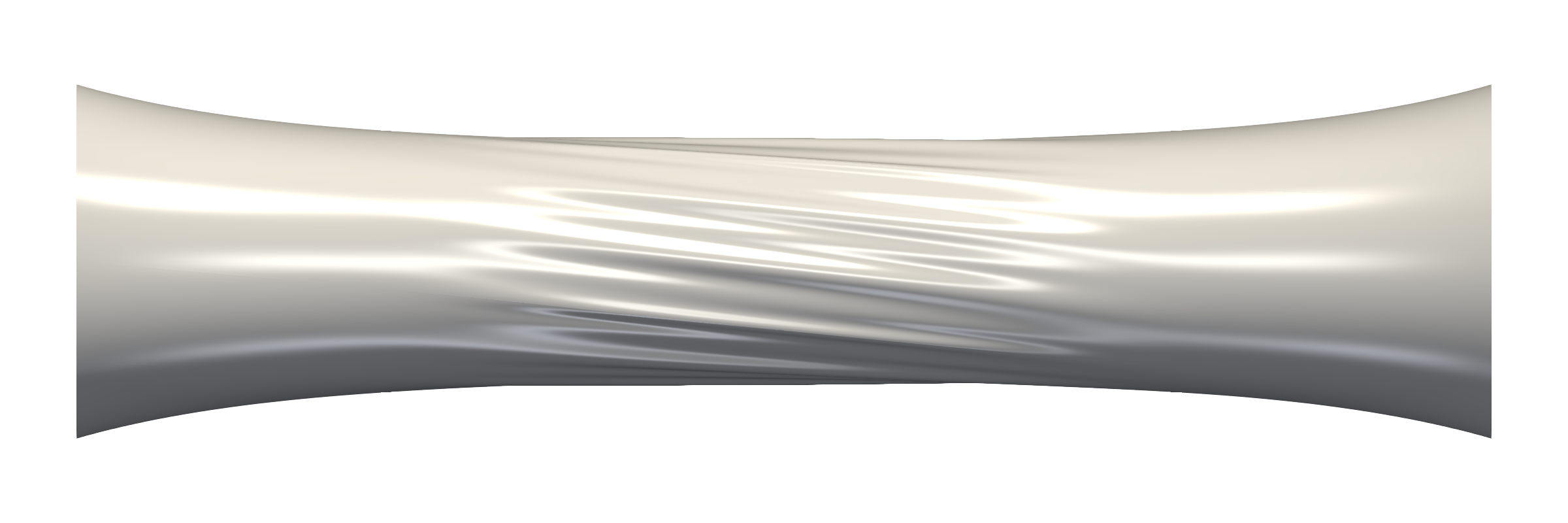}
		\includegraphics[width=\linewidth]{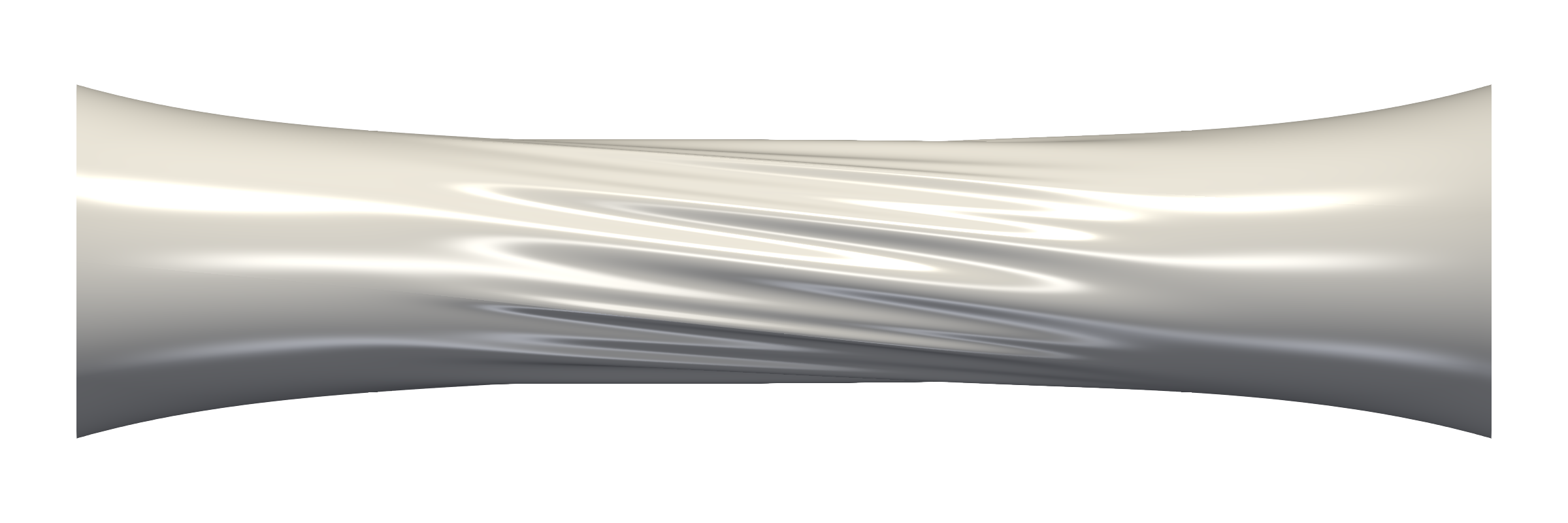}
		\includegraphics[width=\linewidth]{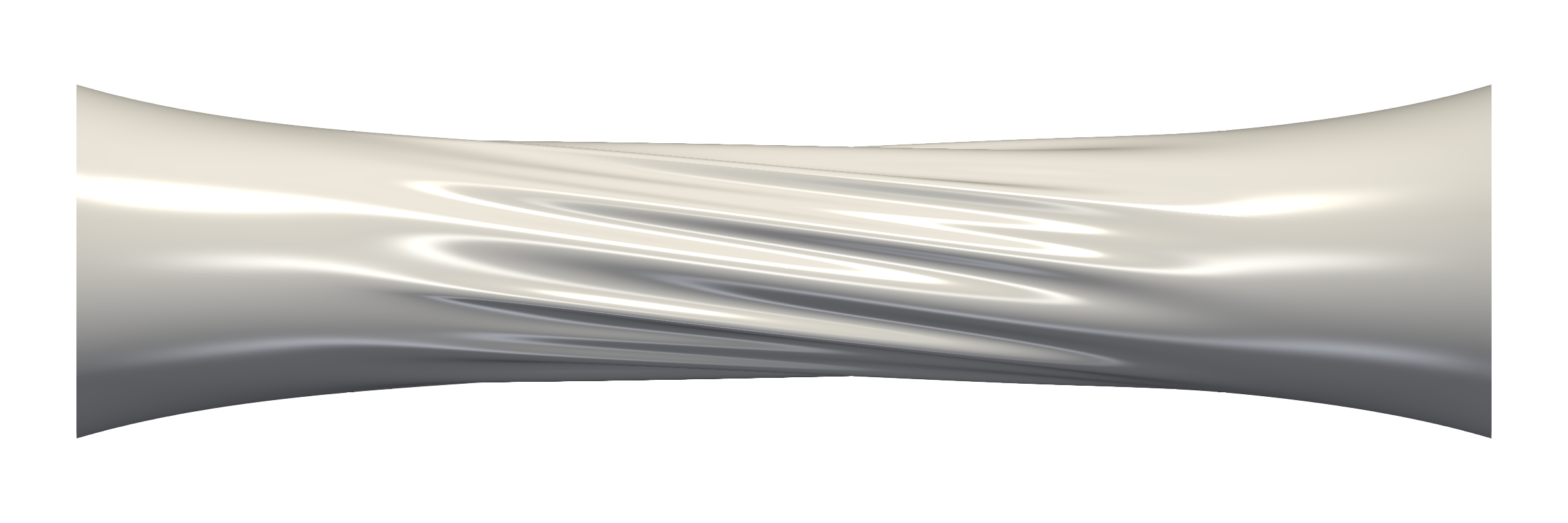}		\includegraphics[width=\linewidth]{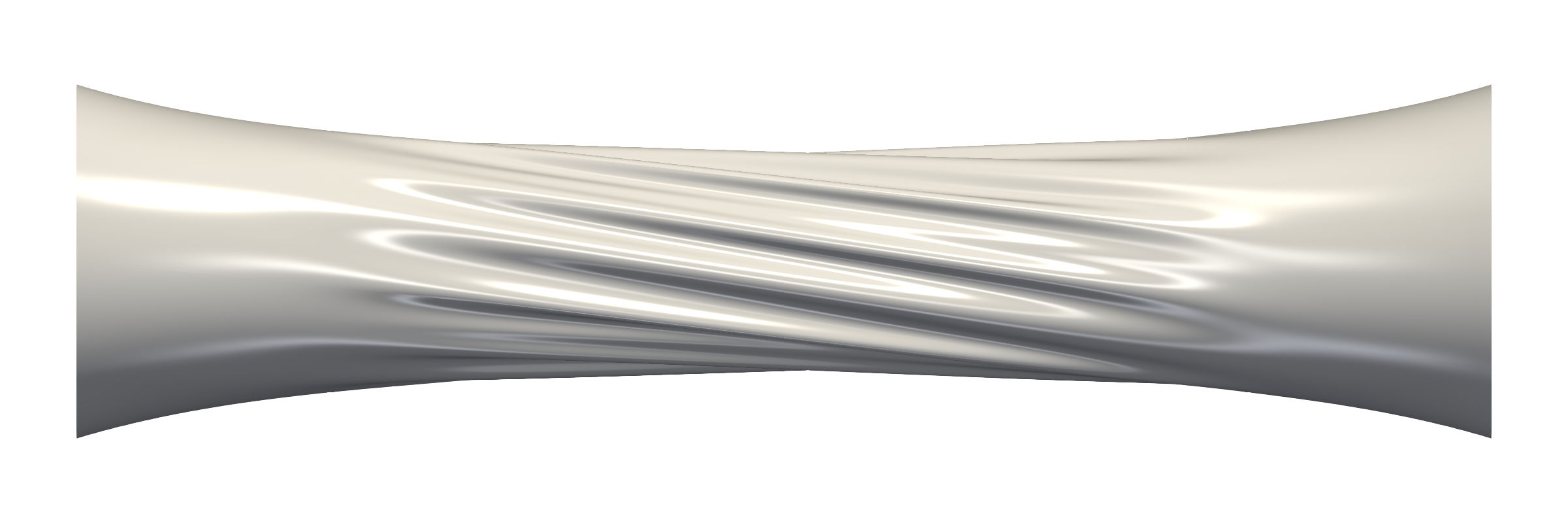}
		\caption{Side view of the deformed cylinder using the Kirchhoff--Love shell model with a rotation of $\{0.35,0.40,0.45,0.50\}\:[\pi\:\text{rad}]$ from top to bottom.}
		\label{fig:benchmarks_cylinder_side}
	\end{subfigure}
	\hfill
	\begin{subfigure}[t]{0.45\linewidth}
		\centering
		\includegraphics[width=\linewidth]{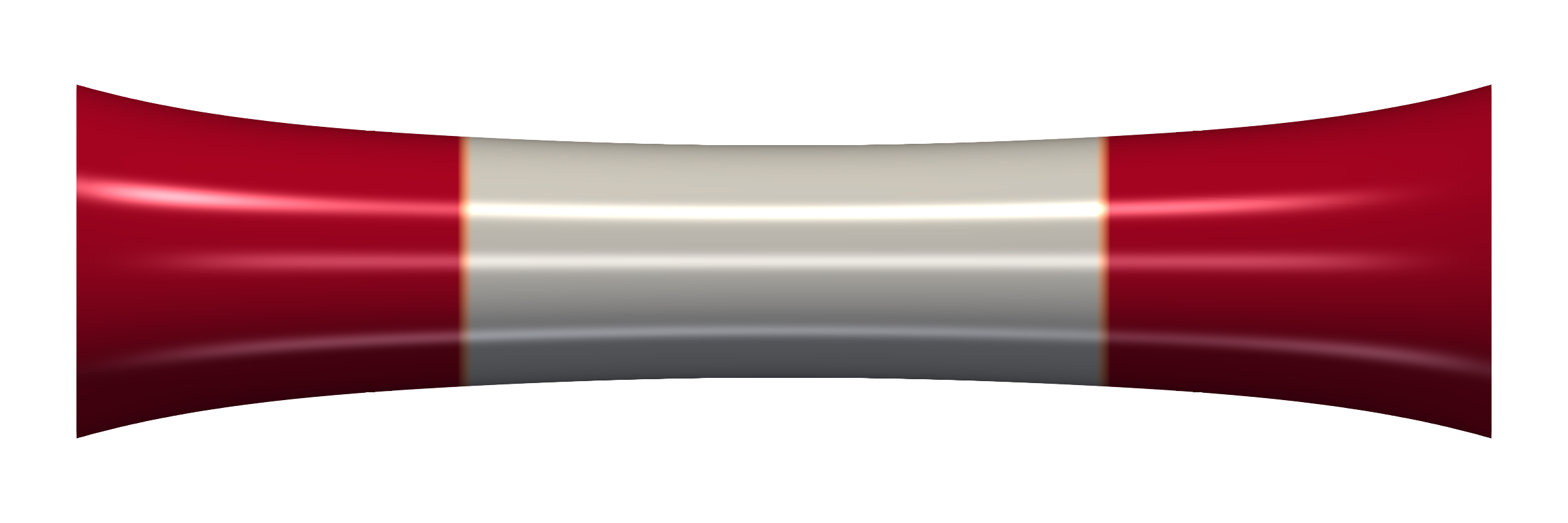}
		\includegraphics[width=\linewidth]{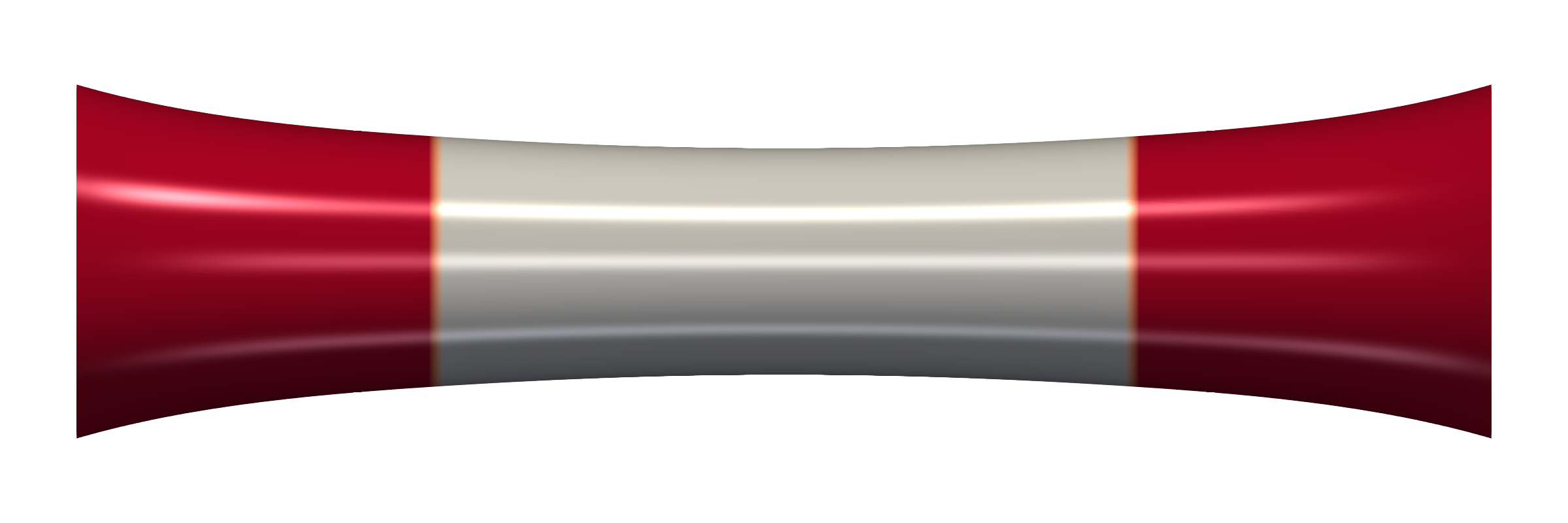}
		\includegraphics[width=\linewidth]{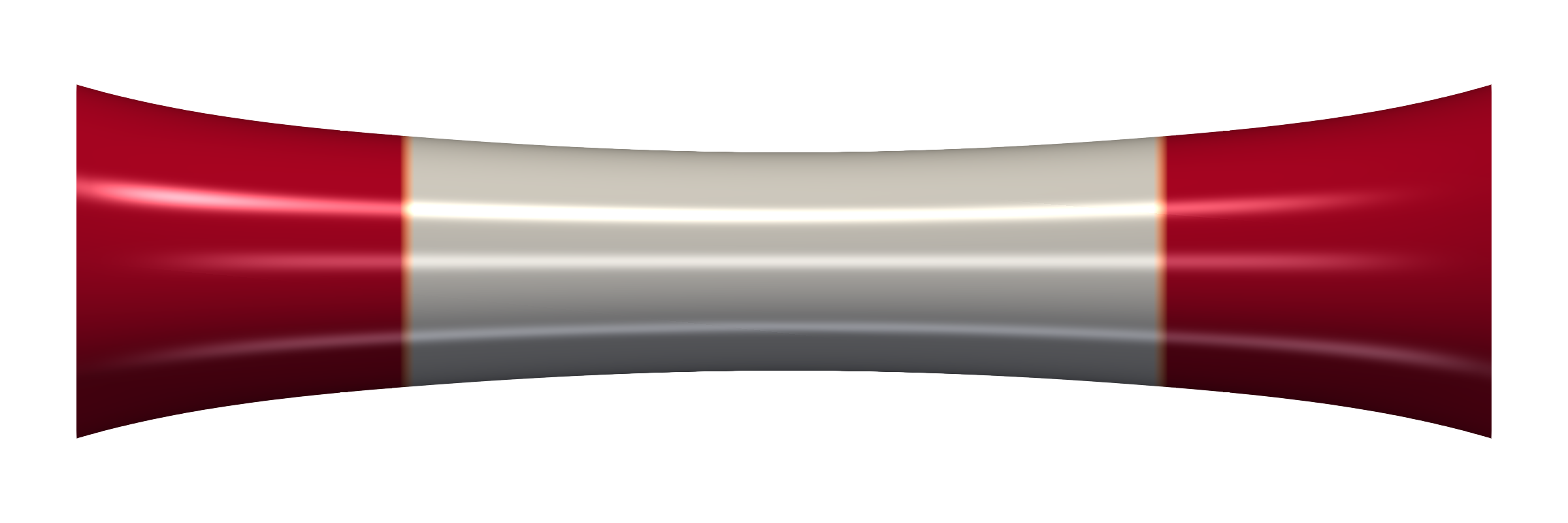}		\includegraphics[width=\linewidth]{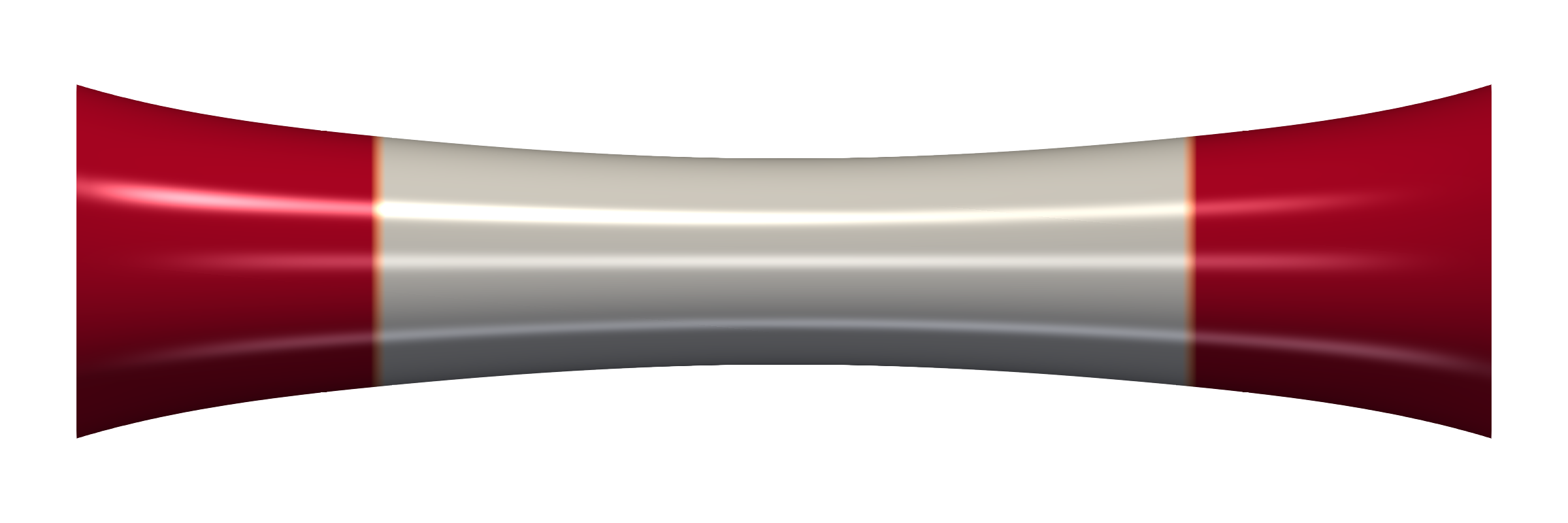}
		\caption{Side view of the deformed cylinder using the tension field theory membrane model with a rotation of $\{0.35,0.40,0.45,0.50\}\:[\pi\:\text{rad}]$ from top to bottom. The red region denotes a taut, and the gray region denotes a wrinkled region.}
		\label{fig:benchmarks_cylinder_side_tensionfield}
	\end{subfigure}
\caption{(Caption on next page)}
\end{figure}

\begin{figure}
\ContinuedFloat
	\centering
	\begin{subfigure}[t]{0.45\linewidth}
		\centering
		\includegraphics[width=\linewidth]{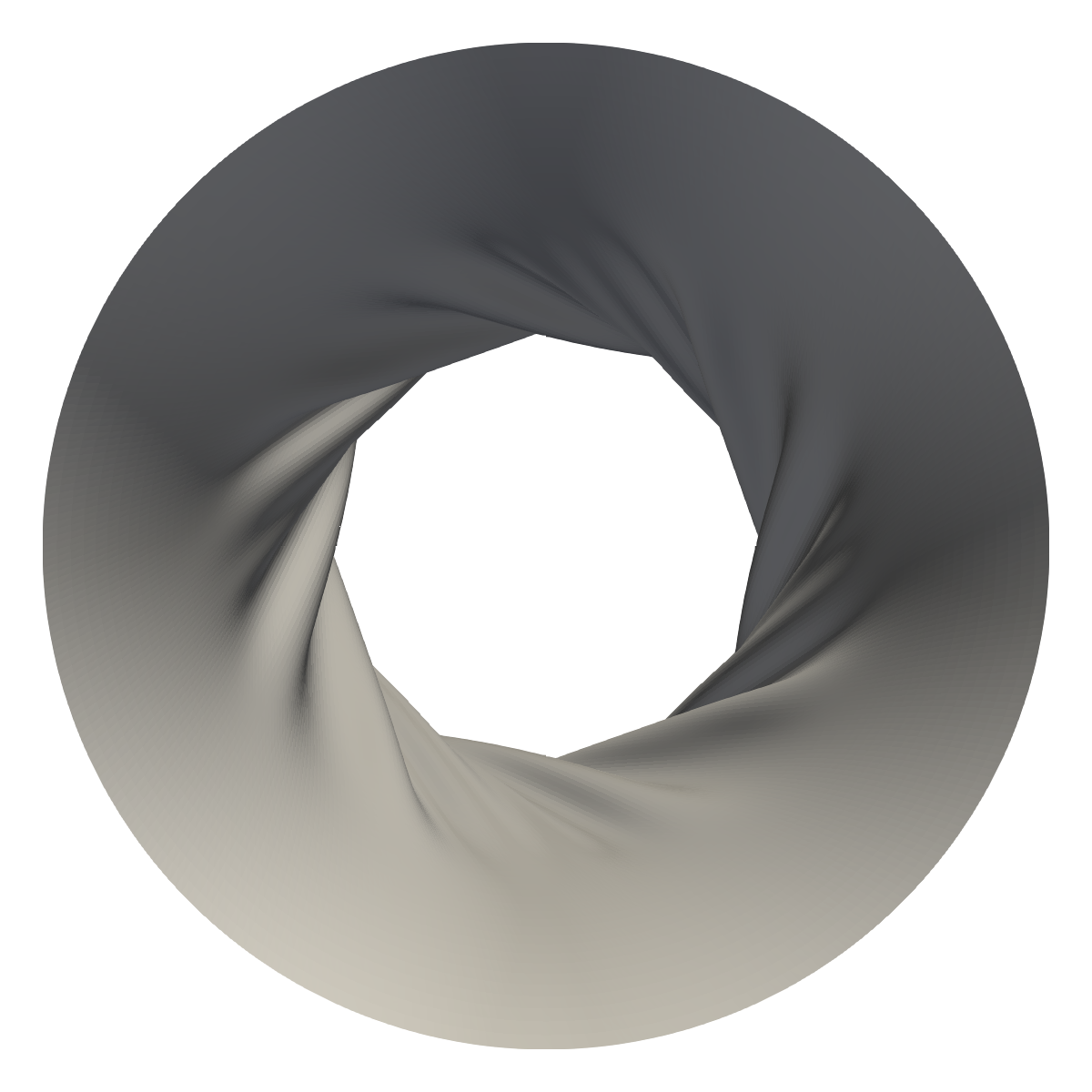}
		\caption{Top view of the deformed cylinder using the Kirchhoff--Love shell model.}
		\label{fig:benchmarks_cylinder_top}
	\end{subfigure}
	\hfill
	\begin{subfigure}[t]{0.45\linewidth}
		\centering
		\includegraphics[width=\linewidth]{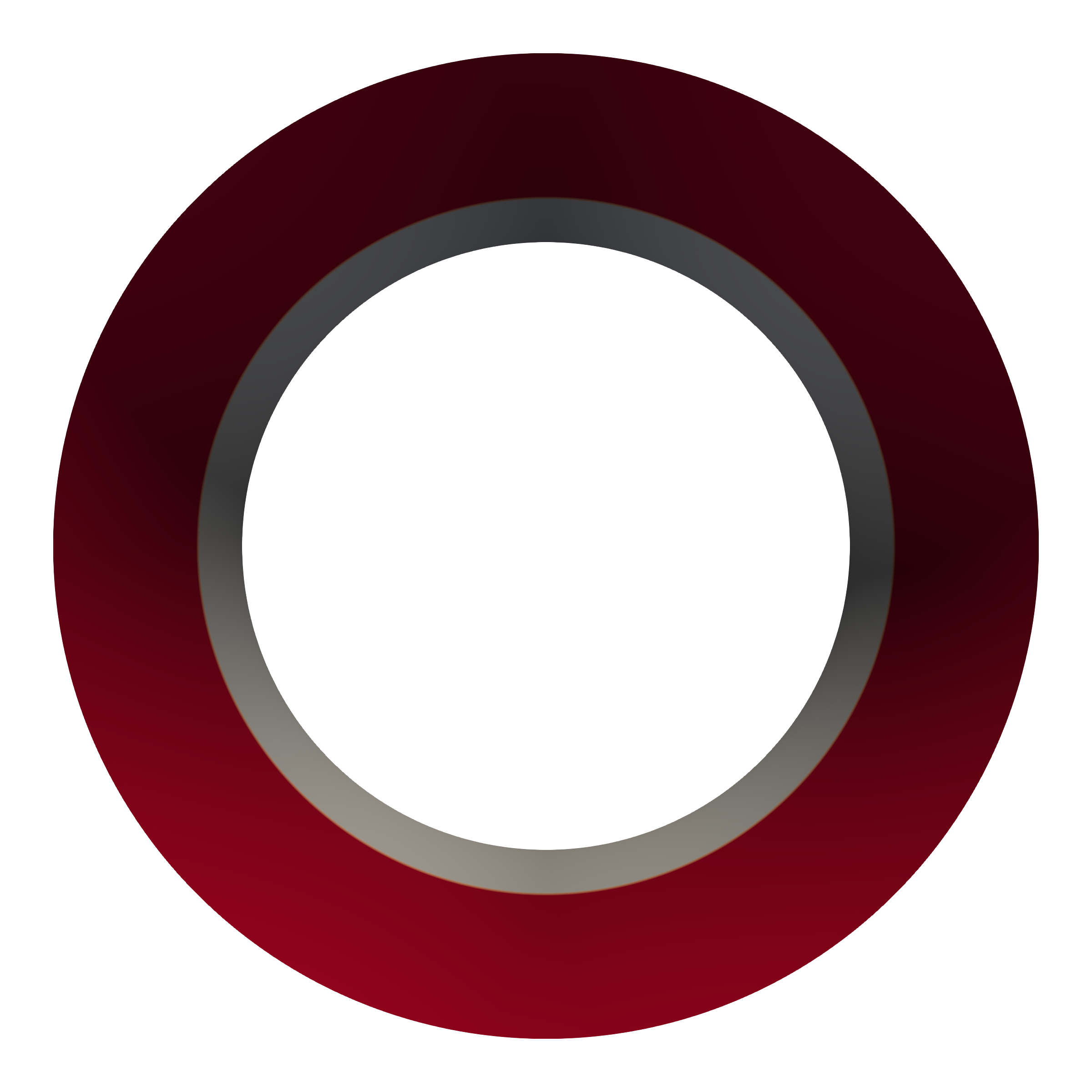}
		\caption{Top view of the deformed cylinder using the tension field theory membrane model. The red region denotes a taut, and the gray region denotes a wrinkled region.}
		\label{fig:benchmarks_cylinder_top_tensionfield}
	\end{subfigure}
	\caption{Results of the example with an cylinder with fixed bottom boundary and with a top boundary subject to a translation and a rotation, see \cref{fig:benchmarks_cylinder_setup}. A side view (\subref{fig:benchmarks_cylinder_side}) and a top view (\subref{fig:benchmarks_cylinder_top}) of the wrinkled membrane using the Kirchhoff--Love shell simulation are provided, as well as the deformed geometry from the tension field theory membrane simulation with the tension field for colouring (\subref{fig:benchmarks_cylinder_top_tensionfield}).}
	\label{fig:benchmarks_cylinder}
\end{figure}

From the results in \cref{fig:benchmarks_cylinder,fig:benchmarks_cylinder_contours}, similar observations as for the previous benchmark problem can be made. From \cref{fig:benchmarks_cylinder_contours} it can be seen that the present model predicts the mid-plane of the wrinkled membrane. In addition, \cref{fig:benchmarks_cylinder_torque} shows great similarity between the load-displacement curves obtained by using the proposed tension field theory membrane model and the reference Kirchhoff--Love shell model, even for the coarsest mesh ($16\times16$) and lowest degree ($p=2$).

\begin{figure}
	\centering
	\begin{subfigure}[t]{0.225\linewidth}
		\centering
		\begin{tikzpicture}
			\begin{axis}
				[
				height=\linewidth,
				width=\linewidth,
				scale only axis,
				hide axis,
				]
				\addplot+[style=style1,solid,opacity=0.5,no markers] table[col sep=comma,header=true,x index={0},y index={1}] {cylinder_e1_r2_line_D=0.700000_u=0.5_dir=1.csv};
				\addplot+[style=style2,solid,opacity=0.5,no markers] table[col sep=comma,header=true,x index={0},y index={1}] {cylinder_e1_r3_line_D=0.700000_u=0.5_dir=1.csv};
				\addplot+[style=style3,solid,opacity=0.5,no markers] table[col sep=comma,header=true,x index={0},y index={1}] {cylinder_e1_r4_line_D=0.700000_u=0.5_dir=1.csv};

				\addplot+[style=style1,dotted,no markers] table[col sep=comma,header=true,x index={0},y index={1}] {cylinder_e0_r2_TFT_line_D=0.700000_u=0.5_dir=1.csv};
				\addplot+[style=style2,dotted,no markers] table[col sep=comma,header=true,x index={0},y index={1}] {cylinder_e0_r3_TFT_line_D=0.700000_u=0.5_dir=1.csv};
				\addplot+[style=style3,dotted,no markers] table[col sep=comma,header=true,x index={0},y index={1}] {cylinder_e0_r4_TFT_line_D=0.700000_u=0.5_dir=1.csv};

				\addplot+[style=style1,dashed,no markers] table[col sep=comma,header=true,x index={0},y index={1}] {cylinder_e1_r2_TFT_line_D=0.700000_u=0.5_dir=1.csv};
				\addplot+[style=style2,dashed,no markers] table[col sep=comma,header=true,x index={0},y index={1}] {cylinder_e1_r3_TFT_line_D=0.700000_u=0.5_dir=1.csv};
				\addplot+[style=style3,dashed,no markers] table[col sep=comma,header=true,x index={0},y index={1}] {cylinder_e1_r4_TFT_line_D=0.700000_u=0.5_dir=1.csv};
			\end{axis}
		\end{tikzpicture}
		\caption*{$0.35\:[\pi\:\text{rad}]$}
	\end{subfigure}
	\hfill
	\begin{subfigure}[t]{0.225\linewidth}
		\centering
		\begin{tikzpicture}
			\begin{axis}
				[
				height=\linewidth,
				width=\linewidth,
				scale only axis,
				hide axis,
				]
				\addplot+[style=style1,solid,opacity=0.5,no markers] table[col sep=comma,header=true,x index={0},y index={1}] {cylinder_e1_r2_line_D=0.800000_u=0.5_dir=1.csv};
				\addplot+[style=style2,solid,opacity=0.5,no markers] table[col sep=comma,header=true,x index={0},y index={1}] {cylinder_e1_r3_line_D=0.800000_u=0.5_dir=1.csv};
				\addplot+[style=style3,solid,opacity=0.5,no markers] table[col sep=comma,header=true,x index={0},y index={1}] {cylinder_e1_r4_line_D=0.800000_u=0.5_dir=1.csv};

				\addplot+[style=style1,dotted,no markers] table[col sep=comma,header=true,x index={0},y index={1}] {cylinder_e0_r2_TFT_line_D=0.800000_u=0.5_dir=1.csv};
				\addplot+[style=style2,dotted,no markers] table[col sep=comma,header=true,x index={0},y index={1}] {cylinder_e0_r3_TFT_line_D=0.800000_u=0.5_dir=1.csv};
				\addplot+[style=style3,dotted,no markers] table[col sep=comma,header=true,x index={0},y index={1}] {cylinder_e0_r4_TFT_line_D=0.800000_u=0.5_dir=1.csv};

				\addplot+[style=style1,dashed,no markers] table[col sep=comma,header=true,x index={0},y index={1}] {cylinder_e1_r2_TFT_line_D=0.800000_u=0.5_dir=1.csv};
				\addplot+[style=style2,dashed,no markers] table[col sep=comma,header=true,x index={0},y index={1}] {cylinder_e1_r3_TFT_line_D=0.800000_u=0.5_dir=1.csv};
				\addplot+[style=style3,dashed,no markers] table[col sep=comma,header=true,x index={0},y index={1}] {cylinder_e1_r4_TFT_line_D=0.800000_u=0.5_dir=1.csv};
			\end{axis}
		\end{tikzpicture}
		\caption*{$0.40\:[\pi\:\text{rad}]$}
	\end{subfigure}
	\hfill
	\begin{subfigure}[t]{0.225\linewidth}
		\centering
		\begin{tikzpicture}
			\begin{axis}
				[
				height=\linewidth,
				width=\linewidth,
				scale only axis,
				hide axis,
				]
				\addplot+[style=style1,solid,opacity=0.5,no markers] table[col sep=comma,header=true,x index={0},y index={1}] {cylinder_e1_r2_line_D=0.900000_u=0.5_dir=1.csv};
				\addplot+[style=style2,solid,opacity=0.5,no markers] table[col sep=comma,header=true,x index={0},y index={1}] {cylinder_e1_r3_line_D=0.900000_u=0.5_dir=1.csv};
				\addplot+[style=style3,solid,opacity=0.5,no markers] table[col sep=comma,header=true,x index={0},y index={1}] {cylinder_e1_r4_line_D=0.900000_u=0.5_dir=1.csv};

				\addplot+[style=style1,dotted,no markers] table[col sep=comma,header=true,x index={0},y index={1}] {cylinder_e0_r2_TFT_line_D=0.900000_u=0.5_dir=1.csv};
				\addplot+[style=style2,dotted,no markers] table[col sep=comma,header=true,x index={0},y index={1}] {cylinder_e0_r3_TFT_line_D=0.900000_u=0.5_dir=1.csv};
				\addplot+[style=style3,dotted,no markers] table[col sep=comma,header=true,x index={0},y index={1}] {cylinder_e0_r4_TFT_line_D=0.900000_u=0.5_dir=1.csv};

				\addplot+[style=style1,dashed,no markers] table[col sep=comma,header=true,x index={0},y index={1}] {cylinder_e1_r2_TFT_line_D=0.900000_u=0.5_dir=1.csv};
				\addplot+[style=style2,dashed,no markers] table[col sep=comma,header=true,x index={0},y index={1}] {cylinder_e1_r3_TFT_line_D=0.900000_u=0.5_dir=1.csv};
				\addplot+[style=style3,dashed,no markers] table[col sep=comma,header=true,x index={0},y index={1}] {cylinder_e1_r4_TFT_line_D=0.900000_u=0.5_dir=1.csv};
			\end{axis}
		\end{tikzpicture}
		\caption*{$0.45\:[\pi\:\text{rad}]$}
	\end{subfigure}
	\hfill
	\begin{subfigure}[t]{0.225\linewidth}
		\centering
		\begin{tikzpicture}
			\begin{axis}
				[
				height=\linewidth,
				width=\linewidth,
				scale only axis,
				hide axis,
				]
				\addplot+[style=style1,solid,opacity=0.5,no markers] table[col sep=comma,header=true,x index={0},y index={1}] {cylinder_e1_r2_line_D=1.000000_u=0.5_dir=1.csv};
				\addplot+[style=style2,solid,opacity=0.5,no markers] table[col sep=comma,header=true,x index={0},y index={1}] {cylinder_e1_r3_line_D=1.000000_u=0.5_dir=1.csv};
				\addplot+[style=style3,solid,opacity=0.5,no markers] table[col sep=comma,header=true,x index={0},y index={1}] {cylinder_e1_r4_line_D=1.000000_u=0.5_dir=1.csv};

				\addplot+[style=style1,dotted,no markers] table[col sep=comma,header=true,x index={0},y index={1}] {cylinder_e0_r2_TFT_line_D=1.000000_u=0.5_dir=1.csv};
				\addplot+[style=style2,dotted,no markers] table[col sep=comma,header=true,x index={0},y index={1}] {cylinder_e0_r3_TFT_line_D=1.000000_u=0.5_dir=1.csv};
				\addplot+[style=style3,dotted,no markers] table[col sep=comma,header=true,x index={0},y index={1}] {cylinder_e0_r4_TFT_line_D=1.000000_u=0.5_dir=1.csv};

				\addplot+[style=style1,dashed,no markers] table[col sep=comma,header=true,x index={0},y index={1}] {cylinder_e1_r2_TFT_line_D=1.000000_u=0.5_dir=1.csv};
				\addplot+[style=style2,dashed,no markers] table[col sep=comma,header=true,x index={0},y index={1}] {cylinder_e1_r3_TFT_line_D=1.000000_u=0.5_dir=1.csv};
				\addplot+[style=style3,dashed,no markers] table[col sep=comma,header=true,x index={0},y index={1}] {cylinder_e1_r4_TFT_line_D=1.000000_u=0.5_dir=1.csv};
			\end{axis}
		\end{tikzpicture}
		\caption*{$0.50\:[\pi\:\text{rad}]$}
	\end{subfigure}
	\caption{Contour line for the deformed geometry in the middle of the cylinder. The solid lines represent the Kirchhoff-Love shell results, and the dashed lines represent the tension field theory membrane results. See \cref{fig:benchmarks_cylinder_contours} for the legend.}
	\label{fig:benchmarks_cylinder_contours}
\end{figure}

\begin{figure}
	\centering
	\begin{tikzpicture}
		\begin{groupplot}
			[
			height=0.2\textheight,
			width=0.38\linewidth,
			enlarge x limits = true,
			enlarge y limits = true,
			scale only axis,
			group style={
				group name = top,
				group size=2 by 1,
				xlabels at=edge bottom,
				x descriptions at=edge bottom,
				vertical sep=0.05\textheight,
				horizontal sep=0.04\linewidth,
			},
			legend pos = south east,
			]
			\nextgroupplot[xlabel = {Rotation $[\pi\:\text{rad}]$},
			ylabel = {Applied torque $M$\:[Nm]}]

			\addplot+[style=style1,solid,opacity=0.5,no markers] table[col sep=comma,header=true,x expr=\thisrowno{0}/2, y index={1}] {cylinder_e1_r2_out.csv};
			\addplot+[style=style2,solid,opacity=0.5,no markers] table[col sep=comma,header=true,x expr=\thisrowno{0}/2, y index={1}] {cylinder_e1_r3_out.csv};
			\addplot+[style=style3,solid,opacity=0.5,no markers] table[col sep=comma,header=true,x expr=\thisrowno{0}/2, y index={1}] {cylinder_e1_r4_out.csv};

			\addplot+[style=style1,no markers,dotted] table[col sep=comma,header=true,x expr=\thisrowno{0}/2, y index={1}] {cylinder_e0_r2_TFT_out.csv};
			\addplot+[style=style2,no markers,dotted] table[col sep=comma,header=true,x expr=\thisrowno{0}/2, y index={1}] {cylinder_e0_r3_TFT_out.csv};
			\addplot+[style=style3,no markers,dotted] table[col sep=comma,header=true,x expr=\thisrowno{0}/2, y index={1}] {cylinder_e0_r4_TFT_out.csv};

			\addplot+[style=style1,no markers,dashed] table[col sep=comma,header=true,x expr=\thisrowno{0}/2, y index={1}] {cylinder_e1_r2_TFT_out.csv};
			\addplot+[style=style2,no markers,dashed] table[col sep=comma,header=true,x expr=\thisrowno{0}/2, y index={1}] {cylinder_e1_r3_TFT_out.csv};
			\addplot+[style=style3,no markers,dashed] table[col sep=comma,header=true,x expr=\thisrowno{0}/2, y index={1}] {cylinder_e1_r4_TFT_out.csv};

			\draw[gray,dashed] (axis cs:0.2,2600) rectangle (axis cs: 0.5,3600);
			\node[anchor=south east,gray] at (axis cs:0.5,2600) {Inset};

			\nextgroupplot[xlabel = {Rotation $[\pi\:\text{rad}]$},
			xmin = 0.2, xmax = 0.5,
			ymin = 2600, ymax = 3600,
			enlargelimits=false,
			draw=gray,
			xlabel style={color=gray},
			x tick label style={color = gray},
			ylabel style={color=gray},
			y tick label style={color = gray},
			ylabel near ticks, yticklabel pos=right]

			\addlegendimage{style=style1,solid,no markers}\addlegendentry{$16\times16$};
			\addlegendimage{style=style2,solid,no markers}\addlegendentry{$32\times32$};
			\addlegendimage{style=style3,solid,no markers}\addlegendentry{$64\times64$};
			\addlegendimage{black,solid}\addlegendentry{KL-shell, $p=3$};
			\addlegendimage{black,dotted}\addlegendentry{TFT-membrane, $p=2$};
			\addlegendimage{black,dashed}\addlegendentry{TFT-membrane, $p=3$};

			\addplot+[style=style1,solid,opacity=0.5,no markers] table[col sep=comma,header=true,x expr=\thisrowno{0}/2, y index={1}] {cylinder_e1_r2_out.csv};
			\addplot+[style=style2,solid,opacity=0.5,no markers] table[col sep=comma,header=true,x expr=\thisrowno{0}/2, y index={1}] {cylinder_e1_r3_out.csv};
			\addplot+[style=style3,solid,opacity=0.5,no markers] table[col sep=comma,header=true,x expr=\thisrowno{0}/2, y index={1}] {cylinder_e1_r4_out.csv};

			\addplot+[style=style1,no markers,dotted] table[col sep=comma,header=true,x expr=\thisrowno{0}/2, y index={1}] {cylinder_e0_r2_TFT_out.csv};
			\addplot+[style=style2,no markers,dotted] table[col sep=comma,header=true,x expr=\thisrowno{0}/2, y index={1}] {cylinder_e0_r3_TFT_out.csv};
			\addplot+[style=style3,no markers,dotted] table[col sep=comma,header=true,x expr=\thisrowno{0}/2, y index={1}] {cylinder_e0_r4_TFT_out.csv};

			\addplot+[style=style1,no markers,dashed] table[col sep=comma,header=true,x expr=\thisrowno{0}/2, y index={1}] {cylinder_e1_r2_TFT_out.csv};
			\addplot+[style=style2,no markers,dashed] table[col sep=comma,header=true,x expr=\thisrowno{0}/2, y index={1}] {cylinder_e1_r3_TFT_out.csv};
			\addplot+[style=style3,no markers,dashed] table[col sep=comma,header=true,x expr=\thisrowno{0}/2, y index={1}] {cylinder_e1_r4_TFT_out.csv};

		\end{groupplot}
	\end{tikzpicture}
	\caption{Load-displacement curves of the torque $M$ applied on the top boundary (vertical axis) versus the rotation of the top boundary (horizontal axis) of the cylinder. The full diagram up to a rotation of $\pi/2$ radians is given in the left figure, while an inset of the final end of the curve is given in the right figure. The results are provided for different mesh sizes for the Kirchhoff--Love shell (KL-shell) and the tension field theory membrane (TFT-membrane) models fir the considered degrees.}
	\label{fig:benchmarks_cylinder_torque}
\end{figure}


%% file: Sections_Conclusions.tex
Aiming for efficient modelling of membrane wrinkling, this paper presents a modification strategy for hyperelastic membranes in an implicit way. The model extends the work of \fullcitelist{Nakashino2005}{Nakashino2005,Nakashino2020}, who presented the modification scheme for linear elastic membranes based on tension fields. The modification scheme for hyperelastic membranes presented in this paper uses the assumptions from \fullcitelist{Roddeman1987}{Roddeman1987,Roddeman1987a}, implying a modification of the deformation tensor rather than the constitutive model. As a consequence, the kinematic equation changes. Assuming a non-linear constitutive relation, the introduction of the wrinkling strain tensor adds contributions to the stress and material tensors, including dependency on the derivative of the material tensor with respect to the strain tensor. Since the latter term can be difficult to obtain due to condensation of the through-thickness strains for incompressible and compressible material models, it is computed through finite differences in the present work. Since the modification scheme is defined through derivatives of the strain energy density function, it can be used for general hyperelastic materials.\\

The present model is verified using a series of benchmarks. The isogeometric Kirchhoff--Love shell formulation is used as a reference for wrinkling computations. The proposed tension field theory-based (TFT) membrane is implemented in isogeometric finite elements as well, although the formulations do not necessarily depend on high geometric continuity. The first benchmark problem involves an uniaxial tension test, verifying the speed of convergence for a fixed tension field. The second benchmark demonstrates application on an inflatable hyperelastic membrane, showing good mesh convergence in terms of displacements compared to the shell model. Third, the TFT model applied to an annular geometry subject to a controlled displacement and rotation of the inner boundary shows a good comparison with the shell model, similar to the linear model from \fullcite{Taylor2014}. The resulting torque on the inner boundary shows great similarity between the TFT model and the shell model, and the TFT model shows consistent results through mesh refinements. Similar to the annular membrane, a new benchmark problem is provided using a cylinder that is elongated and rotated along its central axis. In this case, the TFT membrane model provides great similarity with the shell model both in terms of the computed displacement field and the applied torque on the boundary. In this case, the proposed TFT model also shows consistency over mesh refinements. It should be noted that the use of the Dynamic Relaxation method in combination with the Newton--Raphson method is necessary in some situations due to evolving tension fields during the iterations, making the Newton--Raphson methods diverge with a poor initial guess. Concluding, the benchmarks demonstrate the validity and applicability of the proposed model, in addition they show that the proposed model provides accurate prediction of global structural response even on very coarse meshes.

As for the linear elastic model from \fullcite{Nakashino2005}, the present model relies on a tension field evaluated given a certain deformation. If the tension field is subject to large changes during iterations, the convergence behaviour of the model deteriorates since the changing tension field is not included in the variation of the wrinkling stress. In future work, it is recommended to incorporate the tension field into the derivative of the wrinkling stress. Furthermore, an equivalent for the Dynamic-Relaxation method combined with Newton--Raphson methods for quasi-static simulations could be investigated as well, starting with the explicit arc-length method from \fullcite{Lee2011}.